\documentclass[]{amsart}

\usepackage{amssymb}
\usepackage{url}
\usepackage{amscd}
\usepackage{epsfig}
\usepackage{subfigure}

\theoremstyle{plain}
\newtheorem{theorem}{Theorem}

\newtheorem{corollary}[theorem]{Corollary}
\newtheorem{lemma}[theorem]{Lemma}
\newtheorem{proposition}[theorem]{Proposition}
\newtheorem{definition}{Definition}

\theoremstyle{remark}
\newtheorem{remark}{Remark}
\newtheorem{example}{Example}

\def\R{{\mathbb R}}
\def\C{{\mathbb C}}
\def\H{{\mathbb H}}
\def\E{{\mathbb E}}
\def\S{{\mathbb S}}
\def\L{{\mathbb L}}
\def\N{{\mathbb N}}
\def\bP{{\mathbb P}}
\def\bQ{{\mathbb Q}}
\def\v{{\bf v}}
\def\w{{\bf w}}
\def\u{{\bf u}}
\def\1{{\bf 1}}
\def\i{{\bf i}}
\def\j{{\bf j}}
\def\k{{\bf k}}
\def\d{\partial}
\def\A{{\mathcal A}}

\def\Q{{\mathcal Q}}
\def\G{{\mathcal G}}
\def\F{{\mathcal F}}
\def\K{{\mathcal K}}
\def\M{{\mathcal M}}
\def\U{{\mathcal U}}
\def\V{{\mathcal V}}
\def\s{{\mathfrak s\mathfrak l}}

\def\tilde{\widetilde}
\def\dis{\displaystyle}
\def\tr{\mathop{\rm tr}\nolimits}
\def\cmc{{\mathsf c\mathsf m\mathsf c}}

\input{epsf.tex}

\begin{document}

\date{Received July 14, 2005}

\title[Timelike $\cmc$ $\pm 1$ Surfaces in AdS 3-Space
  $\H^3_1(-1)$]{Timelike Surfaces of Constant Mean Curvature $\pm 1$
  in Anti-de Sitter 3-Space $\H^3_1(-1)$}

\author{Sungwook Lee}

\address{Department of Mathematics\\
University of Southern Mississippi\\
Hattiesburg, MS 39406-5045, U.S.A.}

\email{sunglee@usm.edu}


\dedicatory{This paper is dedicated to my sensei professor
    Kinetsu Abe.}

\begin{abstract}
It is shown that timelike surfaces of constant mean curvature $\pm 1$ in
anti-de Sitter $3$-space ${\mathbb H}^3_1(-1)$ can
be constructed from a pair of Lorentz holomorphic and Lorentz
antiholomorphic null curves in ${\mathbb P}{\rm SL}_2{\mathbb R}$ via
Bryant type representation formulae. These Bryant type representation
formulae are used to investigate an explicit one-to-one correspondence, the
so-called \emph{Lawson-Guichard correspondence}, between timelike
surfaces of constant mean curvature $\pm 1$ in ${\mathbb H}^3_1(-1)$ and timelike
minimal 
surfaces in Minkowski $3$-space ${\mathbb E}^3_1$. The hyperbolic Gau{\ss} map
of timelike surfaces in ${\mathbb H}^3_1(-1)$, which is a close analogue of the
classical Gau{\ss} map is considered. It is discussed that the
hyperbolic Gau{\ss} map plays an important role in the study of
timelike surfaces of constant mean curvature $\pm 1$ in
${\mathbb H}^3_1(-1)$. In particular, the relationship between the Lorentz
holomorphicity of the hyperbolic Gau{\ss} map and timelike surface of
constant mean curvature $\pm 1$ in ${\mathbb H}^3_1(-1)$ is studied.
\end{abstract}

\maketitle

\section{Introduction}
\label{sec:intro}
 It is known that surfaces of constant mean curvature $\pm 1$
surfaces in hyperbolic $3$-space $\H^3(-1)$ can be constructed from
holomorphic null curves in $\bP{\rm SL}_2\C={\rm
  SL}_2\C/\{\pm{\rm id}\}$ (\cite{Br}, \cite{U-Y}), while minimal
  surfaces in Euclidean $3$-space $\E^3$ can be constructed from
  holomorphic null curves in $\C^3$ via well-known Weierstra{\ss}-Enneper
  representation formula. It is also known that spacelike surfaces of
  constant mean curvature $\pm 1$ in de-Sitter $3$-space $\S^3_1(1)$
  can be constructed from holomorphic null curves in $\bP{\rm SL}_2\C$
  (\cite{A-A2}, \cite{Lee}), while spacelike
  maximal surfaces in Minkowski $3$-space $\E^3_1$ can be constructed
  from holomorphic null curves in $\C^3$ via an analogue of
  Weierstra{\ss}-Enneper representation formula (\cite{Mc},
  \cite{Kob}). These are all related by the \emph{Lawson-Guichard
    correspondence} between minimal surfaces in $\E^3$ and surfaces of
  constant mean curvature $\pm 1$ in $\H^3(-1)$ (\cite{Law}) and the
  one between spacelike maximal surfaces in $\E^3_1$ and spacelike
  surfaces of constant mean curvature $\pm 1$ (\cite{Palm}). Note that
  the correspondents (they are usually called the \emph{cousins}) in
  different space forms satisfy the same Gau{\ss} and Mainardi-Codazzi equations.

It is interesting to see that there exists a Lawson-Guichard
correspondence between timelike minimal surfaces in $\E^3_1$ and
timelike surfaces of constant mean curvature $\pm 1$ in anti-de
Sitter $3$-space $\H^3_1(-1)$. See sections
\ref{sec:intsystem}, \ref{sec:corresp}, and \ref{sec:lawson}
(appendix I) for details. In \cite{I-T}, J.~Inoguchi and M.~Toda
show that timelike minimal surfaces can be constructed from a pair
of Lorentz holomorphic and Lorentz antiholomorphic null curves in
$\R^3$ via normalized Weierstra{\ss} formula \eqref{eq:nwf}. Hence,
one might expect a similar construction of timelike surfaces of
constant mean curvature $\pm 1$ in $\H^3_1(-1)$ in terms of Lorentz
holomorphic and Lorentz antiholomorphic null curves. In this paper,
we prove that a pair of Lorentz holomorphic and Lorentz
antiholomorphic null curves in $\bP{\rm SL}_2\R$ gives rise to a
timelike surface of constant mean curvature $\pm 1$ in $\H^3_1(-1)$. Furthermore,
every timelike surface of constant mean curvature $\pm 1$ in
$\H^3_1(-1)$ can be constructed from a pair of Lorentz holomorphic
and Lorentz antiholomorphic null curves in $\bP{\rm SL}_2\R$.

An analogue of the hyperbolic Gau{\ss} map\footnote{The hyperbolic
  Gau{\ss} map was introduced by C.~Epstein in \cite{Ep} and used by
  R.~.L.~Bryant to study $\cmc$ $1$ surfaces in $\H^3(-1)$ of surfaces in
hyperbolic $3$-space $\H^3(-1)$ in \cite{Br}.} can be defined for
timelike surfaces of constant mean curvature in
  $\H^3_1(-1)$ and plays an important role in studying timelike
  surfaces of constant mean curvature $\pm 1$ in $\H^3_1(-1)$. It is
  shown in section \ref{sec:hol} that
\begin{enumerate}
\item The hyperbolic Gau{\ss} map of a Lorentz surface
$\varphi: M\longrightarrow\H^3_1(-1)$ is Lorentz holomorphic if and
only if $\varphi$ satisfies $H=1$ and $Q=0$.\\
\item The hyperbolic Gau{\ss} map of a Lorentz surface
$\varphi: M\longrightarrow\H^3_1(-1)$ is Lorentz antiholomorphic if
and only if $\varphi$ satisfies $H=1$ and $R=0$.
\end{enumerate}
Here, $Q$ and $R$ are coefficients of quadratic differentials, the
so-called \emph{Hopf
  pairs}. They are defined in the following section.

In \cite{Akutagawa} and \cite{Ramanathan}, K.~Akutagawa and
J.~Ramanathan proved independently that

\noindent{\bf Theorem:} \emph{Let $M$ be a complete spacelike
surface in de Sitter $3$-space $\S^3_1(1)$ with constant mean
curvature $H=\pm 1$. Then $M$ is a totally umbilic flat surface.
Moreover, $M$ is a parabolic type surface of revolution.}
\\
This theorem tells us that de Sitter $3$-space $\S^3_1(1)$ admits
horosphere type spacelike surfaces. It is also interesting to see that
anti-de Sitter $3$-space $\H^3_1(-1)$ admits horosphere type
timelike surfaces. See section \ref{sec:hypgauss} for details.
\section{Timelike Surfaces in Anti-de Sitter 3-Space $\H^3_1(-1)$}
\label{sec:timelike}
Let $\E^4_2$ be the semi-Euclidean $4$-space with natural coordinates\\
$(x_0,x_1,x_2,x_3)$ and the semi-Riemannian metric
$\langle\cdot,\cdot\rangle$ of
signature $(-,-,+,+)$ given by the quadratic form
$-(dx_0)^2-(dx_1)^2+(dx_2)^2+(dx_3)^2$.

The anti-de Sitter (abbreviated: AdS) $3$-space $\H^3_1(-1)$ is a Lorentzian
$3$-manifold of sectional curvature $-1$ that can be
realized as the hyperquadric in $\E^4_2$:
$$\H^3_1(-1):=\{(x_0,x_1,x_2,x_3)\in\E^4_2:
-(x_0)^2-(x_1)^2+(x_2)^2+(x_3)^2=-1\}.$$

Let $M$ be a connected orientable $2$-manifold and $\varphi:
M\longrightarrow\H^3_1(-1)$ an immersion. The immersion $\varphi$ is said
to be \emph{ timelike} if the induced metric $I$ on $M$ is
Lorentzian. The induced Lorentzian metric $I$ determines a Lorentzian
conformal structure $\mathcal{C}_I$ on $M$.

Let $(x,y)$ be a Lorentz isothermal coordinate system with respect to
the conformal structure $\mathcal{C}_I$. Then the first fundamental
form
$$I=\langle d\varphi,d\varphi\rangle$$
is given by the matrix
$I=e^\omega\begin{pmatrix}
-1 & 0\\
0 & 1
\end{pmatrix}$. The first fundamental form is also written in
terms of $(x,y)$ as $I=e^\omega\{-(dx)^2+(dy)^2\}$.
Let $u:=x+y$ and $v:=-x+y$. Then $(u,v)$ defines a \emph{null coordinate
system} with respect to the conformal structure $\mathcal{C}_I$. The
first fundamental form $I$ is written in terms of $(u,v)$ as
$$I=e^\omega dudv.$$
In terms of null coordinates $u$ and $v$, the differential operators
$\frac{\d}{\d u}$ and $\frac{\d}{\d v}$ are computed to be
$$\frac{\d}{\d u}=\frac{1}{2}\left(\frac{\d}{\d x}+\frac{\d}{\d
  y}\right),\
\frac{\d}{\d v}=\frac{1}{2}\left(-\frac{\d}{\d x}+\frac{\d}{\d
  y}\right).$$
The
conformality condition is equivalent to
$$\langle\varphi_u,\varphi_u\rangle=\langle\varphi_v,\varphi_v\rangle=0,\
\langle\varphi_u,\varphi_v\rangle=\frac{1}{2}e^\omega.$$
Let $N$ be a unit normal vector field of $M$. Then
$$\langle N,N\rangle=1,\
\langle\varphi,N\rangle=\langle\varphi_u,N\rangle=\langle\varphi_v,N\rangle=0.$$
The mean curvature $H$ is given by $H=2e^{-\omega}<\varphi_{uv},N>$.
Let $Q:=<\varphi_{uu},N>$ and $R:=<\varphi_{vv},N>$. Then the quadratic
differentials $Q^\sharp:=Qdu\otimes du$ and $R^\sharp:=Rdv\otimes dv$
are called \emph{ Hopf pairs}\footnote{In \cite{D-I-T},
  \cite{In}, \cite{I-T}, the quadratic differentials $Q\sharp$ and
  $R\sharp$ are defined as Hopf differentials.} of $M$. The quadratic
differential
$$\Q:=Qdu^2+Rdv^2=Q^\sharp+R^\sharp$$
is called \emph{Hopf differential}\footnote{The definition of Hopf
  differential $\Q$ was suggested to the author by J. Inoguchi
  \cite{In2}.}. This differential is
globally defined on the Lorentz surface $(M,\mathcal{C}_I)$.
The second fundamental form $I\!I$ of $M$ derived from $N$ is defined by
$$I\!I:=-\langle d\varphi,dN\rangle$$
and it is given by the matrix
$$I\!I=\begin{pmatrix}
Q+R-He^\omega & Q-R\\
Q-R & Q+R+He^\omega
\end{pmatrix}$$
with respect to Lorentz isothermal coordinate system $(x,y)$.
The second fundamental form is related to Hopf differential $\Q$ by
\begin{equation}
\label{eq:sff}
I\!I=\Q+HI.
\end{equation}
The \emph{shape operator} ${\mathcal S}$ of $M$ derived from $N$ is
$S:=-dN$. The shape operator ${\mathcal S}$ is related to $I\!I$ by
$$I\!I(X,Y)=\langle{\mathcal S}X,Y\rangle$$
for all vector fields $X,Y$ on $M$. The shape operator ${\mathcal S}$ is also
represented by the matrix $I\!I\cdot I^{-1}$. The mean curvature $H$ of
$M$ is
$$H=\frac{1}{2}\tr{\mathcal S}=\frac{1}{2}\tr(I\!I\cdot
I^{-1})$$
and the Gau{\ss}ian curvature\footnote{This can be easily computed
  from the Gau{\ss} equation \eqref{eqn:gausseq}} $K$ of $M$ is
$$K:=-1+\det{\mathcal S}=-1+\det(I\!I\cdot I^{-1}).$$
The eigenvalues of $S$, i.e., the solutions to the characteristic equation
$$\det({\mathcal S}-\lambda{\mathcal I})=0,\ {\mathcal I}=\mbox{\rm
  identity of}\ TM$$ are called the principal curvatures. Since the
  metric $I$ is indefinite, both principal curvatures may be nonreal complex
numbers. The mean curvature $H$ is the mean of the two principal
curvatures and the Gau{\ss}ian curvature
$K$ is the product of the two principal curvatures minus one.

A point $p\in M$ is said to be an umbilic point if $I\!I$ is
proportional to $I$ at $p$. Equivalently, $p$ is an umbilic point if
and only if the two principal curvatures at $p$ are the same real
number and the corresponding eigenspace is $2$-dimensional. A timelike
surface is said to be a totally umbilic if all the points are
umbilical. The formula \eqref{eq:sff} implies that $p\in M$ is an
umbilic point if and only if $\Q(p)=0$, i.e., $p\in M$ is a common
zero of Hopf pairs $Q$ and $R$.

The Gau{\ss} equation which describes a relationship between $K$, $H$,
$Q$ and $R$ takes the following form:
\begin{equation}
\label{eq:gauss}
H^2-K-1=4e^{-2\omega}QR.
\end{equation}
Note that the condition $QR=0$ does not imply the condition $Q=R=0$
(See \cite{Mi}).

Let $M$ be a simply-connected open and orientable $2$-manifold and
$\varphi: M\longrightarrow\H^3_1(-1)$ a timelike conformal immersion with unit
normal vector field $N$. Then we can define an orthonormal frame field
$\F$ along $\varphi$ by
\begin{equation}
\label{eq:frame1}
\F=(\varphi,e^{-\frac{\omega}{2}}\varphi_x,e^{-\frac{\omega}{2}}\varphi_y,N):
M\longrightarrow{\rm O}^{++}(2,2),
\end{equation}
where ${\rm O}^{++}(2,2)$ denotes the identity component of the
Lorentz group
$${\rm O}(2,2)=\{\A\in {\rm GL}_4\R: <\A{\bf u},\A{\bf v}>=<{\bf
  u},{\bf v}>,\ {\bf u},{\bf v}\in\E^4_2\}.$$
In terms of null coordinates $(u,v)$, $\F$ is defined by
\begin{equation}
\label{eq:frame2}
\F=(\varphi,e^{-\frac{\omega}{2}}(\varphi_u-\varphi_v),e^{-\frac{\omega}{2}}
(\varphi_u+\varphi_v),N): M\longrightarrow{\rm O}^{++}(2,2).
\end{equation}
The semi-Euclidean $4$-space $\E^4_2$ is identified with the linear
space ${\rm M}_2\R$ of all $2\times 2$ real matrices via the correspondence
$${\bf u}=(x_0,x_1,x_2,x_3)\longleftrightarrow\begin{pmatrix}
x_0+x_3  & x_1+x_2\\
-x_1+x_2 & x_0-x_3
\end{pmatrix}.$$
The scalar product of $\E^4_2$ corresponds to the scalar product
\begin{equation}
\label{eq:scalproduct}
<{\bf u},{\bf v}>=\frac{1}{2}\{\tr({\bf u}{\bf v})-\tr({\bf
  u})\tr({\bf v})\},\ {\bf u},{\bf v}\in{\rm M}_2\R.
\end{equation}
Note that $<{\bf u},{\bf u}>=-\det{\bf u}$. The standard basis
$e_0,\ e_1,\ e_2,\ e_3$ for $\E^4_2$ is identified with the
matrices
$$\1=\begin{pmatrix}
1 & 0\\
0 & 1
\end{pmatrix},\ \i=\begin{pmatrix}
0 & 1\\
-1 & 0
\end{pmatrix},\ \j^{'}=\begin{pmatrix}
0 & 1\\
1 & 0
\end{pmatrix},\ \k^{'}=\begin{pmatrix}
1 & 0\\
0 & -1
\end{pmatrix},$$
i.e.,
\begin{equation}
\label{eq:matrix}
(x_0,x_1,x_2,x_3)\longleftrightarrow
x_0\1+x_1\i+x_2\j^{'}+x_3\k^{'}=\begin{pmatrix}
x_0+x_3 & x_1+x_2\\
-x_1+x_2 & x_0-x_3
\end{pmatrix}.
\end{equation}
Note that the $2\times 2$ matrices
$x_0\1+x_1\i+x_2\j^{'}+x_3\k^{'}$ form the algebra ${\mathbb H}^{'}$ of
\emph{split-quaternions}. (For more details, see, for example, \cite{I-T}.)
Under the identification \eqref{eq:matrix}, the group $G$ of timelike
unit vectors corresponds to a special linear group
$${\rm SL}_2\R=\left\{\begin{pmatrix}
a & b\\
c & d
\end{pmatrix}\in{\rm M}_2\R: ad-bc=1\right\}.$$
The metric of $G$ induced by the scalar product \eqref{eq:scalproduct}
is a bi-invariant Lorentz metric of constant curvature $-1$. Hence,
$G$ is identified with $\H^3_1(-1)$.
\section{Cartan's Formalism}
\label{sec:cartan}
Let $\{e_\alpha: \alpha=0,1,2,3\}$ be a frame field of $\E^4_2$, i.e.,
$\{e_\alpha(p): \alpha=0,1,2,3\}$ is a basis for the tangent space
${\rm T}_p\E^4_2$ at each $p\in\E^4_2$. Denote by
$\langle\cdot,\cdot\rangle_p$ the scalar product on the tangent space
${\rm T}_p\E^4_2,\ p\in\E^4_2$. Then
$$<e_\alpha,e_\beta>=\left\{\begin{array}{ccc}
-1 & {\rm if} & \alpha=\beta=0\ {\rm or}\ 1,\\
0  & {\rm if} & \alpha\ne\beta,\\
1  & {\rm if} & \alpha=\beta=2\ {\rm or}\ 3.
\end{array}\right.$$
There exist unique connection $1$-forms $\{\omega_\alpha^\beta:
\alpha, \beta=0,1,2,3\}$ such that
\begin{equation}
\label{eq:conneq}
de_\alpha=\omega_\alpha^\beta e_\beta.
\end{equation}
We use the index range $1\leq i,j,k\leq 3$ and denote by $\omega^i$
the connection form $\omega_0^i$. Then the equation \eqref{eq:conneq}
can be written
\begin{align}
\label{eq:conneq2}
de_0&=\omega^ie_i,\\
de_1&=\omega^0_1e_0+\omega^2_1e_2+\omega^3_1e_3,\\
de_i&=\omega^0_ie_0+\omega^j_ie_j,\ i=2,3.
\end{align}
The connection $1$-forms $\{\omega_\alpha^\beta: \alpha,
\beta=0,1,2,3\}$ satisfy:
\begin{eqnarray*}
\omega^1_0=-\omega^0_1,\ \omega_\alpha^\alpha=0,\ \alpha=0,1,2,3,\\
\omega^i_1=\omega^1_i,\ \omega^i_0=\omega^0_i,\ i=2,3,\\
\omega^j_i=-\omega^i_j,\ i,j=2,3.
\end{eqnarray*}
Differentiating the equation \eqref{eq:conneq2} we get the \emph{first
  structure equation}:
\begin{equation}
\label{eq:steqn1}
d\omega^i=-\omega^i_j\wedge\omega^j.
\end{equation}
Differentiating this first structure equation \eqref{eq:steqn1} we
get the \emph{second structure equation}:
\begin{equation}
\label{eq:steqn2}
d\omega^i_j=-\omega^i_k\wedge\omega^k_j.
\end{equation}
For the frame field $\mathcal F$ of timelike immersion $\varphi:
M\longrightarrow\H^3_1(-1)$, we have
\begin{align}
d\omega^1&=\omega^2\wedge\omega^1_2\ (\mbox{\rm The First Structure
  Equations})\\
d\omega^2&=\omega^1\wedge\omega^1_2\\
\label{eq:symeqn}
0&=\omega^1\wedge\omega^1_3-\omega^2\wedge\omega^2_3\ (\mbox{\rm
  Symmetry Equation})\\
\label{eq:gausseq}
d\omega^1_2&=-\omega^1\wedge\omega^2+\omega^1_3\wedge\omega^2_3\
  (\mbox{\rm Gau{\ss} Equation})\\
d\omega^1_3&=-\omega^1_2\wedge\omega^2_3\ (\mbox{\rm Mainardi-Codazzi Equations})\\
d\omega^2_3&=-\omega^1_2\wedge\omega^1_3.
\end{align}
\begin{proposition}
\label{prop:curv} Let $\varphi: M\longrightarrow\H^3_1(-1)$ be a
timelike immersion. If
$\{e_1=e^{-\frac{\omega}{2}}\varphi_x,e_2=e^{-\frac{\omega}{2}}\varphi_y,e_3=N
\}$ forms an adapted frame field along $\varphi$, then the
Gau{\ss}ian curvature $K$ and mean curvature $H$ of $\varphi$
satisfy the following equations:
\begin{eqnarray}
\omega^1_3\wedge\omega^2_3&=&(K+1)\omega^1\wedge\omega^2,\\
\label{eq:meanc}
\omega^1_3\wedge\omega^2+\omega^1\wedge\omega^2_3&=&-2H\omega^1\wedge\omega^2.
\end{eqnarray}
\end{proposition}
\begin{proof}
From the symmetry equation \eqref{eq:symeqn}, we see that there exist
smooth functions $h_{ij},\ i,j=1,2$ such that
$$\begin{pmatrix}
\omega^1_3\\
\omega^2_3
\end{pmatrix}=\begin{pmatrix}
-h_{11} & h_{12}\\
h_{21} & -h_{22}
\end{pmatrix}\begin{pmatrix}
\omega^1\\
\omega^2
\end{pmatrix}\ {\rm and}\ h_{12}=-h_{21}.$$
Note that $\omega^1$ and $\omega^2$ are the dual $1$-forms of $e_1$
and $e_2$, resp., and so
\begin{align*}
\omega^1_3\wedge\omega^2_3&=(h_{11}h_{22}+h_{12}^2)\omega^1\wedge\omega^2\\
                          &=(K+1)\omega^1\wedge\omega^2,
                          \end{align*}
                          where $K$ is the Gau{\ss}ian curvature of $\varphi$.
Thus, the Gau{\ss}ian equation \eqref{eq:gausseq} can be written as
$d\omega^1_2=K\omega^1\wedge\omega^2$. The mean curvature $H$ of
$\varphi$ is $\dis\frac{h_{11}+h_{22}}{2}$. Hence,
\begin{eqnarray*}
\omega^1_3\wedge\omega^2+\omega^1\wedge\omega^2_3&=&-(h_{11}+h_{22})\omega^1
\wedge\omega^2\\
                                                 &=&-2H\omega^1\wedge\omega^2.
\end{eqnarray*}
\end{proof}
\section{Lie Group Actions $\mu$ and $\nu$ on $\E^4_2$}
\label{sec:action}
The Lie group ${\rm SL}_2\R\times{\rm SL}_2\R$ acts isometrically on
$\E^4_2$ via the group action:
$$\mu: ({\rm SL}_2\R\times{\rm
  SL}_2\R)\times\E^4_2\longrightarrow\E^4_2;\ \mu(g_1,g_2){\bf
  u}=g_1{\bf u}g_2^t.$$
This action is transitive on $\H^3_1(-1)$. The isotropy subgroup of ${\rm
  SL}_2\R\times{\rm SL}_2\R$ at $\1$ is
$\K=\{(g,(g^{-1})^t): g\in{\rm SL}_2\R)\}$ and $\H^3_1(-1)$ is
represented as the Lorentzian symmetric space ${\rm SL}_2\R\times{\rm
  SL}_2\R/\K$. The natural projection $\pi_\mu: {\rm SL}_2\R\times{\rm
  SL}_2\R\longrightarrow\H^3_1(-1)$ is given explicitly by
  $\pi_\mu(g_1,g_2)=g_1g_2^t$.

The Lie group ${\rm SL}_2\R\times{\rm SL}_2\R$ also acts isometrically
on $\E^4_2$ via the diagonal action:
$$\nu: ({\rm SL}_2\R\times{\rm
  SL}_2\R)\times\E^4_2\longrightarrow\E^4_2;\ \mu(g_1,g_2){\bf
  u}=g_1{\bf u}g_2^{-1}.$$
This action is also transitive on $\H^3_1(-1)$. The isotropy subgroup of
  ${\rm SL}_2\R\times{\rm SL}_2\R$ at $\1$ is the diagonal subgroup
  $\Delta$ of ${\rm SL}_2\R\times{\rm SL}_2\R$, that is,
$\Delta=\{(g,g): g\in{\rm SL}_2\R\}$ and $\H^3_1(-1)$ is also
  represented by ${\rm SL}_2\R\times{\rm SL}_2\R/\Delta$ as a
  Lorentzian symmetric space. The natural projection $\pi_\nu: {\rm
  SL}_2\R\times{\rm SL}_2\R\longrightarrow\H^3_1(-1)$ is given
  explicitly by
$$\pi_\nu(g_1,g_2)=g_1g_2^{-1},\ (g_1,g_2)\in{\rm SL}_2\R\times{\rm SL}_2\R.$$
Moreover, ${\rm SL}_2\R$ acts isometrically on $\E^3_1$ via the Ad-action:
$${\rm Ad}: {\rm SL}_2\R\times\E^3_1\longrightarrow\E^3_1;\ {\rm
  Ad}(g){\bf u}=g{\bf u}g^{-1},\ g\in{\rm SL}_2\R,\ {\bf
  u}\in\E^3_1.$$
The actions $\mu$ and $\nu$ both induces a double covering ${\rm
  SL}_2\R\times{\rm SL}_2\R\longrightarrow{\rm O}^{++}(2,2)$ of the
  Lorentz group ${\rm O}^{++}(2,2)$.
\begin{remark}
In \cite{Ho}, J.~Q.~Hong used the group action $\mu$ to study a Bryant
type representation formula for timelike $\cmc$ $1$ surfaces in $\H^3_1(-1)$. In
\cite{A-A}, R.~Aiyama and K.~Akutagawa also used the action $\mu$ to
study Kenmotsu-Bryant type representation formula for spacelike $\cmc$
surfaces in $\H^3_1(-1)$. In this paper, we use both actions $\mu$ and
$\nu$.
\end{remark}
The frame field $\{e_\alpha: \alpha=0,1,2,3\}$ can be parametrized by the
Lie group ${\rm SL}_2\R\times{\rm SL}_2\R$ via the Lie group action
$\mu$: for each $g=(g_1,g_2)\in{\rm SL}_2\R\times{\rm SL}_2\R$,
\begin{align*}
e_0(g)&:=\mu(g)(\1)=g_1\1 g_2^t,\\
e_1(g)&:=\mu(g)(\i)=g_1\i g_2^t,\\
e_2(g)&:=\mu(g)(\j^{'})=g_1\j^{'}g_2^t,\\
e_3(g)&:=\mu(g)(\k^{'})=g_1\k^{'}g_2^t.
\end{align*}
The frame field $\{e_\alpha: \alpha=0,1,2,3\}$ can also be parametrized by the
Lie group ${\rm SL}_2\R\times{\rm SL}_2\R$ via the Lie group action
$\nu$: for each $g=(g_1,g_2)\in{\rm SL}_2\R\times{\rm SL}_2\R$,
\begin{align*}
e_0(g)&:=\nu(g)(\1)=g_1\1 g_2^{-1},\\
e_1(g)&:=\nu(g)(\i)=g_1\i g_2^{-1},\\
e_2(g)&:=\nu(g)(\j^{'})=g_1\j^{'}g_2^{-1},\\
e_3(g)&:=\nu(g)(\k^{'})=g_1\k^{'}g_2^{-1}.
\end{align*}
We need the following two equations in order to do some differential
geometric computations in Sections \ref{sec:bryant} and
\ref{sec:bryant2}.
\begin{lemma}
\label{lem:m-c}
If the frame field $\{e_\alpha: \alpha=0,1,2,3\}$ is parametrized by
${\rm SL}_2\R\times{\rm SL}_2\R$ via the action $\mu$, then the pull
back $g^{-1}dg$ of Maurer-Cartan form $\Omega=(\omega^\beta_\alpha)$
can be written as the following equation in the Lie algebra
$\s_2\R\oplus\s_2\R$:
$$g^{-1}dg=g_1^{-1}dg_1\oplus g_2^{-1}dg_2,$$
where
\begin{equation}
\label{eq:m-c1}
g_1^{-1}dg_1=\frac{1}{2}\begin{pmatrix}
\omega^3+\omega^1_2 & \omega^1+\omega^2-\omega^1_3-\omega^2_3\\
-\omega^1+\omega^2-\omega^1_3+\omega^2_3 & -\omega^3-\omega^1_2
\end{pmatrix}
\end{equation}
and
\begin{equation}
\label{eq:m-c2}
g_2^{-1}dg_2=\frac{1}{2}\begin{pmatrix}
\omega^3-\omega^1_2 & -\omega^1+\omega^2+\omega^1_3-\omega^2_3\\
\omega^1+\omega^2+\omega^1_3+\omega^2_3 & -\omega^3+\omega^1_2
\end{pmatrix}.
\end{equation}
\end{lemma}
\begin{proof}
For simplicity, let $\sigma_0:=\1, \sigma_1:=\i, \sigma_2:=\j^{'},
\sigma_3:=\k^{'}$.
By applying the chain rule,
\begin{align*}
de_\alpha(g)&=d(g_1\sigma_\alpha g_2^t)\\
            &=(dg_1)\sigma_\alpha g_2^t+g_1\sigma_\alpha dg_2^t\\
            &=g_1\{g_1^{-1}(dg_1)\sigma_\alpha+\sigma_\alpha(g_2^{-1}dg_2)^t\}
g_2^t.
\end{align*}
On the other hand,
$$de_\alpha=\omega^\beta_\alpha e_\beta=\omega^\beta_\alpha
g_1\sigma_\beta g_2.$$
Hence, we have the equation
$$(g_1^{-1}dg_1)\sigma_\alpha+\sigma_\alpha(g_2^{-1}dg_2)^t=\omega^\beta_\alpha
\sigma_\beta.$$
The equations \eqref{eq:m-c1} and \eqref{eq:m-c2} follow from this equation.
\end{proof}
\begin{lemma}
If the frame field $\{e_\alpha: \alpha=0,1,2,3\}$ is parametrized by
${\rm SL}_2\R\times{\rm SL}_2\R$ via the action $\nu$, then the pull
back $g^{-1}dg$ of Maurer-Cartan form $\Omega=(\omega^\beta_\alpha)$
can be written as the following equation in the Lie algebra
$\s_2\R\oplus\s_2\R$:
$$g^{-1}dg=g_1^{-1}dg_1\oplus (dg_2^{-1})g_2,$$
where
\begin{equation}
\label{eq:m-c3}
g_1^{-1}dg_1=\frac{1}{2}\begin{pmatrix}
\omega^3+\omega^1_2 & \omega^1+\omega^2-\omega^1_3-\omega^2_3\\
-\omega^1+\omega^2-\omega^1_3+\omega^2_3 & -\omega^3-\omega^1_2
\end{pmatrix}
\end{equation}
and
\begin{equation}
\label{eq:m-c4}
(dg_2^{-1})g_2=\frac{1}{2}\begin{pmatrix}
\omega^3-\omega^1_2 & \omega^1+\omega^2+\omega^1_3+\omega^2_3\\
-\omega^1+\omega^2+\omega^1_3-\omega^2_3 & -\omega^3+\omega^1_2
\end{pmatrix}.
\end{equation}
\end{lemma}
\begin{proof}
Similar to the proof of Lemma \ref{lem:m-c}, we get the equation
$$(g_1^{-1}dg_1)\sigma_\alpha+\sigma_\alpha(dg_2^{-1})g_2=\omega^\beta_\alpha\sigma_\beta$$
and the equations \eqref{eq:m-c3} and \eqref{eq:m-c4} then follow.
\end{proof}
\section{Timelike $\cmc$ Surfaces in AdS $3$-Space $\H^3_1(-1)$ and
  Integrable Systems}
\label{sec:intsystem}
Let $M$ be a simply-connected open and orientable $2$-manifold and
$\varphi: M\longrightarrow\H^3_1(-1)$ a timelike conformal immersion.

By using a double covering induced by the group action $\mu$, we can
find lift $\Phi=(\Phi_1,\Phi_2)$ (called a coordinate frame) of $\F$ to ${\rm
  SL}_2\R\times{\rm SL}_2\R$:
$$\mu(\Phi)(\1,\i,\j^{'},\k^{'})=\F.$$
That is,
the lifted framing $\Phi=(\Phi_1,\Phi_2): M\longrightarrow{\rm
  SL}_2\R\times{\rm SL}_2\R$ satisfies\footnote{Here, we use the same
  $\mu$ for both Lie group action and group representation.}
\begin{align*}
\mu(\Phi)(\1)&=\Phi_1\1 \Phi_2^t=\varphi,\\
\mu(\Phi)(\i)&=\Phi_1\i \Phi_2^t=e^{-\frac{\omega}{2}}\varphi_x,\\
\mu(\Phi)(\j^{'})&=\Phi_1\j^{'}\Phi_2^t=e^{-\frac{\omega}{2}}\varphi_y,\\
\mu(\Phi)(\k^{'})&=\Phi_1\k^{'}\Phi_2^t=N.
\end{align*}
Then
\begin{equation}
\label{eq:basis1}
\varphi_u=e^{\frac{\omega}{2}}\Phi_1\begin{pmatrix}
0 & 1\\
0 & 0
\end{pmatrix}\Phi_2^t
\end{equation}
and
\begin{equation}
\label{eq:basis2}
\varphi_v=e^{\frac{\omega}{2}}\Phi_1\begin{pmatrix}
0 & 0\\
1 & 0
\end{pmatrix}\Phi_2^t.
\end{equation}
Similarly, by using a double covering induced by the group action $\nu$, we can
find lift $\Psi=(\Psi_1,\Psi_2)$ of $\F$ to ${\rm SL}_2\R\times{\rm SL}_2\R$:
$$\nu(F)(\1,\i,\j^{'},\k^{'})=\F.$$
The lifted framing $\Psi=(\Psi_1,\Psi_2): M\longrightarrow{\rm
  SL}_2\R\times{\rm SL}_2\R$ satisfies
\begin{align*}
\nu(\Psi)(\1)&=\Psi_1\1 \Psi_2^{-1}=\varphi,\\
\nu(\Psi)(\i)&=\Psi_1\i \Psi_2^{-1}=e^{-\frac{\omega}{2}}\varphi_x,\\
\nu(\Psi)(\j^{'})&=\Psi_1\j^{'}\Psi_2^{-1}=e^{-\frac{\omega}{2}}\varphi_y,\\
\nu(\Psi)(\k^{'})&=\Psi_1\k^{'}\Psi_3^{-1}=N.
\end{align*}
Then
\begin{equation}
\label{eq:basis3}
\varphi_u=e^{\frac{\omega}{2}}\Psi_1\begin{pmatrix}
0 & 1\\
0 & 0
\end{pmatrix}\Psi_2^{-1}
\end{equation}
and
\begin{equation}
\label{eq:basis4}
\varphi_v=e^{\frac{\omega}{2}}\Psi_1\begin{pmatrix}
0 & 0\\
1 & 0
\end{pmatrix}\Psi_2^{-1}.
\end{equation}

Let $\mathfrak s:=(\varphi,\varphi_u,\varphi_v,N)$. Then $\mathfrak s$
defines a moving frame
on the immersed surface $\varphi$ and satisfy the following
Gau{\ss}-Weingarten equations:
\begin{equation}
\label{eq:g-w}
\mathfrak s_u=\mathfrak s{\mathcal U},\ \mathfrak s_v=\mathfrak
s{\mathcal V},
\end{equation}
where
\begin{equation}
\mathcal U=\begin{pmatrix}
0 & 0 & \frac{1}{2}e^\omega & 0\\
1 & \omega_u & 0 & -H\\
0 & 0 & 0 & -2Qe^{-\omega}\\
0 & Q & \frac{1}{2}e^\omega H & 0
\end{pmatrix},\ \mathcal V=\begin{pmatrix}
0 & \frac{1}{2}e^\omega & 0 & 0\\
0 & 0 & 0 & -2Re^{-\omega}\\
1 & 0 & \omega_v & -H\\
0 & \frac{1}{2}e^\omega H & R & 0
\end{pmatrix}.
\end{equation}
The integrability condition of the Gau{\ss}-Weingarten equation is the
Gau{\ss}-Mainardi-Codazzi equation
$$\mathcal V_u-\mathcal U_v+[\mathcal U,\mathcal V]=0.$$
This Gau{\ss}-Mainardi-Codazzi equation is equivalent to
\begin{eqnarray}
\label{eq:g-c}
\omega_{uv}+\frac{1}{2}e^\omega(H^2-1)-2QRe^{-\omega}=0,\\
\label{eq:meancurv}
H_u=2e^{-\omega}Q_v,\ H_v=2e^{-\omega}R_u.
\end{eqnarray}
\begin{remark}
From the equations \eqref{eq:meancurv}, we see that a
timelike surface $\varphi: M\longrightarrow\H^3_1(-1)$ has constant mean
curvature if and only if $R_u=Q_v=0$. In this case, $R$ is said to be
\emph{Lorentz antiholomorphic} and $Q$ is said to be \emph{Lorentz
  holomorphic}, respectively.
\end{remark}
Each component framing $\Phi_1$ and $\Phi_2$ of $\Phi$ satisfy the
following Lax equations:
\begin{eqnarray*}
(\Phi_1)_u=\Phi_1\U_1, & (\Phi_1)_v=\Phi_1\V_1;\\
(\Phi_2)_u=\Phi_2\U_2, & (\Phi_2)_v=\Phi_2\V_2,
\end{eqnarray*}
where
\begin{eqnarray}
\label{eq:lax1}
\U_1&=&\begin{pmatrix}
\frac{\omega_u}{4} & \frac{1}{2}e^\frac{\omega}{2}(H+1)\\
-e^{-\frac{\omega}{2}}Q & -\frac{\omega_u}{4}
\end{pmatrix},\\
\V_1&=&\begin{pmatrix}
-\frac{\omega_v}{4} & e^{-\frac{\omega}{2}}R\\
-\frac{1}{2}e^\frac{\omega}{2}(H-1) & \frac{\omega_v}{4}
\end{pmatrix},\\
\U_2&=&\begin{pmatrix}
-\frac{\omega_u}{4} & e^{-\frac{\omega}{2}}Q\\
-\frac{1}{2}e^\frac{\omega}{2}(H-1) & \frac{\omega_u}{4}
\end{pmatrix},\\
\V_2&=&\begin{pmatrix}
\frac{\omega_v}{4} & \frac{1}{2}e^\frac{\omega}{2}(H+1)\\
-e^{-\frac{\omega}{2}}R & -\frac{\omega_v}{4}
\end{pmatrix}.
\end{eqnarray}
The compatibility conditions $(\Phi_1)_{uv}=(\Phi_1)_{vu}$ and
$(\Phi_2)_{uv}=(\Phi_2)_{vu}$ give the Maurer-Cartan equations
\begin{equation}
(\U_1)_v-(\V_1)_u-[\U_1,\V_1]=0
\end{equation}
and
\begin{equation}
(\U_2)_v-(\V_2)_u-[\U_2,\V_2]=0.
\end{equation}
Each of these two Maurer-Cartan equations is also equivalent to the
Gau{\ss}-Mainardi-Codazzi equations \eqref{eq:g-c} and \eqref{eq:meancurv}.

Each component framing $\Psi_1$ and $\Psi_2$ of $\Psi$ satisfy the
following Lax equations:
\begin{eqnarray*}
(\Psi_1)_u=\Psi_1\U_1, & (\Psi_1)_v=\Psi_1\V_1;\\
(\Psi_2)_u=\Psi_2\U_2, & (\Psi_2)_v=\Psi_2\V_2,
\end{eqnarray*}
where
\begin{eqnarray}
\label{eq:lax2}
\U_1&=&\begin{pmatrix}
\frac{\omega_u}{4} & \frac{1}{2}e^\frac{\omega}{2}(H+1)\\
-e^{-\frac{\omega}{2}}Q & -\frac{\omega_u}{4}
\end{pmatrix},\\
\V_1&=&\begin{pmatrix}
-\frac{\omega_v}{4} & e^{-\frac{\omega}{2}}R\\
-\frac{1}{2}e^\frac{\omega}{2}(H-1) & \frac{\omega_v}{4}
\end{pmatrix},\\
\U_2&=&\begin{pmatrix}
\frac{\omega_u}{4} & \frac{1}{2}e^\frac{\omega}{2}(H-1)\\
-e^{-\frac{\omega}{2}}Q & -\frac{\omega_u}{4}
\end{pmatrix},\\
\V_2&=&\begin{pmatrix}
-\frac{\omega_v}{4} & e^{-\frac{\omega}{2}}R\\
-\frac{1}{2}e^\frac{\omega}{2}(H+1) & \frac{\omega_v}{4}
\end{pmatrix}.
\end{eqnarray}
The compatibility conditions $(\Psi_1)_{uv}=(\Psi_1)_{vu}$ and
$(\Psi_2)_{uv}=(\Psi_2)_{vu}$ give the Maurer-Cartan equations
\begin{equation}
(\U_1)_v-(\V_1)_u-[\U_1,\V_1]=0
\end{equation}
and
\begin{equation}
(\U_2)_v-(\V_2)_u-[\U_2,\V_2]=0.
\end{equation}
Again, each of these two Maurer-Cartan equations is equivalent to the
Gau{\ss}-Mainardi-Codazzi equations \eqref{eq:g-c} and \eqref{eq:meancurv}.

We now have the following representation formulae for timelike $\cmc$
  surfaces in $\H^3_1(-1)$.
\begin{theorem}
\label{thm:sym}

Let $M$ be a simply-connected region in Minkowski plane
$\E^2_1=(\R^2(u,v),dudv)$.
\begin{enumerate}
\item Let $\Phi=(\Phi_1,\Phi_2): M\longrightarrow{\rm
  SL}_2\R\times{\rm SL}_2\R$ be a solution
to the following Lax equations:
\begin{equation}
\label{eq:lax3}
\begin{array}{cc}
(\Phi_1)_u=\Phi_1\U_1, & (\Phi_1)_v=\Phi_1\V_1;\\
(\Phi_2)_u=\Phi_2\U_2, & (\Phi_2)_v=\Phi_2\V_2,
\end{array}
\end{equation}
where
\begin{eqnarray*}
\U_1=\begin{pmatrix}
\frac{\omega_u}{4} & \frac{1}{2}e^{\frac{\omega}{2}}(H+1)\\
-e^{-\frac{\omega}{2}}Q & -\frac{\omega_u}{4}
\end{pmatrix} & \V_1=\begin{pmatrix}
-\frac{\omega_v}{4} & e^{-\frac{\omega}{2}}R\\
-\frac{1}{2}e^{\frac{\omega}{2}}(H-1) & \frac{\omega_v}{4}
\end{pmatrix},\\
\U_2=\begin{pmatrix}
-\frac{\omega_u}{4} & e^{-\frac{\omega}{2}}Q\\
-\frac{1}{2}e^{\frac{\omega}{2}}(H-1) & \frac{\omega_u}{4}
\end{pmatrix}, & \V_2=\begin{pmatrix}
\frac{\omega_v}{4} & \frac{1}{2}e^{\frac{\omega}{2}}(H+1)\\
-e^{-\frac{\omega}{2}}R & -\frac{\omega_v}{4}
\end{pmatrix}.
\end{eqnarray*}
Then $\varphi:=\mu(\Phi)(\1)=\Phi_1\Phi_2^t:
M\longrightarrow\H^3_1(-1)$ defines a
timelike $\cmc$ $H$ immersion into $\H^3_1(-1)$.
\item Let $\Psi=(\Psi_1,\Psi_2): M\longrightarrow{\rm
  SL}_2\R\times{\rm SL}_2\R$ be a
solution to the following Lax equations:
\begin{equation}
\label{eq:lax4}
\begin{array}{cc}
(\Psi_1)_u=\Psi_1\U_1, & (\Psi_1)_v=\Psi_1\V_1;\\
(\Psi_2)_u=\Psi_2\U_2, & (\Psi_2)_v=\Psi_2\V_2,
\end{array}
\end{equation}
where
\begin{eqnarray*}
\U_1=\begin{pmatrix}
\frac{\omega_u}{4} & \frac{1}{2}e^{\frac{\omega}{2}}(H+1)\\
-e^{-\frac{\omega}{2}}Q & -\frac{\omega_u}{4}
\end{pmatrix}, & \V_1=\begin{pmatrix}
-\frac{\omega_v}{4} & e^{-\frac{\omega}{2}}R\\
-\frac{1}{2}e^{\frac{\omega}{2}}(H-1) & \frac{\omega_v}{4}
\end{pmatrix},\\
\U_2=\begin{pmatrix}
\frac{\omega_u}{4} & \frac{1}{2}e^{\frac{\omega}{2}}(H-1)\\
-e^{-\frac{\omega}{2}}Q & -\frac{\omega_v}{4}
\end{pmatrix}, & \V_2=\begin{pmatrix}
-\frac{\omega_v}{4} & e^{-\frac{\omega}{2}}R\\
-\frac{1}{2}e^{\frac{\omega}{2}}(H+1) & \frac{\omega_v}{4}
\end{pmatrix}.
\end{eqnarray*}
Then $\psi:=\nu(\Psi)(\1)=\Psi_1\Psi_2^{-1}:
M\longrightarrow\H^3_1(-1)$ defines a
timelike $\cmc$ $H$ immersion into $\H^3_1(-1)$.
\end{enumerate}
\end{theorem}

Let $H_e,\ H_s,\ H_h$ be the constant mean curvatures of timelike surfaces in
Minkowski $3$-space $\E^3_1$, de Sitter $3$-space $\S^3_1(1)$ and
anti-de Sitter $3$-space $\H^3_1(-1)$, resp. Then these timelike $\cmc$
surfaces in each space-form satisfy the following Gau{\ss}-Mainardi-Codazzi
equations:
\begin{eqnarray*}
\omega_{uv}+\frac{1}{2}H_e^2e^\omega-2QRe^{-\omega}=0 & \mbox{\rm Lorentz
  $3$-space $\E^3_1$ case},\\
\omega_{uv}+\frac{1}{2}(H_s^2+1)e^\omega-2QRe^{-\omega}=0 & \mbox{\rm
  de Sitter $3$-space $\S^3_1(1)$ case},\\
\omega_{uv}+\frac{1}{2}(H_h^2-1)e^\omega-2QRe^{-\omega}=0 & \mbox{\rm
  anti-de Sitter $3$-space $\H^3_1(-1)$ case}.
\end{eqnarray*}
By comparing these Gau{\ss}-Mainardi-Codazzi equations,
we can deduce the \emph{Lawson-Guichard correspondence} between timelike
$\cmc$ $H_e$ surfaces in $\E^3_1$, timelike $\cmc$ $H_s$ surfaces in
$\S^3_1(1)$ and timelike $\cmc$ $H_h$ surfaces in $\H^3_1(-1)$, which
satisfy the same Gau{\ss}-Mainardi-Codazzi equation. Such $\cmc$ surfaces are
called \emph{cousins} of each other. In particular, we see that there
is a bijective correspondence between timelike minimal surfaces
($H=0$) in $\E^3_1$ and timelike $\cmc$ $\pm 1$ in $\H^3_1(-1)$. For
this reason, we are mainly interested in timelike $\cmc$ $\pm 1$ surfaces in
$\H^3_1(-1)$.  If $H=\pm 1$ then the Gau{\ss}-Mainardi-Codazzi
equations \eqref{eq:g-c} become
\begin{equation}
\label{eq:g-c2}
\left\{
\begin{array}{c}
\omega_{uv}-2e^{-\omega}RQ=0,\\
R_u=Q_v=0.
\end{array}
\right.
\end{equation}
For more details about the Lawson-Guichard corresponce,
please see the appendix I (section \ref{sec:lawson}).
\begin{remark}
Note that one can normalize the Hopf pairs $Q$ and $R$ if $M$ has real
distinct principal curvatures or imaginary principal curvatures
everywhere. For example, $Q=\pm R=1$ reduces the Gau{\ss}-Mainardi-Codazzi
equations \eqref{eq:g-c} to the Liouville equation $\omega_{uv}=\mp
2e^{-\omega}$.
\end{remark}
Since the sign of $H$ depends upon the orientation of a surface
(i.e., the orientation of the unit normal vector field $N$),
hereafter we consider only $H=1$ case.
\begin{corollary}

Let $M$ be a simply-connected $2$-manifold.
\begin{enumerate}
\item Let $\Phi=(\Phi_1,\Phi_2): M\longrightarrow{\rm
  SL}_2\R\times{\rm SL}_2\R$ be solutions
to the following Lax equations:
\begin{eqnarray*}
(\Phi_1)_u=\Phi_1\U_1, & (\Phi_1)_v=\Phi_1\V_1;\\
(\Phi_2)_u=\Phi_2\U_2, & (\Phi_2)_v=\Phi_2\V_2,
\end{eqnarray*}
where
\begin{eqnarray*}
\U_1=\begin{pmatrix}
\frac{\omega_u}{4} & e^{\frac{\omega}{2}}\\
-e^{-\frac{\omega}{2}}Q & -\frac{\omega_u}{4}
\end{pmatrix} & \V_1=\begin{pmatrix}
-\frac{\omega_v}{4} & e^{-\frac{\omega}{2}}R\\
0 & \frac{\omega_v}{4}
\end{pmatrix},\\
\U_2=\begin{pmatrix}
-\frac{\omega_u}{4} & e^{-\frac{\omega}{2}}Q\\
0 & \frac{\omega_u}{4}
\end{pmatrix}, & \V_2=\begin{pmatrix}
\frac{\omega_v}{4} & e^{\frac{\omega}{2}}\\
-e^{-\frac{\omega}{2}}R & -\frac{\omega_v}{4}
\end{pmatrix}.
\end{eqnarray*}
Then $\varphi:=\Phi_1\Phi_2^t: M\longrightarrow\H^3_1(-1)$ defines a
timelike $\cmc$ $1$ immersion into $\H^3_1(-1)$.
\item Let $\Psi=(\Psi_1,\Psi_2): M\longrightarrow{\rm
  SL}_2\R\times{\rm SL}_2\R$ be
solutions to the following Lax equations:
\begin{eqnarray*}
(\Psi_1)_u=\Psi_1\U_1, & (\Psi_1)_v=\Psi_1\V_1;\\
(\Psi_2)_u=\Psi_2\U_2, & (\Psi_2)_v=\Psi_2\V_2,
\end{eqnarray*}
where
\begin{eqnarray*}
\U_1=\begin{pmatrix}
\frac{\omega_u}{4} & e^{\frac{\omega}{2}}\\
-e^{-\frac{\omega}{2}}Q & -\frac{\omega_u}{4}
\end{pmatrix}, & \V_1=\begin{pmatrix}
-\frac{\omega_v}{4} & e^{-\frac{\omega}{2}}R\\
0 & \frac{\omega_v}{4}
\end{pmatrix},\\
\U_2=\begin{pmatrix}
\frac{\omega_u}{4} & 0\\
-e^{-\frac{\omega}{2}}Q & -\frac{\omega_v}{4}
\end{pmatrix}, & \V_2=\begin{pmatrix}
-\frac{\omega_v}{4} & e^{-\frac{\omega}{2}}R\\
e^{-\frac{\omega}{2}} & \frac{\omega_v}{4}
\end{pmatrix}.
\end{eqnarray*}
Then $\psi:=\Psi_1\Psi_2^{-1}: M\longrightarrow\H^3_1(-1)$ defines a
timelike $\cmc$ $1$ immersion into $\H^3_1(-1)$.
\end{enumerate}
\end{corollary}
\section{Timelike Surfaces of Constant Mean Curvature $\pm 1$ in AdS
  $3$-Space $\H^3_1(-1)$ via the Group Action $\mu$}
\label{sec:bryant}
In \cite{Ho}, J.~Q.~Hong gave a Bryant type representation formula
for timelike surfaces of constant mean curvature $1$ in
$\H^3_1(-1)$. In this section, we reproduce
J.~Q.~Hong's Bryant type reprepsentation formula in a more general context.

Let $F: M\longrightarrow{\rm SL}_2\R\times{\rm SL}_2\R$ be a lift of
$\F$ to ${\rm SL}_2\R\times{\rm SL}_2\R$ via the Lie group action
$\mu$, i.e., $\mu(F)(\1,\i,\j^{'},\k^{'})=\F$.
Let $\Omega:=\Omega_1\oplus\Omega_2\in\s_2\R\oplus\s_2\R$, where
$\Omega_i=F_i^{-1}dF_i\in\s_2\R,\ i=1,2$. The Gau{\ss} and Mainardi-Codazzi
equations are equivalent to Maurer-Cartan equation
\begin{equation}
\label{eq:m-c}
d\Omega+\Omega\wedge\Omega=0,
\end{equation}
which is the null curvature (integrability) condition of the
Maurer-Cartan form $\Omega$.
Let $F^*$ denote the pull-back map $F^*:
(\s_2\R)^*\oplus(\s_2\R)^*\longrightarrow T^*M$. The
Maurer-Cartan forms $\Omega_i=F_i^{-1}dF_i\in\s_2\R,\ i=1,2$ can be
written as the following equations:
\begin{eqnarray}
\label{eq:m-c5}
F_1^{-1}dF_1&=&\frac{1}{2}F^*\begin{pmatrix}
\omega^3+\omega_2^1 & \omega^1+\omega^2-\omega_3^1-\omega_3^2\\
-\omega_1+\omega_2-\omega_3^1+\omega_3^2 & -\omega^3-\omega_2^1
\end{pmatrix},\\
\label{eq:m-c6}
F_2^{-1}dF_2&=&\frac{1}{2}F^*\begin{pmatrix}
\omega^3-\omega_2^1 & -\omega^1+\omega^2+\omega_3^1-\omega_3^2\\
\omega^1+\omega^2+\omega_3^1+\omega_3^2 & -\omega^3+\omega_2^1
\end{pmatrix}.
\end{eqnarray}
\begin{definition}
\label{def:null}
Let $M$ be a $2$-manifold. A map $F: M\longrightarrow{\rm SL}_2\R$ is
said to be \emph{null} if $F^*(\phi)=0$, or equivalently
$\det(F^{-1}dF)=0$, where $\phi$ is the quadratic Cartan-Killing form
$\phi=-8\det(g^{-1}dg)$.
\end{definition}
\begin{theorem}[A Bryant type representation formula for timelike $\cmc$
    $\pm 1$ Surfaces in $\H^3_1(-1)$]
\label{thm:bryant}
Let $M$ be an open $2$-manifold and $F=(F_1,F_2): M\longrightarrow{\rm
  SL}_2\R\times{\rm SL}_2\R$ an immersion such that
\begin{enumerate}
\item $F_1$ is Lorentz holomorphic, i.e., $(F_1)_v=0$ and
$F_2$ is Lorentz antiholomorphic, i.e., $(F_2)_u=0$,
\item $F_1$ and $F_2$ are null, i.e., $\det F_1^{-1}dF_1=\det F_2^{-1}dF_2=0$.
\end{enumerate}
Then
\begin{equation}
\label{eq:bryant}
\varphi:=\mu(F)(\1)=F_1F_2^t
\end{equation}
is a smooth conformal timelike immersion into $\H^3_1(-1)$ with $\cmc$
$\pm 1$. Conversely, let $M$ be an oriented and
simply-connected Lorentzian $2$-manifold with globally defined null
coordinates\footnote{In Lorentzian case, the so-called \emph{Riemann
    Mapping Theorem} or \emph{K\"{o}be Uniformization Theorem} does not
  hold. So, globally defined null coordinates do not exist, in
  general, on a simply-connected Lorentzian $2$-manifold.}. If $\varphi:
M\longrightarrow\H^3_1(-1)$ is a smooth conformal timelike immersion
with $\cmc$ $\pm 1$, then there exists an immersion $F=(F_1,F_2):
M\longrightarrow{\rm SL}_2\R\times{\rm SL}_2\R$ such that $F_1,\ F_2$
satisfy the conditions {\rm (1)}, {\rm (2)}, and $\varphi=F_1F_2^t$.
\end{theorem}
\begin{proof}
Let $\omega^+=\omega^1+\omega^2$, $\omega^-=-\omega^1+\omega^2$,
$\pi^+=\omega_3^1+\omega_3^2$, $\pi^-=-\omega_3^1+\omega_3^2$. Also let
\begin{align*}
F^*(\omega^3+\omega_2^1)&=2\alpha_1,\ F^*(\omega^+-\pi^+)=2\beta_1,\
F^*(\omega^-+\pi^-)=2\gamma_1\\
F^*(\omega^3-\omega_2^1)&=2\alpha_2,\ F^*(\omega^--\pi^-)=2\beta_2,\
F^*(\omega^++\pi^+)=2\gamma_2.
\end{align*}
Then the Maurer-Cartan equations \eqref{eq:m-c1} and \eqref{eq:m-c2} become
\begin{align}
F_1^{-1}dF_1&=\frac{1}{2}F^*\begin{pmatrix}
\omega^3+\omega_2^1 & \omega^+-\pi^+\\
\omega^-+\pi^- & -\omega^3-\omega_2^1
\end{pmatrix}=\begin{pmatrix}
\alpha_1 & \beta_1\\
\gamma_1 & -\alpha_1
\end{pmatrix},\\
F_2^{-1}dF_2&=\frac{1}{2}F^*\begin{pmatrix}
\omega^3-\omega_2^1 & \omega^--\pi^-\\
\omega^++\pi^+ & -\omega^3+\omega_2^1
\end{pmatrix}=\begin{pmatrix}
\alpha_2 & \beta_2\\
\gamma_2 & -\alpha_2
\end{pmatrix}.
\end{align}
Note that $F^*(\omega^3)=\alpha_1+\alpha_2$,
$F^*(\omega^+)=\beta_1+\gamma_2$, $F^*(\omega^-)=-(\beta_2+\gamma_1)$.
Since $\det F_1^{-1}dF_1=\det F_2^{-1}dF_2=0$,
$\alpha_1^2+\beta_1\gamma_1=0$ and $\alpha_2^2+\beta_2\gamma_2=0$.

Denote by $ds^2$ the metric in $\H^3_1(-1)$ induced by the canonical
semi-\\Riemannian metric in $\E^4_2$.
If we regard $e_0$ as a map $e_0: {\rm SL}_2\R\times{\rm
  SL}_2\R\longrightarrow\H^3_1(-1)$ given by $e_0(g)=\mu(g)(\1),\
g\in{\rm SL}_2\R\times{\rm
  SL}_2\R$, then
\begin{align*}
e_0^*(ds^2)&=<de_0,de_0>\\
           &=<\omega^1e_1+\omega^2e_2+\omega^3e_3,\omega^1e_1+\omega^2e_2+
\omega^3e_3>\\
           &=-(\omega^1)^2+(\omega^2)^2+(\omega^3)^2
\end{align*}
defines an indefinite metric in the oriented orthonormal frame bundle
of $\H^3_1(-1)$.
Since $\varphi=e_0\circ F$,
\begin{align*}
ds^2_\varphi&=\varphi^*(ds^2)\\
            &=F^*\circ e_0^*(ds^2)\\
            &=F^*((\omega^3)^2+\{-(\omega^1)^2+(\omega^2)^2\})\\
            &=F^*(\omega^3)^2+F^*(\omega^+)F^*(\omega^-)\\
            &=2\alpha_1\alpha_2+\beta_1\beta_2+\gamma_1\gamma_2.
\end{align*}
Since $F$ is an immersion, the last expression defines a metric.

We now show that for the immersion $\varphi:
M\longrightarrow\H^3_1(-1)$, $H\equiv 1$. Let $U\subset M$ be a
simply-connected open set in which there exists a null coordinate
system $(u,v)$ such that $ds^2_\varphi=e^{\omega}dudv$ for some
real-valued function $\omega: U\longrightarrow\R$ defined on $U$.
Clearly, $M$ is covered by such open sets. Let
$\eta=e^{\frac{\omega}{2}}du$ and $\xi=e^{\frac{\omega}{2}}dv$.
There exist functions $A_1,\ A_2,\ B_1,\ B_2,\ C_1,\ C_2$ such that
\begin{align*}
F^*(\omega^3+\omega_2^1)&=2A_1\eta,\ F^*(\omega^3-\omega^1_2)=2A_2\xi,\\
F^*(\omega^+-\pi^+)&=2B_1\eta,\ F^*(\omega^--\pi^-)=2B_2\xi,\\
F^*(\omega^-+\pi^-)&=2C_1\eta,\ F^*(\omega^++\pi^+)=2C_2\xi.
\end{align*}
We, then, have equations
$$A_1^2+B_1C_1=0,\ A_2+B_2C_2=0,\ \rm{and}\ 2A_1A_2+B_1B_2+C_1C_2=1.$$
In the open set $U$,
\begin{eqnarray*}
ds^2_\varphi&=&(2A_1A_2+B_1B_2+C_1C_2)\eta\xi\\
            &=&\eta\xi.
\end{eqnarray*}
Since $A_1^2+B_1C_1=0$ and $A_2+B_2C_2=0$,
there exist smooth functions $p_1,\ p_2,\ q_1,\ q_2$ defined in $U$
(unique up to replacement by $(-p_1,-q_1)$ and $(-p_2,-q_2)$
respectively) such that
$$\begin{array}{cccccc}
A_1&=&p_1q_1, & A_2&=&p_2q_2\\
B_1&=&q_1^2,  & B_2&=&q_2^2\\
C_1&=&-p_1^2, & C_2&=&-p_2^2
\end{array}$$
and
$$2p_1q_1p_2q_2+q_1^2q_2^2+p_1^2p_2^2=(p_1p_2+q_1q_2)^2=1.$$
By the continuity of $p_1p_2+q_1q_2$, either
$p_1p_2+q_1q_2=1$ in $U$
or $p_1p_2+q_1q_2=-1$ in $U$. Without loss of
generality, we may
assume that $p_1p_2+q_1q_2=1$. Now define a map $h:
U\longrightarrow{\rm SL}_2\R$ by
$h=\begin{pmatrix}
p_1 & -q_2\\
q_1 & p_2
\end{pmatrix}$.
Then $e_0(F_1h,F_2(h^{-1})^t)=e_0(F_1,F_2)=F_1F_2^t.$
\begin{align*}
(F_1h)^{-1}d(F_1h)&=(h^{-1}F_1^{-1})((dF_1)h+F_1dh)\\
                  &=h^{-1}(F_1^{-1}dF_1)h+h^{-1}dh\\
                  &=h^{-1}\begin{pmatrix}
p_1q_1 & -p_1^2\\
q_1^2 & -p_1q_1
\end{pmatrix}\eta h+h^{-1}dh\\
                  &=\begin{pmatrix}
p_2dp_1+q_2dq_1 & -p_2dq_2+q_2dp_2-\eta\\
-q_1dp_1+p_1dq_1 & q_1dq_2+p_1dp_2
\end{pmatrix}.
\end{align*}
Similarly,
$$(F_2(h^{-1})^t)^{-1}d(F_2(h^{-1})^t)=\begin{pmatrix}
p_1dp_2+q_1dq_2 & -p_1dq_1+q_1dp_1-\xi\\
-q_2dp_2+p_2dq_2 & q_2dq_1+p_2dp_1
\end{pmatrix}.$$
It then follows that $(Fh)^*(\omega^3)=0$, $(Fh)^*(\omega^+)=-\eta$
and $(Fh)^*(\omega^-)=-\xi$, i.e., $(F_1h,F_2(h^{-1})^t):
U\longrightarrow{\rm SL}_2\R\times{\rm SL}_2\R$ is an oriented framing
in $U$ along the immersion $\varphi=e_0\circ
F=e_0(F_1h,F_2(h^{-1})^t)=F_1F_2^t$.

The $1$-form $\dis\frac{1}{2}(Fh)^*(\omega^-+\pi^-)=p_1dq_1-q_1dp_1$
can be written
$$p_1dq_1-q_1dp_1=\left\{\begin{array}{ccc}
p_1^2d\left(\frac{q_1}{p_1}\right) & {\rm where} & p_1\ne 0\\
-q_1^2d\left(\frac{p_1}{q_1}\right) & {\rm where} & q_1\ne 0
\end{array}\right..$$ Hence, $\frac{1}{2}(Fh)^*(\omega^-+\pi^-)$ is a
$1$-form of type $(1,0)$, i.e., a multiple of the $1$-form
$\eta$. Similarly, $(Fh)^*(\omega^++\pi^+)$ is a $1$-form of type
$(0,1)$, i.e., a multiple of the $1$-form $\xi$.

Since $(Fh)^*(\omega^+)=-\eta$
and $(Fh)^*(\omega^-)=-\xi$, by the equation \eqref{eq:meanc}, one can
easily see that:
\begin{enumerate}
\item If $\varphi$ satisfies $H=1$ in $U$, then
$(Fh)^*(\omega^-+\pi^-)$ is a $1$-form of type $(1,0)$ and
$(Fh)^*(\omega^++\pi^+)$ is a $1$-form of type $(0,1)$.
\item If $(Fh)^*(\omega^-+\pi^-)$ is a $1$-form of type
$(1,0)$ or $(Fh)^*(\omega^++\pi^+)$ is a $1$-form of type $(0,1)$,
then $\varphi$ satisfies $H=1$ in $U$.
\end{enumerate}
Therefore, we conclude that $H=1$ in $U$.

Conversely, let $M$ be an oriented and simply-connected Lorentzian
$2$-\\manifold with globally defined null coordinates $(u,v)$. Let
$\varphi: M\longrightarrow\H^3_1(-1)$ be a smooth conformal timelike
immersion into $\H^3_1(-1)$ with $\cmc$ $1$. Then
$ds^2_\varphi=e^{\omega}dudv$ for some real-valued function $\omega:
M\longrightarrow\R$. Let $\eta=e^{\frac{\omega}{2}}du$ and
$\xi=e^{\frac{\omega}{2}}dv$. Then one can choose a lifting
$g=(g_1,g_2): M\longrightarrow{\rm SL}_2\R\times{\rm SL}_2\R$ such
that the associated frame field $\{e_0(g)\}$ is adapted with
$g^*(\omega^+)=-\eta$ and $g^*(\omega^-)=-\xi$. Since
$g^*(\omega^3)=0$,
\begin{align*}
g_1^{-1}dg_1&=\frac{1}{2}g^*\begin{pmatrix}
\omega^1_2 & \omega^+-\pi^+\\
\omega^-+\pi^- & -\omega^1_2
\end{pmatrix}\\
            &=\frac{1}{2}g^*\begin{pmatrix}
\omega^1_2 & -\omega^+-\pi^+\\
\omega^-+\pi^- & -\omega^1_2
\end{pmatrix}+\eta\begin{pmatrix}
0 &-1\\
0 & 0
\end{pmatrix},\\
g_2^{-1}dg_2&=\frac{1}{2}g^*\begin{pmatrix}
-\omega^1_2 & \omega^--\pi^-\\
\omega^++\pi^+ & \omega^1_2
\end{pmatrix}\\
            &=\frac{1}{2}g^*\begin{pmatrix}
-\omega^1_2 & -\omega^--\pi^-\\
\omega^++\pi^+ & \omega^1_2
\end{pmatrix}+\xi\begin{pmatrix}
0 & -1\\
0 & 0
\end{pmatrix}.
\end{align*}
Let $\zeta=\frac{1}{2}g^*\begin{pmatrix}
\omega^1_2 & -\omega^+-\pi^+\\
\omega^-+\pi^- & -\omega^1_2
\end{pmatrix}\in\s_2\R$. Then $d\zeta=-\zeta\wedge\zeta$. The equation
$d\zeta=-\zeta\wedge\zeta$ satisfies the integrability condition;
hence, by the Frobenius Theorem, there exists a smooth map $h:
M\longrightarrow{\rm SL}_2\R$ such that $\zeta=h^{-1}dh$. Since
$h\in{\rm SL}_2\R$, it can be written
$$h=\begin{pmatrix}
p_1 & -q_2\\
q_1 & p_2
\end{pmatrix},\ p_1,p_2,q_1,q_2\in\R.$$
Set $F_1:=g_1h^{-1}$. Then
\begin{align*}
F_1^{-1}dF_1&=(g_1h^{-1})^{-1}d(g_1h^{-1})\\
            &=h(g_1^{-1}dg_1)h^{-1}+hdh^{-1}\\
            &=h\left[\zeta+\begin{pmatrix}
0 & -1\\
0 & 0
\end{pmatrix}\eta\right]h^{-1}+hdh^{-1}\\
            &=h\begin{pmatrix}
0 & -1\\
0 & 0
\end{pmatrix}h^{-1}\eta\\
            &=\begin{pmatrix}
p_1q_1 & -p_1^2\\
q_1^2 & -p_1q_1
\end{pmatrix}\eta.
\end{align*}
The differential $d$ can be written as $d=\partial^{'}+\partial{''}$,
where $\partial^{'}$ is the Lorentz holomorphic part and
$\partial^{''}$ is the Lorentz antiholomorphic part. Since
$\partial^{''}F_1=\dis\frac{\partial F_1}{\partial v}dv=0$, $F_1$ is
Lorentz holomorphic.

Set $F_2=g_2h^t$. Then
\begin{align*}
F_2^{-1}dF_2&=(g_2h^t)^{-1}d(g_2h^t)\\
            &=(h^t)^{-1}(g_2^{-1}dg_2)h^t+(h^t)^{-1}dh^t\\
            &=(h^t)^{-1}\left[-\zeta^t+\begin{pmatrix}
0 & -1\\
0 & 0
\end{pmatrix}\xi\right]h^t+(h^t)^{-1}dh^t\\
            &=(h^t)^{-1}\begin{pmatrix}
0 & -1\\
0 & 0
\end{pmatrix}h^t\xi\\
            &=\begin{pmatrix}
p_2q_2 & -p_2^2\\
q_2^2 & -p_2q_2
\end{pmatrix}\xi.
\end{align*}
Since $\partial^{'}F_2=\dis\frac{\partial F_2}{\partial u}du=0$, $F_2$
is Lorentz antiholomorphic. Finally,
$$F_1F_2^t=g_1h^{-1}(g_2h^t)^t=g_1g_2^t=\varphi.$$
\end{proof}
\begin{remark}
\label{rem:bryant} Note that, in Theorem \ref{thm:bryant},
$\varphi=F_1F_2^t$ has $\cmc$ $1$ ($\cmc$ $-1$) if the framing $F$
is orientation preserving (orientation reversing). In order to prove
Theorem \ref{thm:bryant} for orientation reversing framing $F$, one
needs to take
\begin{align*}
F_1^{-1}dF_1&=\begin{pmatrix}
p_1q_1 & p_1^2\\
-q_1^2 & -p_1q_1
\end{pmatrix}\eta,\ F_2^{-1}dF_2=\begin{pmatrix}
p_2q_2 & p_2^2\\
-q_2^2 & -p_2q_2
\end{pmatrix}\xi,\ {\rm and}\\
h&=\begin{pmatrix}
q_2 & -p_1\\
p_2 & q_1
\end{pmatrix}\in{\rm SL}_2\R\ \mbox{\rm in the proof.}
\end{align*}
\end{remark}
\begin{remark}
\begin{align*}
ds^2_\varphi&=\varphi^*(ds^2)\\
            &=<d\varphi,d\varphi>\\
            &=<d(F_1F_2^t),d(F_1F_2)^t>\\
            &=-\det\{F_1^{-1}dF_1+(F_2^{-1}dF_2)^t\}.
\end{align*}
So, $\varphi$ does not assume degenerate points if and only if
$$\det\{F_1^{-1}dF_1+(F_2^{-1}dF_2)^t\}\ne 0.$$
\end{remark}
\section{Timelike Surfaces of Constant Mean Curvature $\pm 1$ in AdS
  $3$-Space $\H^3_1(-1)$ via the Group Action $\nu$}
\label{sec:bryant2}
Let $F: M\longrightarrow{\rm SL}_2\R\times{\rm SL}_2\R$ be a lift of
$\F$ to ${\rm SL}_2\R\times{\rm SL}_2\R$ via the Lie group action
$\nu$, i.e., $\nu(F)(\1,\i,\j^{'},\k^{'})=\F$.
Let $\Omega:=\Omega_1\oplus\Omega_2\in\s_2\R\oplus\s_2\R$, where
$\Omega_1=F_1^{-1}dF_1,\ \Omega_2=(dF_2^{-1})F_2\in\s_2\R$. The
Gau{\ss} and Mainardi-Codazzi
equations are equivalent to Maurer-Cartan equation
\begin{equation}
d\Omega+\Omega\wedge\Omega=0,
\end{equation}
which is the null curvature (integrability) condition of the
Maurer-Cartan form $\Omega$.
Let $F^*$ denote the pull-back map $F^*:
(\s_2\R)^*\oplus(\s_2\R)^*\longrightarrow T^*M$. The
Maurer-Cartan forms $\Omega_1=F_1^{-1}dF_1,\ \Omega_2=(dF_2^{-1})F_2$ can be
written as the following equations:
\begin{eqnarray}
\label{eq:m-c7}
F_1^{-1}dF_1&=&\frac{1}{2}F^*\begin{pmatrix}
\omega^3+\omega_2^1 & \omega^1+\omega^2-\omega_3^1-\omega_3^2\\
-\omega_1+\omega_2-\omega_3^1+\omega_3^2 & -\omega^3-\omega_2^1
\end{pmatrix},\\
\label{eq:m-c8}
(dF_2^{-1})F_2&=&\frac{1}{2}F^*\begin{pmatrix}
\omega^3-\omega_2^1 & \omega^1+\omega^2+\omega_3^1+\omega_3^2\\
-\omega^1+\omega^2+\omega_3^1-\omega_3^2 & -\omega^3+\omega_2^1
\end{pmatrix}.
\end{eqnarray}
\begin{theorem}[A Bryant type representation formula for timelike $\cmc$
    $\pm 1$ Surfaces in $\H^3_1(-1)$]
\label{thm:bryant2}
Let $M$ be an open $2$-manifold
and $F=(F_1,F_2): M\longrightarrow{\rm SL}_2\R\times{\rm SL}_2\R$ an
immersion such that
\begin{enumerate}
\item $F_1$ is Lorentz holomorphic, i.e., $(F_1)_v=0$ and
$F_2$ is Lorentz antiholomorphic, i.e., $(F_2)_u=0$,
\item $\det F_1^{-1}dF_1=\det(dF_2^{-1})F_2=0$.
\end{enumerate}
Then
\begin{equation}
\label{eq:bryant2}
\psi:=\nu(F)(\1)=F_1F_2^{-1}
\end{equation}
is a smooth conformal timelike immersion
into $\H^3_1(-1)$ with $\cmc$ $\pm 1$. Conversely, let $M$ be an oriented and
simply-connected Lorentzian $2$-manifold with globally defined null
coordinates. If $\psi:
M\longrightarrow\H^3_1(-1)$ is a smooth conformal timelike immersion
with $\cmc$ $\pm 1$, then there exists an immersion $F=(F_1,F_2):
M\longrightarrow{\rm SL}_2\R\times{\rm SL}_2\R$ such that $F_1,\ F_2$
satisfy the conditions {\rm (1)}, {\rm (2)}, and $\psi=F_1F_2^{-1}$.
\end{theorem}
\begin{proof}
We use the same $\omega^+,\omega^-,\pi^+,\pi^-$ and
$\alpha_i,\beta_i,\gamma_i,\ i=1,2$ as defined in the proof of Theorem
\ref{thm:bryant}. Then
$$F_1^{-1}dF_1=\begin{pmatrix}
\alpha_1 & \beta_1\\
\gamma_1 & -\alpha_1
\end{pmatrix}\ {\rm and}\ (dF_2^{-1})F_2=\begin{pmatrix}
\alpha_2 & \gamma_2\\
\beta_2 & -\alpha_2
\end{pmatrix}.$$
Here,
$$ds^2_\psi=F^*(\omega^3)^2+F^*(\omega^+)F^*(\omega^-)=2\alpha_1\alpha_2+\beta_1\beta_2+\gamma_1\gamma_2$$
defines an induced metric of $\psi$ since
$F$ is an immersion.

Let $U\subset M$ be a simply-connected open set in which there exists
a null coordinate system $(u,v)$ such that $ds^2_\psi=e^{\omega}dudv$,
for some real-valued function $\omega: U\longrightarrow\R$. $M$ is
covered by such open sets. Let $\eta=e^{\frac{\omega}{2}}du$ and
$\xi=e^{\frac{\omega}{2}}dv$. Then, by exactly the same argument in the
proof of Theorem \ref{thm:bryant}, there exists a smooth map $h:
U\longrightarrow{\rm SL}_2\R$ given by $h=\begin{pmatrix}
p_1 & -q_2\\
q_1 & p_2
\end{pmatrix}$
and
$$e_0\circ Fh=e_0(F_1h,F_2h)=F_1h(F_2h)^{-1}=F_1F_2^{-1}=e_0\circ F.$$
\begin{align*}
(F_1h)^{-1}d(F_1h)&=h^{-1}(F_1^{-1}dF_1)h+h^{-1}dh\\
                  &=\begin{pmatrix}
p_2dp_1+q_2dq_1 & -p_2dq_2+q_2dp_2-\eta\\
-q_1dp_1+p_1dq_1 & q_1dq_2+p_1dp_2
\end{pmatrix}
\end{align*}
and
\begin{align*}
d(F_2h)^{-1}F_2h&=(dh^{-1})h+h^{-1}[(dF_2^{-1})F_2]h\\
                &=\begin{pmatrix}
p_1dp_2+q_1dq_2 & -q_2dp_2+p_2dq_2\\
-p_1dq_1+q_1dp_1-\xi & q_2dq_1+p_2dp_1
\end{pmatrix}.
\end{align*}
It then follows that $(Fh)^*(\omega^3)=0$, $(Fh)^*(\omega^+)=-\eta$
and $(Fh)^*(\omega^-)=-\xi$. Hence, $Fh: U\longrightarrow{\rm
  SL}_2\R\times{\rm SL}_2\R$ is an oriented framing in $U$ along the
immersion $\psi=F_1F_2^{-1}$ and $\psi$ satisfies $H=1$ in $U$.

Conversely, let $M$ be an oriented and simply-connected Lorentzian
$2$-\\manifold with globally defined null coordinates. Let $\psi: M
\longrightarrow\H^3_1(-1)$ be a smooth conformal timelike immersion
with $\cmc$ $1$. There exist a null coordinate system $(u,v)$ in $M$
such that $ds^2_\psi=e^{\omega}dudv$ for some real-valued function
$\omega: M\longrightarrow\R$. Let $\eta=e^{\frac{\omega}{2}}du$ and
$\xi=e^{\frac{\omega}{2}}dv$. Then one can choose a lifting $g:
M\longrightarrow{\rm SL}_2\R\times{\rm SL}_2\R$ such that the
associated frame field $\{e_0(g)\}$ is adapted with
$g^*(\omega^+)=-\eta$ and $g^*(\omega^-)=-\xi$. Since
$g^*(\omega^3)=0$,
\begin{align*}
g_1^{-1}dg_1&=\frac{1}{2}g^*\begin{pmatrix}
\omega^1_2 & \omega^+-\pi^+\\
\omega^-+\pi^- & -\omega^1_2
\end{pmatrix}\\
            &=\frac{1}{2}g^*\begin{pmatrix}
\omega^1_2 & -\omega^+-\pi^+\\
\omega^-+\pi^- & -\omega^1_2
\end{pmatrix}+\eta\begin{pmatrix}
0 & -1\\
0 & 0
\end{pmatrix}
\end{align*}
and
\begin{align*}
(dg_2^{-1})g_2&=\frac{1}{2}g^*\begin{pmatrix}
-\omega^1_2 & \omega^++\pi^+\\
\omega^--\pi^- & \omega^1_2
\end{pmatrix}\\
               &=\frac{1}{2}g^*\begin{pmatrix}
-\omega^1_2 & \omega^++\pi^+\\
-\omega^--\pi^- & \omega^1_2
\end{pmatrix}+\xi\begin{pmatrix}
0 & 0\\
-1 & 0
\end{pmatrix}.
\end{align*}
Let $\zeta=\dis\frac{1}{2}g^*\begin{pmatrix}
\omega^1_2 & -\omega^+-\pi^+\\
\omega^-+\pi^- & -\omega^1_2
\end{pmatrix}\in\s_2\R$. Then $d\zeta=-\zeta\wedge\zeta$ and this
equation satisfies the integrability condition; hence, by the
Frobenius Theorem, there exists a smooth map $h: M\longrightarrow{\rm
  SL}_2\R$ such that $\zeta=h^{-1}dh$. Since $h\in{\rm SL}_2\R$, it
can be written $h=\begin{pmatrix}
p_1 & -q_2\\
q_1 & p_2
\end{pmatrix},\ p_1,p_2,q_1,q_2\in\R$.
Set $F_1=g_1h^{-1}$. Then
$$F_1^{-1}dF_1=\begin{pmatrix}
p_1q_1 & -p_1^2\\
q_1^2 & -p_1q_1
\end{pmatrix}\eta$$
and $F_1$ is Lorentz holomorphic.
Set $F_2=g_2h^{-1}$. Then
\begin{align*}
(dF_2^{-1})F_2&=(dh)h^{-1}+h[(dg_2^{-1})g_2]h^{-1}\\
              &=(dh)h^{-1}+h\left[-\zeta+\begin{pmatrix}
0 & 0\\
-1 & 0
\end{pmatrix}\xi\right]h^{-1}\\
              &=h\begin{pmatrix}
0 & 0\\
-1 & 0
\end{pmatrix}h^{-1}\xi\\
              &=\begin{pmatrix}
p_2q_2 & q_2^2\\
-p_2^2 & -p_2q_2
\end{pmatrix}\xi
\end{align*}
and $F_2$ is Lorentz antiholomorphic.
Finally,
$$F_1F_2^{-1}=g_1h^{-1}(g_2h^{-1})^{-1}=g_1g_2^{-1}=\psi.$$
\end{proof}
\begin{remark}
Note
that
\begin{align*}
ds^2_\psi&=<d\psi,d\psi>\\
         &=<d(F_1F_2^{-1}),d(F_1F_2^{-1})>\\
         &=-\det\{F_1^{-1}dF_1+(dF_2^{-1})F_2\}.
\end{align*}
So, $\psi$ does not assume degenerate points if and only if
$$\det\{F_1^{-1}dF_1+(dF_2^{-1})F_2\}\ne 0.$$
\end{remark}
\section{Timelike Minimal Surfaces in $\E^3_1$ and the Classical
  Gau{\ss} Map}
\label{sec:tlminimal}
Recall that the Lie group $G\cong{\rm SL}_2\R$ acts isometrically on Lorentz
$3$-space $\E^3_1$ via the Ad-action. The ${\rm Ad}(G)$-orbit of
$\k^{'}=\begin{pmatrix}
1 & 0\\
0 & -1
\end{pmatrix}$ is a pseudosphere or de Sitter $2$-space:
$$\S^2_1(1)=\{(x_1,x_2,x_3)\in\E^3_1: -x_1^2+x_2^2+x_3^2=1\}.$$
The Ad-action of $G$ on $\S^2_1(1)$ is transitive as well. The
isotropy subgroup of $G$ at $\k^{'}$ is the indefinite orthogonal
group ${\rm SO}_1(2)=\{x_0\1+x_3\k^{'}: x_0^2-x_3^2=1\}$. Thus,
$\S^2_1(1)$ can be identified with the symmetric space
$${\rm SL}_2\R/{\rm SO}_1(2)=\left\{h\begin{pmatrix}
1 & 0\\
0 & -1
\end{pmatrix}h^{-1}: h\in{\rm SL}_2\R\right\}.$$ The orthogonal
subgroup ${\rm SO}_1(2)$ can
be regarded as the hyperbola $H^1_0$ in a Lorentz plane
$\E^2_1(x_0,x_3)$. (This is a Lorentz analogue of $S^1\subset\E^2$.)
Note that the group $H^1_0$ is isomorphic to the multiplicative group
$\R^\times=(\R\setminus\{0\},\times)$.

Let $N=(0,0,1)\ {\rm and}\ S=(0,0,-1)\in\S^2_1(1)$ be the \emph{north}
and \emph{south} pole of $\S^2_1(1)$. Let $\wp_+:
\S^2_1(1)\setminus\{x_3=-1\}\longrightarrow\E^2_1\setminus H^1_0$ be
the stereographic projection from the south pole $S=(0,0,-1)$, where
$H^1_0=\{(x_1,x_2)\in\E^2_1: -x_1^2+x_2^2=-1\}$. Then
\begin{align*}
\wp_+(x_1,x_2,x_3)&=\left(\frac{x_1}{1+x_3},\frac{x_2}{1+x_3}\right)\\
                  &\cong\left(\frac{x_1+x_2}{1+x_3},\frac{-x_1+x_2}{1+x_3}\right)\in\E^2_1(u,v).
\end{align*}
Let $\wp_-:
\S^2_1(1)\setminus\{x_3=1\}\longrightarrow\E^2_1\setminus H^1_0$ be
the stereographic projection from the north pole $S=(0,0,1)$. Then
\begin{align*}
\wp_-(x_1,x_2,x_3)&=\left(\frac{x_1}{1-x_3},\frac{x_2}{1-x_3}\right)\\
                  &\cong\left(\frac{x_1+x_2}{1-x_3},\frac{-x_1+x_2}{1-x_3}\right)\in\E^2_1(u,v).
\end{align*}
Note that the classical Gau{\ss} map (i.e., the unit normal vector field)
$N$ of timelike surfaces in $\E^3_1$ is mapped into de Sitter $2$-space
$\S^2_1(1)$. Thus, the image of classical Gau{\ss} map can be
represented by the matrices $h\begin{pmatrix}
1 & 0\\
0 & -1
\end{pmatrix}h^{-1}$, $h\in{\rm SL}_2\R$. If the timelike surface
preserves the orientation, then $h=\begin{pmatrix}
p_1 & -q_2\\
q_1 & p_2
\end{pmatrix}\in{\rm SL}_2\R$. So,
$$h\begin{pmatrix}
1 & 0\\
0 & -1
\end{pmatrix}h^{-1}=\begin{pmatrix}
p_1p_2-q_1q_2 & 2p_1q_2\\
2p_2q_1 & -p_1p_2+q_1q_2
\end{pmatrix}\in\S^2_1(1)$$
and $\wp_-(h\begin{pmatrix}
1 & 0\\
0 & -1
\end{pmatrix}h^{-1})=\dis\left(\frac{p_1}{q_1},\frac{p_2}{q_2}\right)\in
\E^2_1(u,v)$. If the timelike surface reverses the orientation, then
$h=\begin{pmatrix}
q_2 & -p_1\\
p_2 & q_1
\end{pmatrix}\in{\rm SL}_2\R$. So,
$$h\begin{pmatrix}
1 & 0\\
0 & -1
\end{pmatrix}h^{-1}=\begin{pmatrix}
-p_1p_2+q_1q_2 & 2p_1q_2\\
2p_2q_1 & p_1p_2-q_1q_2
\end{pmatrix}\in\S^2_1(1)$$
and $\wp_+(h\begin{pmatrix}
1 & 0\\
0 & -1
\end{pmatrix}h^{-1})=\dis\left(\frac{p_1}{q_1},\frac{p_2}{q_2}\right)\in
\E^2_1(u,v)$.

Here, we recall the following \emph{Weierstra{\ss} formula}
for a timelike minimal surface $\psi: M\longrightarrow\E^3_1$ with
data $(q,f(q))$ and $(r,g(v))$:
\begin{equation}
\label{eq:wf}
\psi_u=(\frac{1}{2}(1+q^2),-\frac{1}{2}(1-q^2),-q)f(u),\
\psi_v=(-\frac{1}{2}(1+r^2),-\frac{1}{2}(1-r^2),-r)g(v).
\end{equation}
The induced metric of $\psi$ is
$$ds^2_\psi=(1+qr)^2f(u)g(v)dudv.$$
\begin{remark}
\label{rem:gaussmap}
The ordered pair $(q,r)$ coincides with the projected Gau{\ss} map
$\wp_-\circ N$ of a timelike minimal surface with data $(q,f(u))$ and
$(r,g(v))$.
\end{remark}
\begin{remark}
In \cite{I-T}, J.~Inoguchi and M.~Toda studied the construction of
timelike minimal surfaces via loop group method. Their
\emph{normalized Wierstra{\ss}
formula} for a timelike minimal surface $\psi: M\longrightarrow\E^3_1$
with data $(q,r)$ is
\begin{equation}
\label{eq:nwf}
\psi_u=(\frac{1}{2}(1+q^2),-\frac{1}{2}(1-q^2),-q),\
\psi_v=(-\frac{1}{2}(1+r^2),-\frac{1}{2}(1-r^2),-r)
\end{equation}
and the induced metric of $\psi$ is
$$ds^2_\psi=(1+qr)^2dudv.$$

In \cite{I-T}, the signs of coordinate functions in
  $\psi_u$ and $\psi_v$ are different. The reason is, in \cite{I-T},
$(x_1,x_2,x_3)\in\E^3_1$ is identified with the matrix\\
$\begin{pmatrix}
-x_3 & -x_1+x_2\\
x_1+x_2 & x_3
\end{pmatrix}$, while in this paper it is identified with\\
$\begin{pmatrix}
x_3 & x_1+x_2\\
-x_1+x_2 & -x_3
\end{pmatrix}$.

Originally, this formula was obtained by M. A. Magid in \cite{Ma}. However, in
\cite{Ma}, the geometric meaning of the data $(q,r)$ is not
clarified. In \cite{I-T}, the data $(q,r)$ are retrieved from the
\emph{normalized potential} in their construction. Moreover, $q$ and
$r$ are the primitive functions of the coefficients $Q$ and $R$, resp., of Hopf
pairs. Note that this is locally true. In general,
\begin{align}
q_u&=\frac{Q}{f(u)},\\
r_v&=\frac{R}{g(v)}.
\end{align}
As is mentioned in Remark \ref{rem:gaussmap}, $(q,r)$ is the projected
Gau{\ss} map $\wp_-\circ N$ of a timelike minimal surface given by the
Weierstra{\ss} formula \eqref{eq:nwf}.
\end{remark}
\section{Lawson-Guichard Correspondence between Timelike $\cmc$ $\pm
  1$ Surfaces in $\H^3_1(-1)$ and Timelike Minimal Surfaces in $\E^3_1$}
\label{sec:corresp}
In Section \ref{sec:intsystem}, we discussed the Lawson-Guichard
correspondence between timelike $\cmc$ surfaces in three different
semi-Riemannian space forms $\E^3_1$, $\S^3_1(1)$ and $\H^3_1(-1)$. In
particular, there is a one-to-one correspondence between timelike $\cmc$
$\pm 1$ surfaces in $\H^3_1(-1)$ and timelike minimal surfaces in $\E^3_1$. In
this section, we give such bijective correspondence explicitly using
the Bryant type representation formulae in Sections \ref{sec:bryant}
and \ref{sec:bryant2}.

Let $\varphi: M\longrightarrow\H^3(-1)$ be a timelike $\cmc$ $-1$
surface. Then, by Theorem \ref{thm:bryant} (or by Theorem
\ref{thm:bryant2}), there exists a smooth immersion $F=(F_1,F_2):
M\longrightarrow{\rm SL}_2\R\times{\rm SL}_2\R$ such that
\begin{enumerate}
\item $F_1$ is Lorentz holomorphic and $F_2$ is Lorentz
antiholomorphic.
\item $\det(F_1^{-1}dF_1)=\det(F_2^{-1}dF_2)=0$\\
(or ($2^{'}$) $\det(F_1^{-1}dF_1)=\det((dF_2^{-1})F_2)=0$).
\item $\varphi=F_1F_2^t$ (or ($3^{'}$) $\varphi=F_1F_2^{-1}$).
\end{enumerate}

As we have seen in the proof of Theorem \ref{thm:bryant} (or Theorem
\ref{thm:bryant2}), locally in an open set $U\subset M$,
\begin{align*}
F_1^{-1}dF_1&=\begin{pmatrix}
p_1q_1 & -p_1^2\\
q_1^2 & -p_1q_1
\end{pmatrix}\eta\\
            &=\begin{pmatrix}
\frac{p_1}{q_1} & -\frac{p_1^2}{q_1^2}\\
1 & -\frac{p_1}{q_1}
\end{pmatrix}q_1^2\eta
\end{align*}
and similarly,
$$F_2^{-1}dF_2=\begin{pmatrix}
\frac{p_2}{q_2} & -\frac{p_2^2}{q_2^2}\\
1 & -\frac{p_2}{q_2}
\end{pmatrix}q_2^2\xi$$
or
$$(dF_2^{-1})F_2=\begin{pmatrix}
\frac{p_2}{q_2} & 1\\
-\frac{p_2^2}{q_2^2} & -\frac{p_2}{q_2}
\end{pmatrix}q_2^2\xi.$$
Let $q:=\dis\frac{p_1}{q_1}$, $f(u):=q_1^2$,
$r:=\dis\frac{p_2}{q_2}$, and $g(v):=q_2^2$. Then the Weierstra{\ss} formula
\eqref{eq:wf} defines a timelike
minimal surface $\psi: U\longrightarrow\E^3_1$.
\begin{align*}
ds_\varphi^2&=\eta\xi\\
            &=e^\omega dudv\\
            &=e^\omega(1+qr)^2f(u)g(v)dudv.
\end{align*}
So, the induced metric $ds^2_\varphi$ of the timelike $\cmc$ $-1$ surface
$\varphi$ is conformal to $ds^2_\psi$.

Conversely, assume that a timelike minimal surface $\psi:
M\longrightarrow\E^3_1$ is given by the Weierstra{\ss}
formula \eqref{eq:wf} with data $(q,f(u))$ and $(r,g(v))$. Consider
the following system of differential equations:
\begin{equation}
\label{eq:BUY}
F_1^{-1}dF_1=\begin{pmatrix}
q & -q^2\\
1 & -q
\end{pmatrix}f(u)du,\ F_2^{-1}dF_2=\begin{pmatrix}
r & -r^2\\
1 & -r
\end{pmatrix}g(v)dv.
\end{equation}
Since these equations satisfy the integrability condition, there exists
a solution $(F_1,F_2): M\longrightarrow{\rm SL}_2\R\times{\rm SL}_2\R$
satisfying the conditions (1) and (2). By Theorem \ref{thm:bryant},
$\varphi:=F_1F_2^t: M\longrightarrow\H^3_1(-1)$ defines a timelike $\cmc$
$-1$ surface in $\H^3_1(-1)$. Similarly, the system of differential
equations:
\begin{equation}
\label{eq:BUY2}
F_1^{-1}dF_1=\begin{pmatrix}
q & -q^2\\
1 & -q
\end{pmatrix}f(u)du,\ (dF_2^{-1})F_2=\begin{pmatrix}
r & 1\\
-r^2 & -r
\end{pmatrix}g(v)dv
\end{equation}
satisfies the integrability condition; hence there exists a solution
$(F_1,F_2): M\longrightarrow{\rm SL}_2\R\times{\rm SL}_2\R$ satisfying
conditions (1) and ($2^{'}$). By Theorem \ref{thm:bryant2},
$\varphi:=F_1F_2^{-1}$ defines a timelike $\cmc$ $-1$ surface in
$\H^3_1(-1)$.
\section{The Hyperbolic Gau{\ss} Map of Timelike Surfaces in $\H^3_1(-1)$}
\label{sec:hypgauss}

Let $\varphi: M\longrightarrow\H^3_1(-1)$ be an oriented timelike surface in
$\H^3_1(-1)$. At each base point $e_0=\varphi(m)\in\H^3_1(-1)$,
$e_3\in{\rm T}_{e_0}\H^3_1(-1)$ is an oriented unit normal vector to
the tangent plane $\varphi_*({\rm T}_mM)$. The oriented normal
geodesic in $\H^3_1(-1)$ emanating from $e_0$, which is tangent to the
normal vector $e_3(\varphi(m))$ asymptotically approaches to the
\emph{null cone} $\N=\{\u\in\E^4_2: <\u,\u>=0\}$ at exactly two points
$[e_0+e_3],\ [e_0-e_3]\in\N$. The orientation allows us to name
$[e_0+e_3]$ the initial point and $[e_0-e_3]$ the terminal point.

Define a map $G: M\longrightarrow\N$ by $G(m)=[e_0+e_3](m)$ for each
$m\in M$. This map is an analogue of the hyperbolic Gau{\ss}
map\footnote{The hyperbolic Gau{\ss} map was introduced by C. Epstein
  (\cite{Ep}) and was used by R. L. Bryant in his study of $\cmc$ $1$
  surfaces in hyperbolic $3$-space $\H^3(-1)$ (\cite{Br}).} of
surfaces in hyperbolic $3$-space $\H^3_1(-1)$. The map will still be
called the \emph{hyperbolic Gau{\ss} map} here.

Let $\varphi: M\longrightarrow\H^3_1(-1)$ be a timelike surface in
$\H^3_1(-1)$ and $d\sigma^2$ denote the induced metric on $\N^3$. Then
$$d\sigma^2_\varphi:=(e_0+e_3)^*(d\sigma^2)=<d(e_0+e_3),d(e_0+e_3)>=-(\omega^1+
\omega^1_3)^2+(\omega^2+\omega^2_3)^2.$$
\begin{proposition}
\label{prop:hypgauss}
The hyperbolic Gau{\ss} map $[e_0+e_3]:
  M\longrightarrow\N$ ($[e_0-e_3]:
  M\longrightarrow\N$) of a timelike surface $\varphi:
  M\longrightarrow\H^3_1(-1)$ is conformal if and only if $\varphi$
  satisfies $H=1$ ($H=-1$) or $\varphi$ is totally umbilic.
\end{proposition}
\begin{proof}
We will assume the same settings in the proof of Proposition
\ref{prop:curv}. Then
\begin{eqnarray*}
d\sigma^2_\varphi&=&-(\omega^1+\omega^1_3)^2+(\omega^2+\omega^2_3)^2\\
                 &=&-((1-h_{11})\omega^1+h_{12}\omega^2)^2+(-h_{12}\omega^1+
(1-h_{22})\omega^2)^2\\
                 &=&(2H(H-1)-K)ds^2_\varphi-\\
                 &{}&(H-1)((h_{11}-h_{22})(\omega^1)^2-4h_{12}
\omega^1\otimes\omega^2+(h_{11}-h_{22})(\omega^2)^2).
\end{eqnarray*}
Thus, $[e_0+e_3]$ is conformal, i.e., $d\sigma^2_\varphi$ is a
multiple of $ds^2_\varphi$ if and only if
$(H-1)(h_{11}-h_{22})=(H-1)h_{12}=0$.

If $H=1$, then $d\sigma^2_\varphi=-Kds^2_\varphi$. Suppose $H\ne 1$.
Let $U=\{m\in M: H(m)\ne 1\}$. Then $U$ is clearly open in $M$.
Since $H\ne 1$ on $U$, $h_{11}=h_{22}$ and $h_{12}=0$ on $U$. The
second fundamental form $I\!I$ is then
\begin{align*}
I\!I&=-h_{11}(\omega^1)^2+2h_{12}\omega^1\otimes\omega^2+h_{22}(\omega^2)^2\\
  &=h_{11}(-(\omega^1)^2+(\omega^2)^2)\\
  &=HI\ {\rm on}\ U.
\end{align*}
By comparing with the equation \eqref{eq:sff}, we see that the Hopf
differential $\Q=0$ on $U$, i.e., $\varphi(U)$ is totally umbilic.
Note that $H$ must be constant on $U$, since $Q=R=0$. Let $V$ be a
connected component of $U$. Since $H$ is constant on $V$ and $H$ is
continuous on $M$, $H$ is constant on $\bar V$. This implies that
$H\ne 1$ on $\bar V$ and so $\bar V\subset U$. Since $V$ is
connected, so is $\bar V$. However, $V$ is a connected component;
thus, $V=\bar V$. It then follows from the connectedness of $M$ that
$M=V$. Therefore, $\varphi$ is totally umbilic on $M$. The converse
is trivial.
\end{proof}
\begin{remark}
Note that $[e_0+e_3]: M\longrightarrow\N$ ($[e_0-e_3]:
M\longrightarrow\N$) is the hyperbolic Gau{\ss} map of a timelike
surface $\varphi: M\longrightarrow\H^3_1(-1)$ if $\varphi$ preserves
(reverses) the orientation.
\end{remark}
\begin{remark}
If $\varphi$ satisfies $\pm 1$ \emph{and} is totally umbilic, then
the hyperbolic Gau{\ss} map is constant, since
$d\sigma^2_\varphi=0$. It is also shown in Section \ref{sec:hol} that
if $\varphi$ satisfies $\pm 1$ and is totally umbilic, then the
(projected) hyperbolic Gau{\ss} map is both Lorentz holomorphic and
antiholomorphic; hence it is constant. By the equation
\eqref{eq:gauss}, the Gau{\ss}ian curvature
$K=H^2-1=0$. Thus, $\varphi$ may be regarded as a \emph{horosphere}
type surface in $\H^3_1(-1)$.
\end{remark}
The null cone $\N$ satisfies the quadric equation
$$-(x_0)^2-(x_1)^2+(x_2)^2+(x_3)^2=0.$$
If $x_0\ne 0$, then the above equation can be written
$$-\left(\frac{x_1}{x_0}\right)^2+\left(\frac{x_2}{x_0}\right)^2+
\left(\frac{x_3}{x_0}\right)^2=1,$$
i.e., $\N$ can be locally identified with de Sitter $2$-space
$\S^2_1(1)$. With this identification, the hyperbolic Gau{\ss} map can
be mapped into de Sitter $2$-space $\S^2_1(1)$. Thus, we may be able
to relate the hyperbolic Gau{\ss} map of a timelike $\cmc$ $\pm 1$ surface in
$\H^3_1(-1)$ and the Gau{\ss} map of corresponding timelike minimal
surface in $\E^3_1$. This relationship is discussed in the next section.

Let $(x_0,x_1,x_2,x_3)\in\N$. Denote by $[x_0,x_1,x_2,x_3]$ the null
line generated by the null vector $(x_0,x_1,x_2,x_3)$. By using
nonhomogeneous coordinates,
\begin{align*}
[x_0,x_1,x_2,x_3]&=\left[1,\frac{x_1}{x_0},\frac{x_2}{x_0},\frac{x_3}{x_0}
\right]\ {\rm provided}\ x_0\ne 0\\
                 &\cong\left(\frac{x_1}{x_0},\frac{x_2}{x_0},\frac{x_3}{x_0}
\right)\in\S^2_1(1)\\
                 &\cong\left(\frac{x_1}{x_0+x_3},\frac{x_2}{x_0+x_3}\right)\in
\E^2_1\ \mbox{\rm via the projection}\ \wp_+
\end{align*}
or
$$[x_0,x_1,x_2,x_3]\cong\left(\frac{x_1}{x_0-x_3},\frac{x_2}{x_0-x_3}\right)\in\E^2_1\
\mbox{\rm via the projection}\ \wp_-.$$
Finally, we have the identification:
\begin{equation}
\label{eq:nullline}
[x_0,x_1,x_2,x_3]\cong\left(\frac{x_1+x_2}{x_0+x_3},\frac{-x_1+x_2}{x_0+x_3}
\right)\in\E^2_1(u,v)
\end{equation}
or
\begin{equation}
\label{eq:nullline2}
[x_0,x_1,x_2,x_3]\cong\left(\frac{x_1+x_2}{x_0-x_3},\frac{-x_1+x_2}{x_0-x_3}
\right)\in\E^2_1(u,v).
\end{equation}
\section{The Hyperbolic Gau{\ss} Map and the Secondary Gau{\ss} Map}
\label{sec:secgauss}
Let $\varphi: M\longrightarrow\H^3_1(-1)$ be a timelike $\cmc$ $-1$
surface in $\H^3_1(-1)$. Then by Theorem \ref{thm:bryant} (or by Theorem
\ref{thm:bryant2}), there exists a smooth immersion $F=(F_1,F_2):
M\longrightarrow{\rm SL}_2\R\times{\rm SL}_2\R$ such that
\begin{enumerate}
\item $F_1$ is Lorentz holomorphic and $F_2$ is Lorentz
antiholomorphic.
\item $\det(F_1^{-1}dF_1)=\det(F_2^{-1}dF_2)=0$\\
(or ($2^{'}$) $\det(F_1^{-1}dF_1)=\det((dF_2^{-1})F_2)=0$).
\item $\varphi=F_1F_2^t$ (or ($3^{'}$) $\varphi=F_1F_2^{-1}$).
\end{enumerate}

Let $F_1=\begin{pmatrix}
F_{11} & F_{12}\\
F_{13} & F_{14}
\end{pmatrix}$ and $F_2=\begin{pmatrix}
F_{21} & F_{22}\\
F_{23} & F_{24}
\end{pmatrix}$. Then
\begin{align*}
(e_0+e_3)(F)&=F_1(\1+\k^{'})F_2^t\\
            &=2\begin{pmatrix}
F_{11} & F_{12}\\
F_{13} & F_{14}
\end{pmatrix}\begin{pmatrix}
1 & 0\\
0 & 0
\end{pmatrix}\begin{pmatrix}
F_{21} & F_{23}\\
F_{22} & F_{24}
\end{pmatrix}\\
            &=2\begin{pmatrix}
F_{11}F_{21} & F_{11}F_{23}\\
F_{13}F_{21} & F_{13}F_{23}
\end{pmatrix}.
\end{align*}
By the identification \eqref{eq:nullline2},
\begin{align*}
[(e_0+e_3)(F)]&=\left[\begin{pmatrix}
F_{11}F_{21} & F_{11}F_{23}\\
F_{13}F_{21} & F_{13}F_{23}
\end{pmatrix}\right]\\
              &\cong\left(\frac{F_{11}F_{23}}{F_{13}F_{23}},\frac{F_{13}F_{21}}{F_{13}F_{23}}\right)\\
              &=\left(\frac{F_{11}}{F_{13}},\frac{F_{21}}{F_{23}}\right)\in\E^2_1(u,v).
\end{align*}
Similarly,
\begin{align*}
[(e_0-e_3)(F)]&=[F_1(\1-\k^{'})F_2^t]\\
              &\cong\left(\frac{F_{12}}{F_{14}},\frac{F_{22}}{F_{24}}\right)\in\E^2_1(u,v).
\end{align*}
If $(e_0\pm e_3)(F)=F_1(\1\pm\k^{'})F_2^{-1}$, then
\begin{align*}
[(e_0+e_3)(F)]&\cong\left(\frac{F_{11}}{F_{13}},-\frac{F_{24}}{F_{22}}\right)\in\E^2_1(u,v),\\
[(e_0-e_3)(F)]&\cong\left(\frac{F_{12}}{F_{14}},-\frac{F_{23}}{F_{21}}\right)\in\E^2_1(u,v).
\end{align*}
Locally,
\begin{align*}
(e_0+e_3)(Fh)&=F_1h(\1+\k^{'})(F_2(h^{-1})^t)^t\\
             &=2F_1h\begin{pmatrix}
1 & 0\\
0 & 0
\end{pmatrix}h^{-1}F_2^t\\
             &=2F_1\begin{pmatrix}
p_1 & -q_2\\
q_1 & p_2
\end{pmatrix}\begin{pmatrix}
1 & 0\\
0 & 0
\end{pmatrix}\begin{pmatrix}
p_2 & q_2\\
-q_1 & p_1
\end{pmatrix}\begin{pmatrix}
F_{21} & F_{23}\\
F_{22} & F_{24}
\end{pmatrix}.
\end{align*}
The last expression is simplified to the matrix:
$$2\begin{pmatrix}
(F_{11}p_1+F_{12}q_1)(F_{21}p_2+F_{22}q_2) &
  (F_{11}p_1+F_{12}q_1)(F_{23}p_2+F_{24}q_2)\\
(F_{13}p_1+F_{14}q_1)(F_{21}p_2+F_{22}q_2) &
  (F_{13}p_1+F_{14}q_1)(F_{23}p_2+F_{24}q_2)
\end{pmatrix}.$$
Thus, by the identification \eqref{eq:nullline2},
\begin{align*}
[(e_0+e_3)(Fh)]&=[F_1h(\1+\k^{'})(F_2(h^{-1})^t)^t]\\
               &\cong\left(\frac{F_{11}p_1+F_{12}q_1}{F_{13}p_1+F_{14}q_1},
\frac{F_{21}p_2+F_{22}q_2}{F_{23}p_2+F_{24}q_2}\right)\in\E^2_1(u,v).
\end{align*}
Note that the hyperbolic Gau{\ss} map $[e_0+e_3]$ is orientation
preserving while $[e_0-e_3]$ is orientation reversing. So,
\begin{align*}
[(e_0-e_3)(Fh)]&=[F_1h(\1-\k^{'})(F_2(h^{-1})^t)^t]\\
               &\cong\left(\frac{F_{11}p_1-F_{12}q_1}{F_{13}p_1-F_{14}q_1},\frac{F_{21}p_2-F_{22}q_2}{F_{23}p_2-F_{24}q_2}\right)\in\E^2_1(u,v),
\end{align*}
where $h=\begin{pmatrix}
q_2 & -p_1\\
p_2 & q_1
\end{pmatrix}$. (See Remark \ref{rem:bryant}.)

Similarly, if $(e_0\pm e_3)(Fh)=F_1h(\1\pm\k^{'})(F_2h)^{-1}$, then
\begin{align*}
[(e_0+e_3)(Fh)]&\cong\left(\frac{F_{11}p_1+F_{12}q_1}{F_{13}p_1+F_{14}q_1},
-\frac{F_{24}p_2-F_{23}q_2}{F_{22}p_2-F_{21}q_2}\right)\in\E^2_1(u,v),\\
[(e_0-e_3)(Fh)]&\cong\left(\frac{F_{11}p_1-F_{12}q_1}{F_{13}p_1-F_{14}q_1},-\frac{F_{24}p_2+F_{23}q_2}{F_{22}p_2+F_{21}q_2}\right)\in\E^2_1(u,v).
\end{align*}

Let $q:=\dis\frac{p_1}{q_1}$ and $r:=\dis\frac{p_2}{q_2}$. The
ordered pair $(q,r)$ is called the \emph{secondary Gau{\ss}
  map}\footnote{As seen in section \ref{sec:tlminimal}, this secondary
Gau{\ss} map coincides with the projected Gau{\ss} map of a corresponding timelike minimal surface in $\E^3_1$.}. In
terms of the secondary Gau{\ss} map $(q,r)$, locally in an open
set $U\subset M$,
\begin{align*}
F_1^{-1}dF_1&=\begin{pmatrix}
q & -q^2\\
1 & -q
\end{pmatrix}f(u)du=\begin{pmatrix}
q & -q^2\\
1 & -q
\end{pmatrix}\eta,\\
F_2^{-1}dF_2&=\begin{pmatrix}
r & -r^2\\
1 & -r
\end{pmatrix}g(v)dv=\begin{pmatrix}
r & -r^2\\
1 & -r
\end{pmatrix}\xi,\\
(dF_2^{-1})F_2&=\begin{pmatrix}
r & 1\\
-r^2 & -r
\end{pmatrix}g(v)dv=\begin{pmatrix}
r & 1\\
-r^2 & -r
\end{pmatrix}\xi.
\end{align*}
Thus, we have the following equations:
\begin{align*}
dF_1&=\begin{pmatrix}
dF_{11} & dF_{12}\\
dF_{13} & dF_{14}
\end{pmatrix}=\begin{pmatrix}
F_{11} & F_{12}\\
F_{13} & F_{14}
\end{pmatrix}\begin{pmatrix}
q & -q^2\\
1 & -q
\end{pmatrix}\eta\\
    &=\begin{pmatrix}
F_{11}q+F_{12} & -(F_{11}q+F_{12}q)\\
F_{13}q+F_{14} & -(F_{13}q+F_{14}q)
\end{pmatrix}\eta,
\end{align*}
\begin{align*}
dF_2&=\begin{pmatrix}
dF_{21} & dF_{22}\\
dF_{23} & dF_{24}
\end{pmatrix}=\begin{pmatrix}
F_{21} & F_{22}\\
F_{23} & F_{24}
\end{pmatrix}\begin{pmatrix}
r & -r^2\\
1 & -r
\end{pmatrix}\xi\\
    &=\begin{pmatrix}
F_{21}r+F_{22} & -(F_{21}r+F_{22})r\\
F_{23}r+F_{24} & -(F_{23}r+F_{24})r
\end{pmatrix}\xi
\end{align*}
and
\begin{align*}
dF_2^{-1}&=\begin{pmatrix}
dF_{24} & -dF_{22}\\
-dF_{23} & dF_{21}
\end{pmatrix}=\begin{pmatrix}
r & 1\\
-r^2 & -r
\end{pmatrix}\begin{pmatrix}
F_{24} & -F_{22}\\
-F_{23} & F_{21}
\end{pmatrix}\xi\\
         &=\begin{pmatrix}
F_{24}r-F_{23} & -(F_{22}r-F_{21})\\
-(F_{24}r-F_{23})r & (F_{22}r-F_{21})r
\end{pmatrix}\xi,
\end{align*}
i.e.,
$dF_2$ is also given by
$$dF_2=\begin{pmatrix}
(F_{22}r-F_{21})r & F_{22}r-F_{21}\\
(F_{24}r-F_{23})r & F_{24}r-F_{23}
\end{pmatrix}\xi.$$
Hence, the hyperbolic Gau{\ss} map can be written:
\begin{align*}
[(e_0+e_3)(Fh)]&=[F_1h(\1+\k^{'})(F_2(h^{-1})^t)^t]\\
               &\cong\left(\frac{F_{11}q+F_{12}}{F_{13}q+F_{14}},
\frac{F_{21}r+F_{22}}{F_{23}r+F_{24}}\right)\\
               &=\left(\frac{dF_{11}}{dF_{13}},\frac{dF_{21}}{dF_{23}}\right)\\
               &=\left(\frac{dF_{12}}{dF_{14}},\frac{dF_{22}}{dF_{24}}\right)
\in\E^2_1(u,v),\\
[(e_0+e_3)(Fh)]&=[F_1h(\1+\k^{'})(F_2h)^{-1}]\\
               &=\left(\frac{F_{11}q+F_{12}}{F_{13}q+F_{14}},-\frac{F_{24}r-
F_{23}}{F_{22}r-F_{21}}\right)\\
               &=\left(\frac{dF_{11}}{dF_{13}},-\frac{dF_{23}}{dF_{21}}\right)
\\
               &=\left(\frac{dF_{12}}{dF_{14}},-\frac{dF_{24}}{dF_{22}}\right)
\in\E^2_1(u,v)
\end{align*}
and
\begin{align*}
[(e_0-e_3)(Fh)]&=[F_1h(\1-\k^{'})(F_2(h^{-1})^t)^t]\\
               &\cong\left(\frac{F_{11}q-F_{12}}{F_{13}q-F_{14}},\frac{F_{21}r-F_{22}}{F_{23}r-F_{24}}\right)\\
               &=\left(\frac{dF_{11}}{dF_{13}},\frac{dF_{21}}{dF_{23}}\right)\\
               &=\left(\frac{dF_{12}}{dF_{14}},\frac{dF_{22}}{dF_{24}}\right)\in\E^2_1(u,v),\\
[(e_0-e_3)(Fh)]&=[F_1h(\1-\k^{'})(F_2h)^{-1}]\\
               &\cong\left(\frac{F_{11}q-F_{12}}{F_{13}q-F_{14}},-\frac{F_{24}r+F_{23}}{F_{22}r+F_{21}}\right)\\
               &=\left(\frac{dF_{11}}{dF_{13}},-\frac{dF_{23}}{dF_{21}}\right)\\
               &=\left(\frac{dF_{12}}{dF_{14}},-\frac{dF_{24}}{dF_{22}}\right)\in\E^2_1(u,v).
\end{align*}
\section{The Generalized Gau{\ss} Map and the Hyperbolic Gau{\ss} Map}
\label{sec:gengauss}
Let $G(2,\E^4_2)$ be the Grassmannian manifold of oriented timelike
$2$-planes in $\E^4_2$. Let $\varphi: M\longrightarrow\H^3_1(-1)$ be
an oriented timelike surface in $\H^3_1(-1)$. At each point $p\in M$,
there is a (timelike) tangent plane to the surface $\varphi$:
$$\varphi_*({\rm T}_p
M)=[\varphi_x\wedge\varphi_y]_p=[(\varphi_x+\varphi_y)\wedge(-\varphi_x+
\varphi_y)]_p$$
spanned by a timelike vector $\varphi_x$ and a spacelike vector
$\varphi_y$ or equivalently, by two null vectors
$\varphi_x+\varphi_y$ and $-\varphi_x+\varphi_y$. Define a map
$$\G: M\longrightarrow G(2,\E^4_2); p\in
M\stackrel{\G}{\longmapsto}[\varphi_x\wedge
\varphi_y]_p=[(\varphi_x+\varphi_y)\wedge(-\varphi_x+\varphi_y)]_p.$$
This map is called the \emph{generalized Gauss map} of a timelike
surface $\varphi: M\longrightarrow\H^3_1(-1)$.

Let $\u$ be a null vector, i.e., $\langle\u,\u\rangle=0$ and $[\u]$ denote the
null line generated by $\u$. Let $\bQ^2_0:=\{[\u]\in\R P^3_2:
\langle\u,\u\rangle=0\}$. Then there is an embedding
$$\Gamma: G(2,\E^4_2)\longrightarrow\bQ^2_0\times\bQ^2_0;
[\v\wedge\w]\stackrel{\Gamma}{\longmapsto}([\v+\w],[-\v+\w]),$$
where $\v$ is a timelike vector and $\w$ is a spacelike vector. That
is, $G(2,\E^4_2)\cong{\rm Im}\Gamma\subset\bQ^2_0\times\bQ^2_0$. The
Lie group ${\rm SL}_2\R\times{\rm SL}_2\R$ acts on $G(2,\E^4_2)$
transitively via the actions:
\begin{align*}
&\mu: ({\rm SL}_2\R\times{\rm SL}_2\R)\times
G(2,\E^4_2)\longrightarrow G(2,\E^4_2);\\
&\mu(g,[v\wedge\w]):=[\mu(g,\v)\wedge\mu(g,\w)]=[(g_1\v
    g_2^t)\wedge(g_1\w g_2^t)]
\end{align*}
and
\begin{align*}
&\nu: ({\rm SL}_2\R\times{\rm SL}_2\R)\times
G(2,\E^4_2)\longrightarrow G(2,\E^4_2);\\
&\nu(g,[v\wedge\w]):=[\nu(g,\v)\wedge\nu(g,\w)]=[(g_1\v
    g_2^{-1})\wedge(g_1\w g_2^{-1})]
\end{align*}
for $g=(g_1,g_2)\in{\rm SL}_2\R\times{\rm
    SL}_2\R$ and $[\v,\w]\in G(2,\E^4_2)$.
Note that
\begin{align*}
\mu(g,[\v\wedge\w])&\cong([g_1(\v+\w)g_2^t],[g_1(-\v+\w)g_2^t])\in\bQ^2_0
\times\bQ^2_0,\\
\nu(g,[\v\wedge\w])&\cong([g_1(\v+\w)g_2^{-1}],[g_1(-\v+\w)g_2^{-1}])\in
\bQ^2_0\times\bQ^2_0.
\end{align*}
The isotropy subgroup of ${\rm SL}_2\R\times{\rm
    SL}_2\R$ with the actions $\mu$ and $\nu$ at
\begin{align*}
[e_1\wedge e_2]&=[(e_1+e_2)\wedge(-e_1+e_2)]\\
               &=\left[\begin{pmatrix}
0 & 1\\
0 & 0
\end{pmatrix}\wedge\begin{pmatrix}
0 & 0\\
1 & 0
\end{pmatrix}\right]\\
               &\cong\left(\left[\begin{pmatrix}
0 & 1\\
0 & 0
\end{pmatrix}\right],\left[\begin{pmatrix}
0 & 0\\
1 & 0
\end{pmatrix}\right]\right)
\end{align*}
is $\R_+\times\R_+$, where $\R_+:=\left\{\begin{pmatrix}
r & 0\\
0 & \frac{1}{r}
\end{pmatrix}: r\in\R\setminus\{0\}\right\}$.
Thus, the Grassmannian manifold $G(2,\E^4_2)$ can be represented as a
symmetric space
\begin{align*}
&\frac{{\rm SL}_2\R\times{\rm
  SL}_2\R}{\R_+\times\R_+}\cong\\
&\left\{\left(\left[g_1\begin{pmatrix}
0 & 1\\
0 & 0
\end{pmatrix}g_2^t\right],\left[g_1\begin{pmatrix}
0 & 0\\
1 & 0
\end{pmatrix}g_2^t\right]\right)\in\bQ^2_0\times\bQ^2_0: g_1,
g_2\in{\rm SL}_2\R\right\}
\end{align*} or
\begin{align*}
&\frac{{\rm SL}_2\R\times{\rm SL}_2\R}{\R_+\times\R_+}\cong\\
&\left\{\left(\left[g_1\begin{pmatrix}
0 & 1\\
0 & 0
\end{pmatrix}g_2^{-1}\right],\left[g_1\begin{pmatrix}
0 & 0\\
1 & 0
\end{pmatrix}g_2^{-1}\right]\right)\in\bQ^2_0\times\bQ^2_0: g_1,
g_2\in{\rm SL}_2\R\right\}.
\end{align*}
Denote by $G(2,\E^4_2)^-$ the Grassmannian manifold of \emph{negatively}
  oriented timelike $2$-planes in $\E^4_2$. Then
\begin{align*}
&G(2,\E^4_2)^-=\frac{{\rm SL}_2\R\times{\rm
  SL}_2\R}{\R_+\times\R_+}\cong\\
&\left\{\left(\left[g_1\begin{pmatrix}
0 & 0\\
1 & 0
\end{pmatrix}g_2^t\right],\left[g_1\begin{pmatrix}
0 & 1\\
0 & 0
\end{pmatrix}g_2^t\right]\right)\in\bQ^2_0\times\bQ^2_0: g_1,
g_2\in{\rm SL}_2\R\right\}
\end{align*}
or
\begin{align*}
&G(2,\E^4_2)^-=\frac{{\rm SL}_2\R\times{\rm
  SL}_2\R}{\R_+\times\R_+}\cong\\
&\left\{\left(\left[g_1\begin{pmatrix}
0 & 0\\
1 & 0
\end{pmatrix}g_2^{-1}\right],\left[g_1\begin{pmatrix}
0 & 1\\
0 & 0
\end{pmatrix}g_2^{-1}\right]\right)\in\bQ^2_0\times\bQ^2_0: g_1,
g_2\in{\rm SL}_2\R\right\}.
\end{align*}
Let us write $g_1=\begin{pmatrix}
g_{11} & g_{12}\\
g_{21} & g_{22}
\end{pmatrix}$ and $g_2=\begin{pmatrix}
g_{21} & g_{22}\\
g_{23} & g_{24}
\end{pmatrix}$.
Define a projection map
$$\phi=(\phi_1,\phi_2):
G(2,\E^4_2)\longrightarrow\E^2_1(u,v)\times\E^2_1(u,v)$$
by
$$\left(\left[g_1\begin{pmatrix}
0 & 1\\
0 & 0
\end{pmatrix}g_2^t\right],\left[g_1\begin{pmatrix}
0 & 0\\
1 & 0
\end{pmatrix}g_2^t\right]\right)\stackrel{(\phi_1,\phi_2)}{\longmapsto}\left(
\frac{g_{11}}{g_{13}},\frac{g_{21}}{g_{23}}\right)$$
or
$$\left(\left[g_1\begin{pmatrix}
0 & 1\\
0 & 0
\end{pmatrix}g_2^{-1}\right],\left[g_1\begin{pmatrix}
0 & 0\\
1 & 0
\end{pmatrix}g_2^{-1}\right]\right)\stackrel{(\phi_1,\phi_2)}{\longmapsto}
\left(\frac{g_{11}}{g_{13}},-\frac{g_{24}}{g_{22}}\right).$$
Similarly, we also define a projection map
$$\phi^-=(\phi^-_1,\phi^-_2):
G(2,\E^4_2)^-\longrightarrow\E^2_1(u,v)\times\E^2_1(u,v)$$
by
$$\left(\left[g_1\begin{pmatrix}
0 & 0\\
1 & 0
\end{pmatrix}g_2^t\right],\left[g_1\begin{pmatrix}
0 & 1\\
0 & 0
\end{pmatrix}g_2^t\right]\right)\stackrel{(\phi^-_1,\phi^-_2)}{\longmapsto}
\left(\frac{g_{12}}{g_{14}},\frac{g_{22}}{g_{24}}\right)$$
or
$$\left(\left[g_1\begin{pmatrix}
0 & 0\\
1 & 0
\end{pmatrix}g_2^{-1}\right],\left[g_1\begin{pmatrix}
0 & 1\\
0 & 0
\end{pmatrix}g_2^{-1}\right]\right)\stackrel{(\phi^-_1,\phi^-_2)}{\longmapsto}
\left(\frac{g_{12}}{g_{14}},-\frac{g_{23}}{g_{21}}\right).$$
Let $\varphi: M\longrightarrow\H^3_1(-1)$ be a timelike surface from
an oriented and simply-connected open $2$-manifold $M$ into $\H^3_1(-1)$ with
$ds^2_\varphi=e^\omega(-dx^2+dy^2)=e^\omega dudv$.
Then there exists an adapted framing $F: M\longrightarrow{\rm
  SL}_2\R\times{\rm SL}_2\R$ of $\varphi$ such that
$e_1\circ F=e^{-\omega}\varphi_x$ and $e_2\circ
F=e^{-\omega}\varphi_y$. The generalized Gau{\ss} map $\G$ of
$\varphi$ can be written
\begin{align*}
\G&=[(e_1\circ F)\wedge(e_2\circ F)]\\
  &=[(e_1+e_2)(F)\wedge(-e_1+e_2)(F)]\\
  &\cong\left(\left[F_1\begin{pmatrix}
0 & 1\\
0 & 0
\end{pmatrix}F_2^t\right],\left[F_1\begin{pmatrix}
0 & 0\\
1 & 0
\end{pmatrix}F_2^t\right]\right): M\longrightarrow G(2,\E^4_2).
\end{align*}
Let $\G_1=\phi_1\circ\G$ and $\G_2=\phi_2\circ\G$. Then
$\G_1=\dis\frac{F_{11}}{F_{13}},\ \G_2=\dis\frac{F_{21}}{F_{23}}$.
Note that the ordered pair
$(\G_1,\G_2)=\left(\dis\frac{F_{11}}{F_{13}},\frac{F_{21}}{F_{23}}\right)$
is the same as the hyperbolic Gau{\ss} map
$[(e_0+e_3)(F)]=[F_1(\1+\k^{'})F_2^t]$.

If
$\G\cong\left(\left[F_1\begin{pmatrix}
0 & 1\\
0 & 0
\end{pmatrix}F_2^{-1}\right],\left[F_1\begin{pmatrix}
0 & 0\\
1 & 0
\end{pmatrix}F_2^{-1}\right]\right)$, then
$$(\G_1,\G_2)=\left(\dis\frac{F_{11}}{F_{13}},-\frac{F_{24}}{F_{22}}\right)=[
F_1(\1+\k^{'})F_2^{-1}]=[(e_0+e_3)(F)].$$

Let us define $\G^-: M\longrightarrow
G(2,\E^4_2)^-$ by
\begin{align*}
\G^-&=[(-e_1+e_2)(F)\wedge(e_1+e_2)(F)]\\
  &\cong\left(\left[F_1\begin{pmatrix}
0 & 0\\
1 & 0
\end{pmatrix}F_2^t\right],\left[F_1\begin{pmatrix}
0 & 1\\
0 & 0
\end{pmatrix}F_2^t\right]\right).
\end{align*}
Let $\G^-_1:=\phi^-_1\circ \G^-,\
\G^-_2:=\phi^-_2\circ \G^-$
Then
$$(\G^-_1,\G^-_2)=\left(\frac{F_{12}}{F_{14}},
\frac{F_{22}}{F_{24}}\right)=[F_1(\1-\k^{'})F_2^t]=[(e_0-e_3)(F)].$$
If $\G^-\cong\left(\left[F_1\begin{pmatrix}
0 & 0\\
1 & 0
\end{pmatrix}F_2^{-1}\right],\left[F_1\begin{pmatrix}
0 & 1\\
0 & 0
\end{pmatrix}F_2^{-1}\right]\right)$, then
$$(\G^-_1,\G^-_2)=\left(\frac{F_{12}}{F_{14}},
-\frac{F_{23}}{F_{21}}\right)=[F_1(\1-\k^{'})F_2^{-1}]=[(e_0-e_3)(F)].$$
\section{The Lorentz Holomorphicity of Hyperbolic Gau{\ss} Map and
  Timelike $\cmc$ $\pm 1$ Surfaces in $\H^3_1(-1)$}
\label{sec:hol}
In this section, we study the relationship between
Lorentz holomorphicity of (projected) hyperbolic Gau{\ss} map and
timelike $\cmc$ $\pm 1$ surfaces in $\H^3_1(-1)$. Their relationship is
summarized as the following theorem. Here, we assume that
$\varphi:=\Phi_1\Phi_2^t: M\longrightarrow\H^3_1(-1)$ and
$\psi:=\Psi_1\Psi_2^{-1}: M\longrightarrow\H^3_1(-1)$ are timelike surfaces
in $\H^3_1(-1)$, where $\Phi:=(\Phi_1,\Phi_2): M\longrightarrow{\rm
  SL}_2\R\times{\rm SL}_2\R$ and $\Psi:=(\Psi_1,\Psi_2): M\longrightarrow{\rm
  SL}_2\R\times{\rm SL}_2\R$ are
solutions of Lax equations \eqref{eq:lax3} and \eqref{eq:lax4}, resp.,
in Theorem \ref{thm:sym}.
Let us regard $[(e_0+e_3)(\Phi)]=[\Phi_1(\1+\k^{'})\Phi_2^t]$
($[(e_0+e_3)(\Psi)]=[\Psi_1(\1+\k^{'})\Psi_2^{-1}]$) as the
projected hyperbolic Gau{\ss}
map $\left(\dis\frac{\Phi_{11}}{\Phi_{13}},\frac{\Phi_{21}}{\Phi_{23}}
\right)\in\E^2_1(u,v)$
($\left(\dis\frac{\Psi_{11}}{\Psi_{13}},-\frac{\Psi_{24}}{\Psi_{22}}\right)\in\E^2_1(u,v)$).
Also, regard
$[(e_0-e_3)(\Phi)]=[\Phi_1(\1-\k^{'})\Phi_2^t]$
($[(e_0-e_3)(\Psi)]=[\Psi_1(\1-\k^{'})\Psi_2^{-1}]$) as the projected
hyperbolic Gau{\ss} map
$\left(\dis\frac{\Phi_{12}}{\Phi_{14}},\frac{\Phi_{22}}
{\Phi_{24}}\right)\in\E^2_1(u,v)$
($\left(\dis\frac{\Psi_{12}}{\Psi_{14}},-\frac{\Psi_{23}}{\Psi_{21}}
\right)\in\E^2_1(u,v)$). Then we have the following theorem holds.
\begin{theorem}
\begin{enumerate}
\item $[(e_0+e_3)(\Phi)]$ $([(e_0+e_3)(\Psi))$ is
  Lorentz antiholomorphic if and only if
$\varphi$ $(\psi)$ satisfies $H=1$ and $Q=0$.
\item $[(e_0+e_3)(\Phi)]$ $([(e_0+e_3)(\Psi))$ is
  Lorentz holomorphic if and only if
$\varphi$ $(\psi)$ satisfies $H=1$ and $R=0$.
\item $[(e_0-e_3)(\Phi)]$ $([(e_0-e_3)(\Psi))$ is
  Lorentz antiholomorphic if and only if
$\varphi$ $(\psi)$ satisfies $H=-1$ and $Q=0$.
\item $[(e_0-e_3)(\Phi)]$ $([(e_0-e_3)(\Psi))$ is
  Lorentz holomorphic if and only if
$\varphi$ $(\psi)$ satisfies $H=-1$ and $R=0$.
\end{enumerate}
\end{theorem}
\begin{proof}
We prove only part (1). The rest can be proved similarly. Since
$$[(e_0+e_3)(\Phi)]=\left(\dis\frac{\Phi_{11}}{\Phi_{13}},\frac{\Phi_{21}}
{\Phi_{23}}\right),$$
$[(e_0+e_3)(\Phi)]$ is Lorentz antiholomorphic, i.e.,
$[(e_0+e_3)(\Phi)]_u=0$ if and only if
$(\Phi_{11})_u\Phi_{13}-\Phi_{11}(\Phi_{13})_u=0$ and
$(\Phi_{21})_u\Phi_{23}-\Phi_{21}(\Phi_{23})_u=0$.

On the other hand, from the Lax equations \eqref{eq:lax3},
\begin{align*}
\begin{pmatrix}
(\Phi_{11})_u & (\Phi_{12})_u\\
(\Phi_{13})_u & (\Phi_{14})_u
\end{pmatrix}&=\begin{pmatrix}
\Phi_{11} & \Phi_{12}\\
\Phi_{13} & \Phi_{14}
\end{pmatrix}\begin{pmatrix}
\frac{\omega_u}{4} & \frac{1}{2}e^{\frac{\omega}{2}}(H+1)\\
-e^{-\frac{\omega}{2}}Q & -\frac{\omega_u}{4}
\end{pmatrix}\\
             &=\begin{pmatrix}
\frac{1}{4}\Phi_{11}\omega_u-\Phi_{12}e^{\frac{-\omega}{2}}Q
& \frac{1}{2}\Phi_{11}e^{\frac{\omega}{2}}(H+1)-\frac{1}{4}\Phi_{12}\omega_u\\
\frac{1}{4}\Phi_{13}\omega_u-\Phi_{14}e^{\frac{-\omega}{2}}Q
& \frac{1}{2}\Phi_{13}e^{\frac{\omega}{2}}(H+1)-\frac{1}{4}\Phi_{14}\omega_u
\end{pmatrix}
\end{align*}
and
\begin{align*}
\begin{pmatrix}
(\Phi_{21})_u & (\Phi_{22})_u\\
(\Phi_{23})_u & (\Phi_{24})_u
\end{pmatrix}&=\begin{pmatrix}
\Phi_{21} & \Phi_{22}\\
\Phi_{23} & \Phi_{24}
\end{pmatrix}\begin{pmatrix}
-\frac{\omega_u}{4} & e^{-\frac{\omega}{2}}Q\\
-\frac{1}{2}e^{\frac{\omega}{2}}(H-1) & \frac{\omega_u}{4}
\end{pmatrix}\\
             &=\begin{pmatrix}
-\frac{1}{4}\Phi_{21}\omega_u-\frac{1}{2}\Phi_{22}e^{\frac{\omega}{2}}(H-1) &
\Phi_{21}e^{-\frac{\omega}{2}}Q+\frac{1}{4}\Phi_{22}\omega_u\\
-\frac{1}{4}\Phi_{23}\omega_u-\frac{1}{2}\Phi_{24}e^{\frac{\omega}{2}}(H-1) &
\Phi_{23}e^{-\frac{\omega}{2}}Q+\frac{1}{4}\Phi_{24}\omega_u
\end{pmatrix}.
\end{align*}
Thus,
\begin{align*}
(\Phi_{11})_u\Phi_{13}-\Phi_{11}(\Phi_{13})_u&=e^{-\frac{\omega}{2}}Q,\\
(\Phi_{21})_u\Phi_{23}-\Phi_{21}(\Phi_{23})_u&=\frac{1}{2}e^{\frac{\omega}{2}}(H-1).
\end{align*}
It then follows immediately that $[(e_0+e_3)(\Phi)]$ is
Lorentz antiholomorphic if and only if $\varphi$ satisfies $H=1$ and
$Q=0$.
\end{proof}
\begin{corollary}
\begin{enumerate}
\item $[(e_0+e_3)(\Phi)]$ $([(e_0+e_3)(\Psi))$ is
  constant if and only if $\varphi$ $(\psi)$ satisfies $H=1$ and is
  totally umbilic ($\Q=0$, i.e., $Q=R=0$).
\item $[(e_0-e_3)(\Phi)]$ $([(e_0-e_3)(\Psi))$ is
 constant if and only if
$\varphi$ $(\psi)$ satisfies $H=-1$ and is totally umbilic.
\end{enumerate}
\end{corollary}
\section{Appendix I: The Lawson-Guichard Correspondence between Timelike $\cmc$
 Surfaces in Different Semi-Riemannian Space Forms}
\label{sec:lawson}
In section \ref{sec:intsystem}, we discussed the \emph{Lawson-Guichard
correspondence} or simply \emph{Lawson correspondence} between
timelike $\cmc$ surfaces in semi-Riemannian space forms $\E^3_1$,
$\S^3_1(1)$ and $\H^3_1(-1)$. In fact, this Lawson correspondence was
already known to A.~Fujioka and J.~Inoguchi (\cite{F-I} and
\cite{F-I2}). In this appendix, we study the Lawson correspondence in
a more general setting.

Let $\bar M$ be a semi-Riemannian manifold and $M\subset\bar M$ a
hypersurface with the sectional curvatures $\bar K$ and $K$, resp. Let
${\mathcal S}$ be the shape operator derived from the unit normal vector field
$N$ on the hypersurface $M$. If $X,Y$ span a
nondegenerate tangent plane on $M$, then the Gau{\ss} equation is given
by
\begin{equation}
\label{eqn:gausseq}
K(X,Y)=\bar K(X,Y)+\epsilon\frac{\langle{\mathcal S}(X),X\rangle\langle
  {\mathcal S}(Y),Y\rangle-\langle{\mathcal S}(X),Y\rangle^2}{\langle X,X\rangle\langle
  Y,Y\rangle-\langle X,Y\rangle^2},
\end{equation}
where $\epsilon=\langle N,N\rangle$. (See B. O'Neill \cite{Oneill} on p.~107.)
We begin with the following theorem which can be found, for example,
in T.~Weinstein \cite{We} on p.~158.
\begin{theorem}[Fundamental Theorem of Surface Theory: Lorentzian
    Version]
Given a simply-connected Lorentz surface with global null coordinate
system, there exists a timelike immersion with the first and the
second fundamental forms $I$ and $I\!I$ if and only if $I$ and $I\!I$
satisfy the Gau{\ss} and Mainardi-Codazzi equations.
\end{theorem}

Let $\M^3(\bar K)$ be the semi-Riemannian $3$-manifold with constant
sectional curvature $\bar K$. For example, $\M^3(-1)=\H^3_1(-1)$,
$\M^3(0)=\E^3_1$, and $\M^3(1)=\S^3_1(1)$. For a conformal timelike
immersion $\varphi: M\longrightarrow\M^3(\bar K)$ with induce metric
$\langle\cdot, \cdot\rangle$, Levi-Civita connection $\nabla$,
Gau{\ss}ian curvature $K$ and shape operator ${\mathcal S}$, the
Gau{\ss} and Mainardi-Codazzi\footnote{Since $\M^3(\bar K)$ has a
constant
  sectional curvature and $M$ is a hypersurface immersed into
  $\M^3(\bar K)$, the Mainardi-Codazzi equation becomes
  \eqref{eqn:codazzi}. See B. O'Neill \cite{Oneill} on p.~115
  for more details.} equations are satisfied:
\begin{align}
\label{eqn:gauss}
K-\bar K&=\det{\mathcal S}\ \mbox{\rm (Gau{\ss} Equation)}\\
\label{eqn:codazzi}
(\nabla_X{\mathcal S})(Y)&=(\nabla_Y{\mathcal S})(X)\ \mbox{\rm
  (Mainardi-Codazzi Equation)}
\end{align}
for all smooth vector fields $X, Y, Z\in T(M)$.

Assume that $\varphi$ has a constant mean curvature
$H=\frac{1}{2}\tr({\mathcal S})$. For any $c\in\R$, define
$$\tilde{\mathcal S}=S+c{\mathcal I},\ \tilde{K}=\bar K-c\tr({\mathcal
  S})-c^2,$$
where ${\mathcal I}$ is the identity transformation. Then the Gau{\ss}
  equation \eqref{eqn:gauss} and the Mainardi-Codazzi equation
  \eqref{eqn:codazzi} still hold when ${\mathcal S}$ and $\bar K$ are
  replaced by $\tilde{\mathcal S}$ and $\tilde{K}$, resp.:
\begin{align*}
K-\tilde{K}&=K-(\bar K-c\tr({\mathcal S})-c^2)\\
           &=K-\bar K+c\tr({\mathcal S})+c^2\\
           &=\det({\mathcal S})+c\tr({\mathcal S})+c^2\\
           &=\det({\mathcal S}+c{\mathcal I})\\
           &=\det(\tilde{\mathcal S}).
\end{align*}
Note that the Gau{\ss}ian curvature $K$ of $M$ is intrinsic and does
not change.
\begin{align*}
(\nabla_X{\mathcal S})(Y)&=\nabla_X{\mathcal S}(Y)-{\mathcal
    S}(\nabla_XY),\\
(\nabla_Y{\mathcal S})(X)&=\nabla_Y{\mathcal S}(X)-{\mathcal
    S}(\nabla_YX).
\end{align*}
Since $(\nabla_X{\mathcal S})(Y)=(\nabla_Y{\mathcal S})(X)$,
\begin{align*}
{\mathcal S}([X,Y])&={\mathcal S}(\nabla_XY-\nabla_YX)\\
                   &={\mathcal S}(\nabla_XY)-{\mathcal S}(\nabla_YX)\\
                   &=\nabla_X{\mathcal S}(Y)-\nabla_Y{\mathcal S}(X).
\end{align*}
Now,
\begin{align*}
\tilde{\mathcal S}&={\mathcal S}([X,Y])+c[X,Y]\\
                  &=\nabla_XY{\mathcal S}(Y)-\nabla_Y{\mathcal
                  S}(X)+c[X,Y]\\
                  &=\nabla_X({\mathcal S}+c{\mathcal
                  I})(Y)-\nabla_Y({\mathcal S}+c{\mathcal I})(X)\\
                  &=\nabla_X\tilde{\mathcal
                  S}(Y)-\nabla_Y\tilde{\mathcal S}(X).
\end{align*}
Therefore, there exists an immersion $\tilde{\varphi}:
M\longrightarrow\M^3(\tilde{K})$ with induced metric $\langle\cdot,
\cdot\rangle$ and shape operator $\tilde{\mathcal S}$, and
$\tilde{\varphi}(M)$ is isometric to $\varphi(M)$. The mean curvature
$\tilde{H}$ of $\tilde{\varphi}(M)$ is
$$\tilde{H}=\frac{1}{2}\tr(\tilde{\mathcal
  S})=\frac{1}{2}\tr({\mathcal S})+c=H+c$$
and this shows the Lawson-Guichard correspondence between timelike
$\cmc$ $H$ surfaces in $\M^3(\bar K)$ and timelike $\cmc$ $(H+c)$
surfaces in $\M^3(\tilde{K})=\M^3(\bar K-2cH-c^2)$. In particular,
  when $H=\bar K=0$ and $c=1$, we have the Lawson-Guichard
  correspondence between timelike minimal surfaces in $\E^3_1$ and
  timelike $\cmc$ $1$ surfaces in $\H^3_1(-1)$. If $H=0$ and $\bar
  K=c=1$, then we have the Lawson-Guichard correspondence between
  timelike minimal surfaces in de Sitter $3$-space $\S^3_1(1)$ and
  timelike $\cmc$ $1$ surfaces in $\E^3_1$.
\section{Appendix II: Some Examples of timelike $\cmc$ $\pm 1$ surfaces
  in $\H^3_1(-1)$}

In this appendix, we present some examples of timelike $\cmc$ $1$ surfaces
in $\H^3_1(-1)$.

Let us consider the following stereographic projections in order to
view the isometric images of timelike $\cmc$ $\pm 1$ surfaces in
$\H^3_1(-1)$ into the interior ${\rm Int}\S^2_1(1)=\{(x_1,x_2,x_3)\in\E^3_1:
-(x_1)^2+(x_2)^2+(x_3)^2<1\}$ of de Sitter
$2$-space $\S^2_1(1)$.

Let $\wp_+:
\H^3_1(-1)\setminus\{x_0=-1\}\longrightarrow\E^3_1\setminus\S^2_1(1)$
be the stereographic projection from $-e_0=(-1,0,0,0)$. Then
\begin{equation}
\label{eqn:stproj1}
\wp_+(x_0,x_1,x_2,x_3)=\left(\frac{x_1}{1+x_0},\frac{x_2}{1+x_0},\frac{x_3}{1+x_3}\right).
\end{equation}
Let $\wp_-:
\H^3_1(-1)\setminus\{x_0=1\}\longrightarrow\E^3_1\setminus\S^2_1(1)$
be the stereographic projection from $e_0=(1,0,0,0)$. Then
\begin{equation}
\label{eqn:stproj2}
\wp_-(x_0,x_1,x_2,x_3)=\left(\frac{x_1}{1-x_0},\frac{x_2}{1-x_0},\frac{x_3}{1-x_3}\right).
\end{equation}
Cut $\H^3_1(-1)$ into two halves by the hyperplane
$x_0=0$. Denote by $\H^3_1(-1)_+$ ($\H^3_1(-1)_-$) the half containing
$e_0=(1,0,0,0)$ ($-e_0=(-1,0,0,0)$). Then $\wp_+:
\H^3_1(-1)_+\longrightarrow{\rm Int}\S^2_1(1)$ and $\wp_-:
\H^3_1(-1)_-\longrightarrow{\rm Int}\S^2_1(1)$.

\begin{example}[Timelike Enneper Cousin in $\H^3_1(-1)$ (Isothermic
    Type)]
Let $(q,r)=(u,v)$. Then using the Bryant-Umehara-Yamada type
representation \eqref{eq:BUY}, we set up the following initial value
problem:
$$
F_1^{-1}dF_1=\begin{pmatrix}
u & -u^2\\
1 & -u
\end{pmatrix}du,\ F_2^{-1}dF_2=\begin{pmatrix}
v & -v^2\\
1 & -v
\end{pmatrix}dv$$
with the initial condition $F_1(0)=F_2(0)=\1$.
This initial value problem has a unique solution
\begin{align*}
F_1(u,v)&=\begin{pmatrix}
\cosh u & \sinh u-u\cosh u\\
\sinh u & \cosh u-u\sinh u
\end{pmatrix},\\
F_2(u,v)&=\begin{pmatrix}
\cosh v & \sinh v-v\cosh v\\
\sinh v & \cosh v-v\sinh v
\end{pmatrix}
\end{align*}
which are Lorentz holomorphic and Lorentz antiholomorphic null
immersions into ${\rm SL}_2\R$. The Bryant type representation formula
\eqref{eq:bryant} then yields a timelike $\cmc$ $1$ surface in
$\H^3_1(-1)$. The resulting surface is a correspondent of isothermic
type\footnote{For details about isothermic and anti-isothermic
      timelike surfaces, please see \cite{F-I2} or \cite{I-T}.}
timelike Enneper surface in $\E^3_1$ under the Lawson-Guichard
correspondence. For this reason, the resulting surface is called
\emph{isothermic type timelike Enneper cousin} in $\H^3_1(-1)$.

Figure \ref{fig:tlenneper} shows different views of isothermic
type timelike Enneper cousin in $\H^3_1(-1)$ projected via $\wp_+$
into the interior of the boundary $\S^2_1(1)$.

\begin{figure}[ht]
\centering
\mbox{\subfigure[]{\epsfig{figure=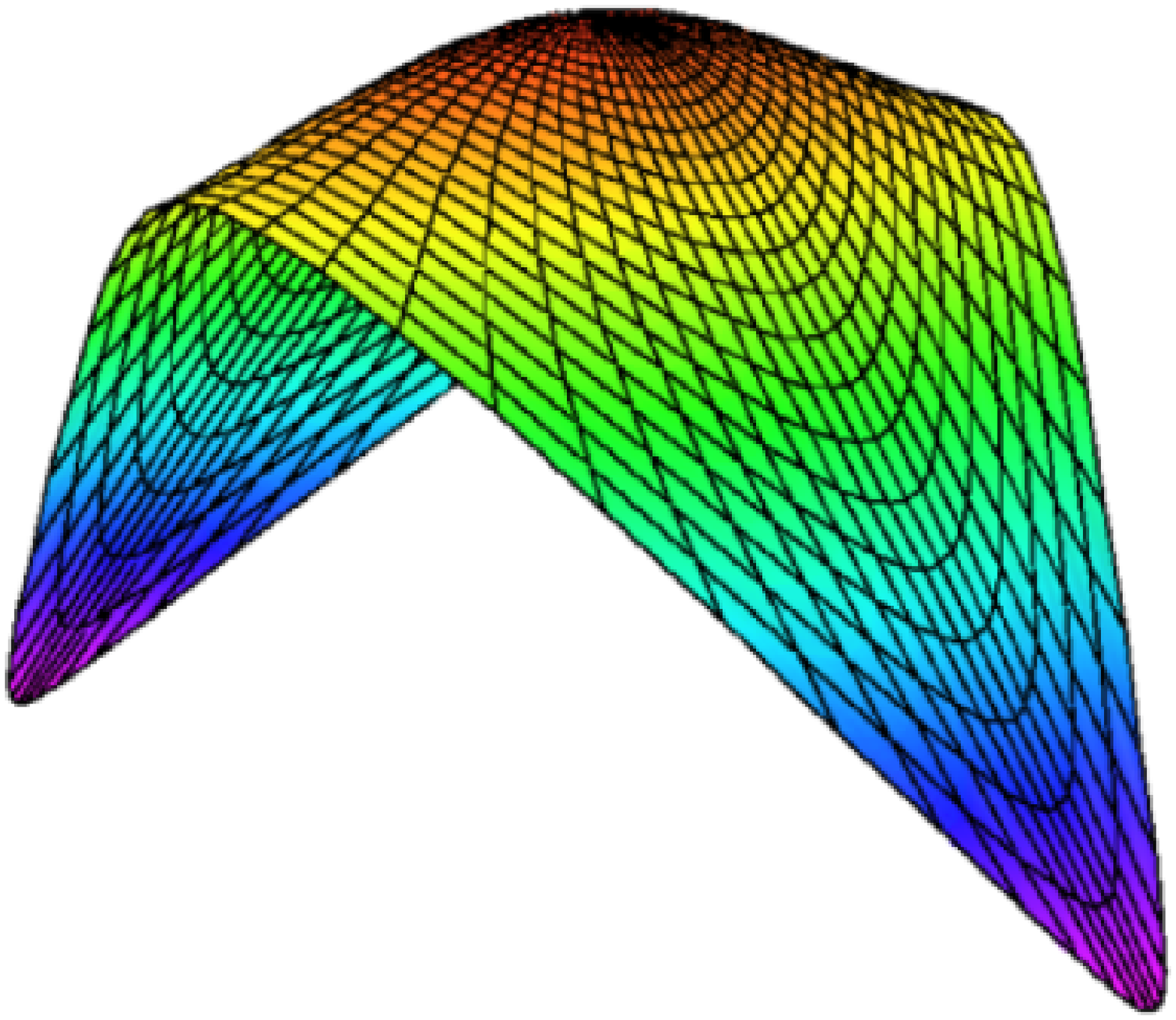,width=.48\textwidth}}\quad
    \subfigure[]{\epsfig{figure=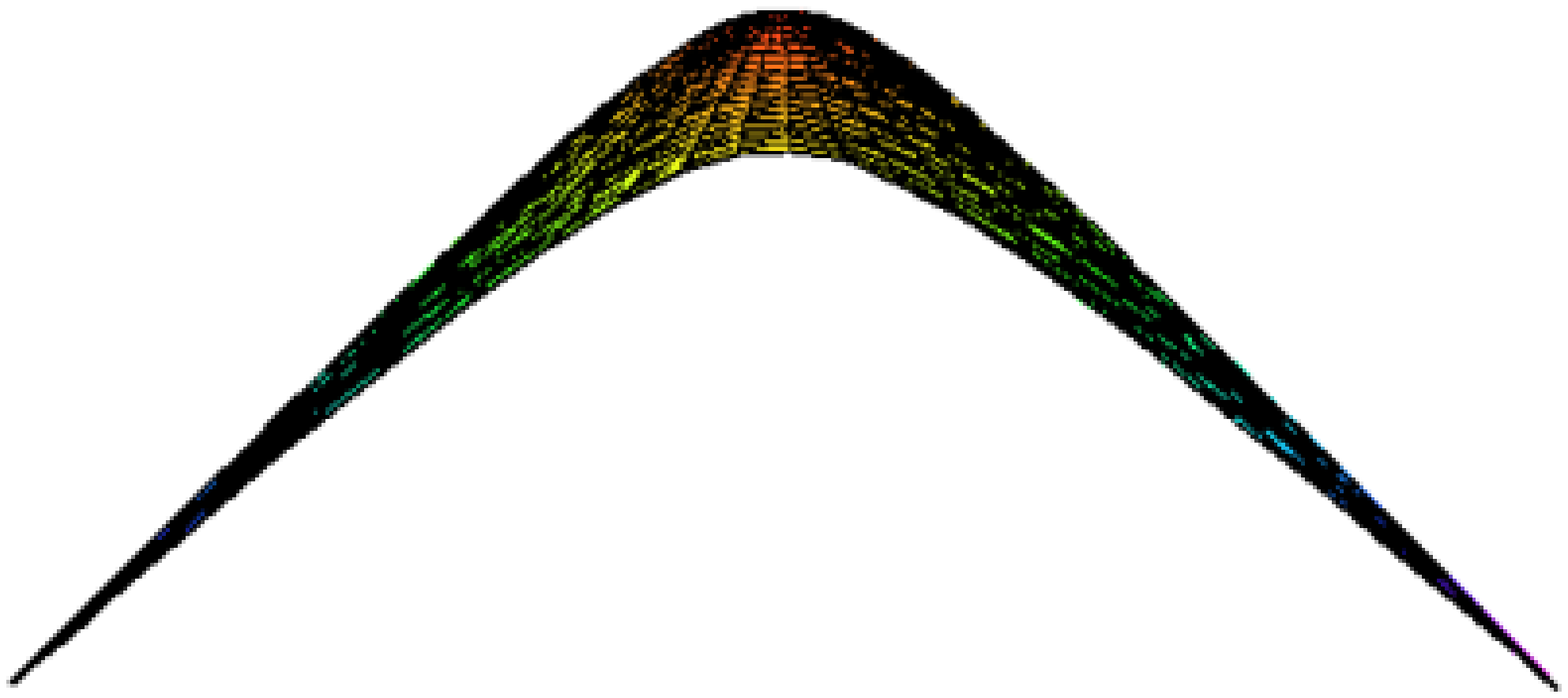,width=.48\textwidth}}}
\mbox{\subfigure[]{\epsfig{figure=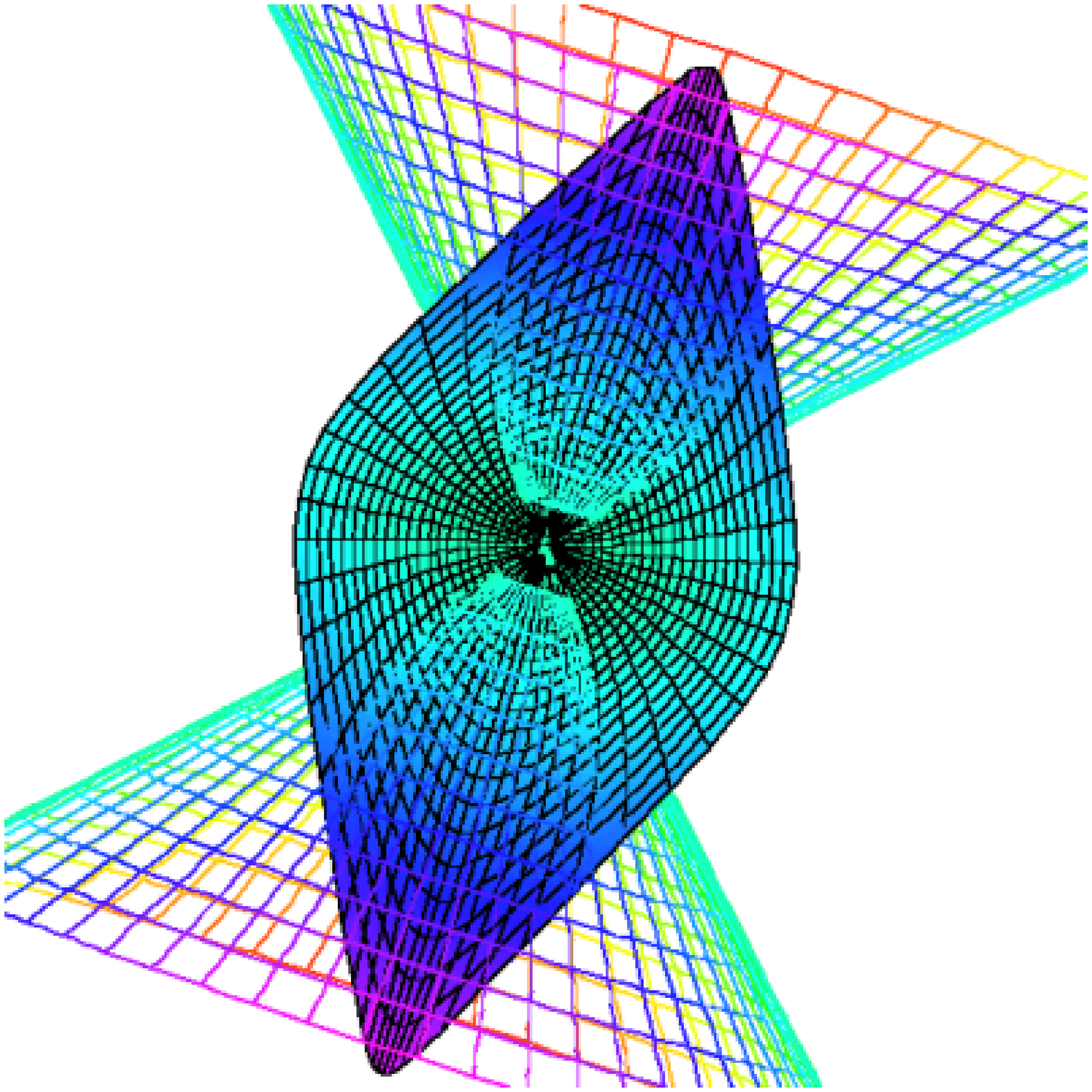,width=.48\textwidth}}\quad
    \subfigure[]{\epsfig{figure=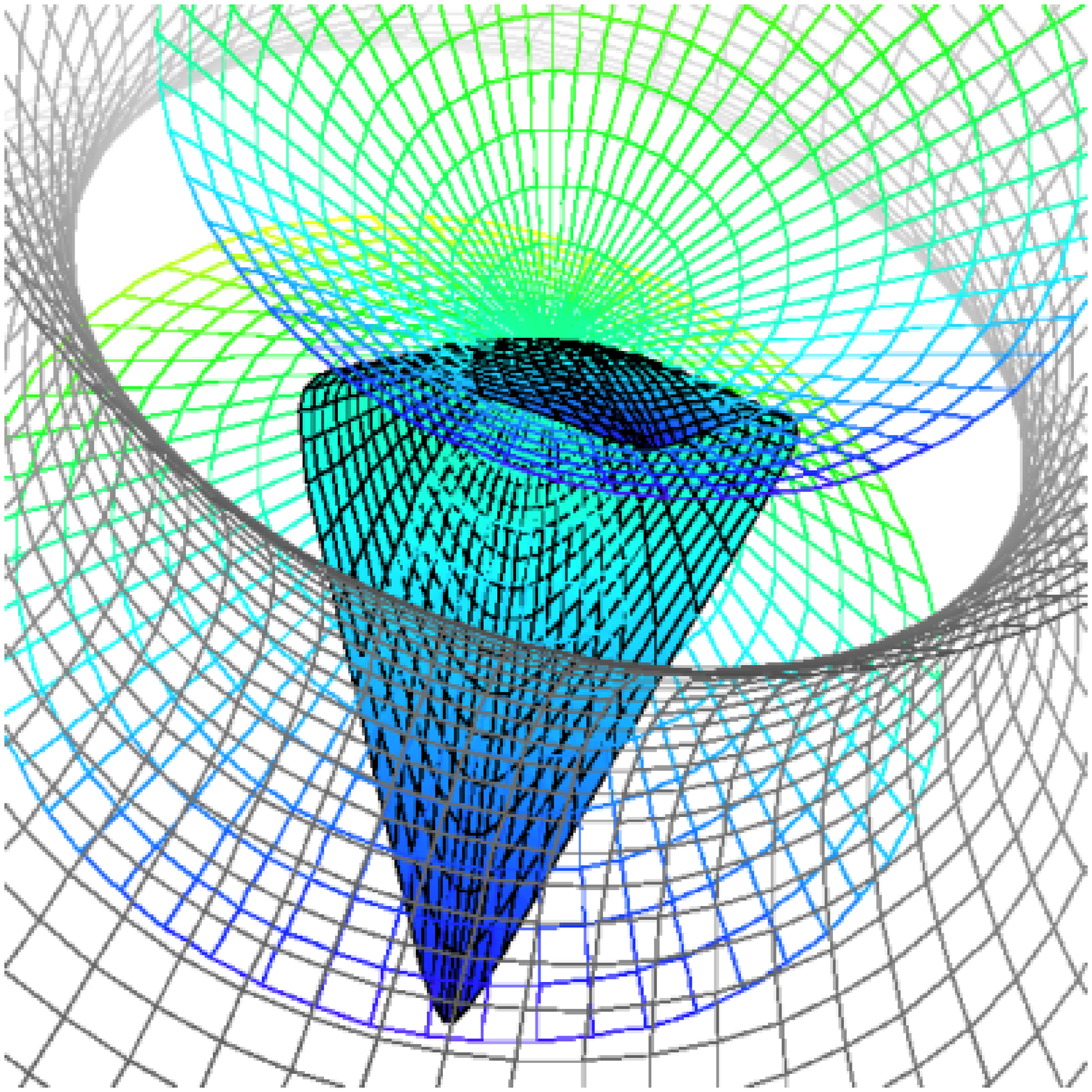,width=.48\textwidth}}}
\caption{Isothermic type timelike Enneper cousin projected into ${\rm
  Int}\S^2_1(1)$ via
  $\wp_+$ with light cone and the boundary $\S^2_1(1)$ in
  $\E^3_1$\label{fig:tlenneper}}
\end{figure}
Figure \ref{fig:tlminenneper} shows different views of isothermic type
timelike Enneper surface in $\E^3_1$.

\begin{figure}[ht]
\centering
\mbox{\subfigure[]{\epsfig{figure=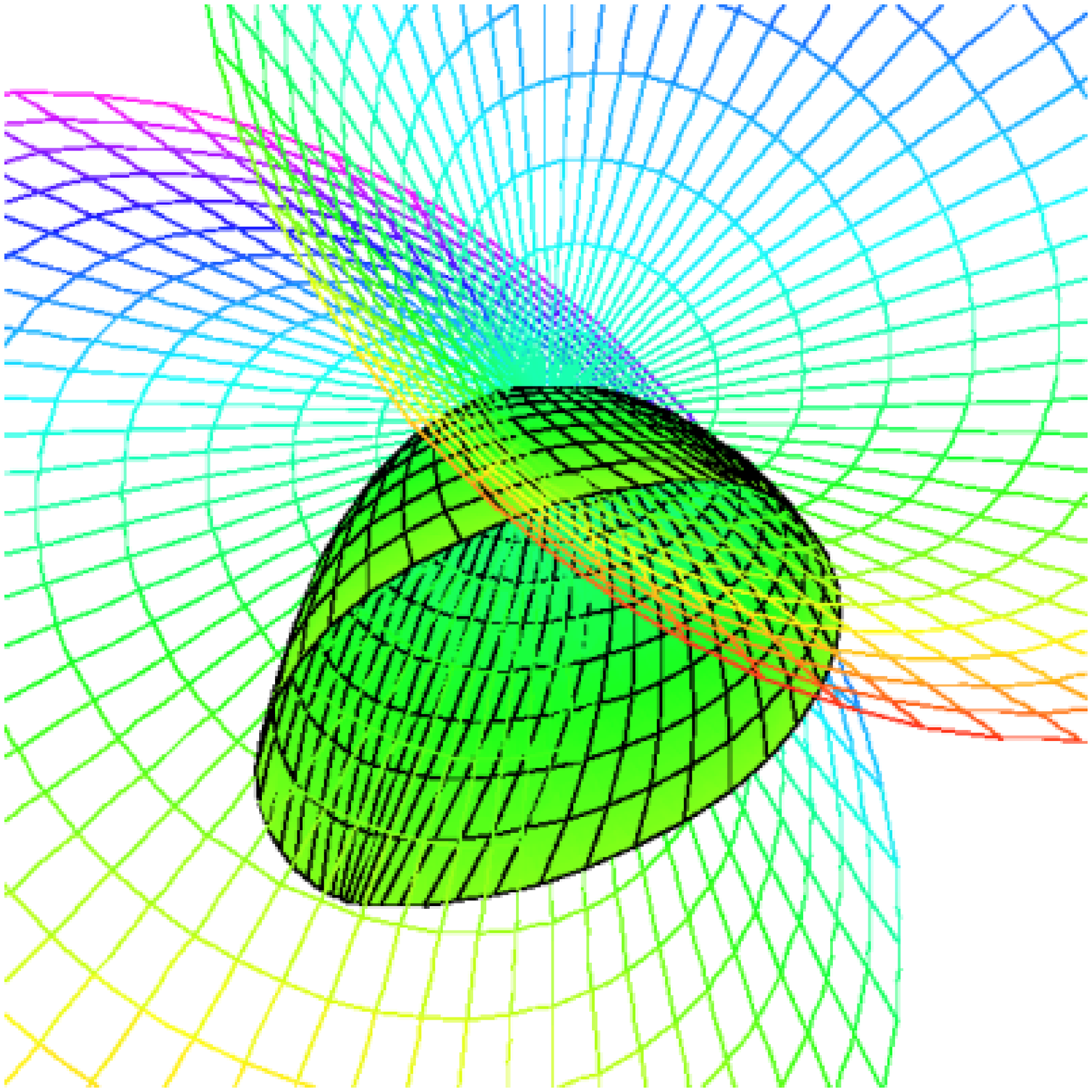,width=.48\textwidth}}\quad
    \subfigure[]{\epsfig{figure=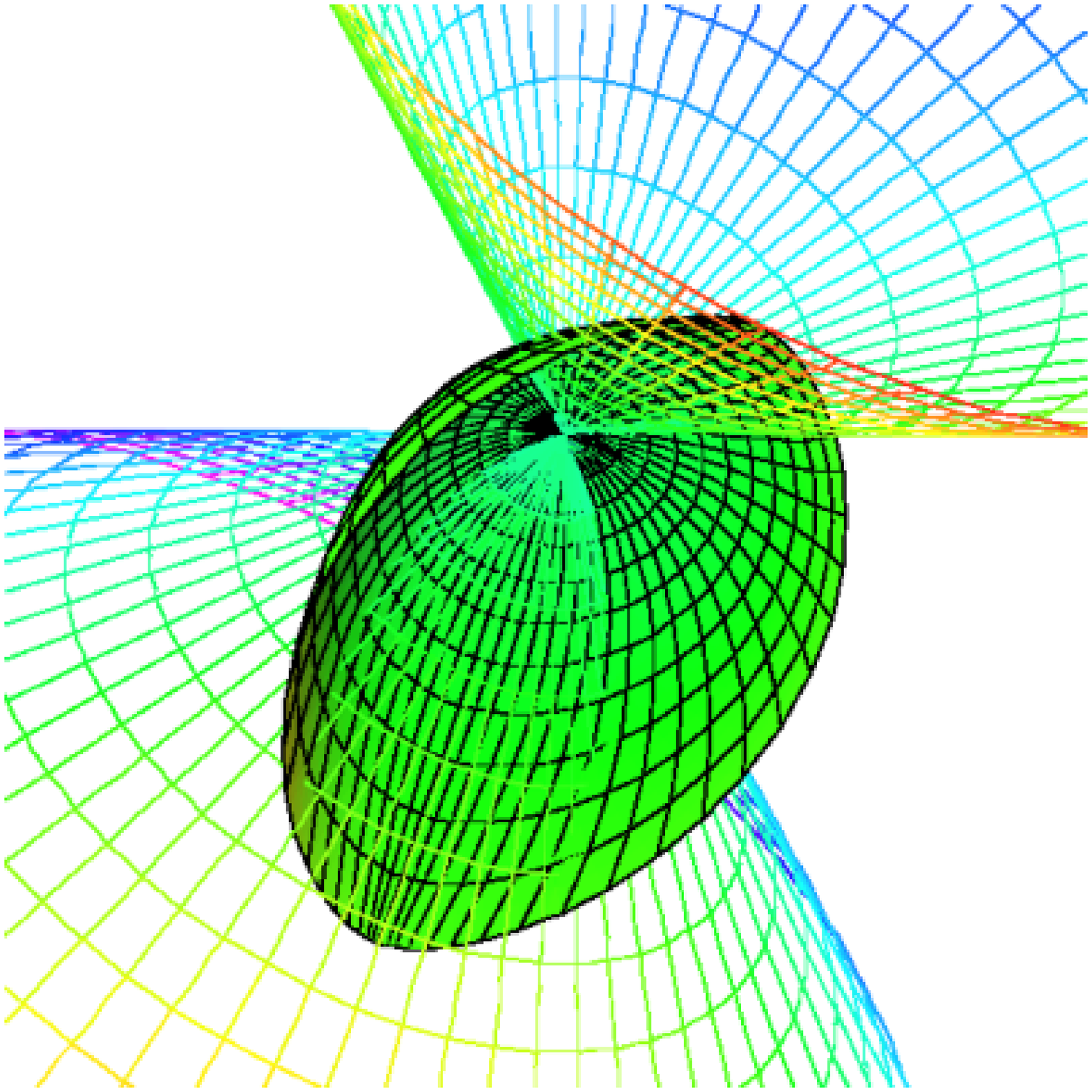,width=.48\textwidth}}}
\mbox{\subfigure[]{\epsfig{figure=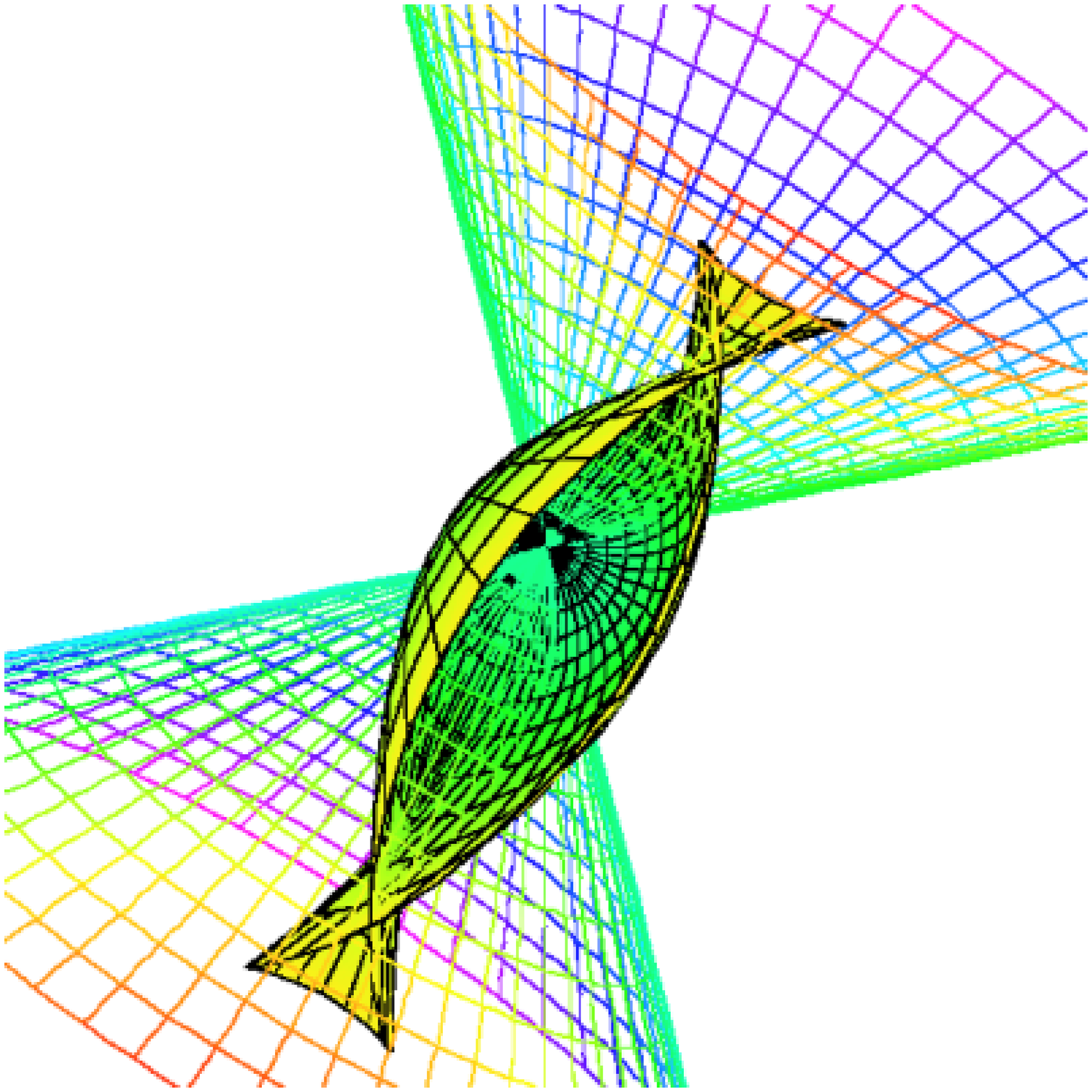,width=.48\textwidth}}\quad
    \subfigure[]{\epsfig{figure=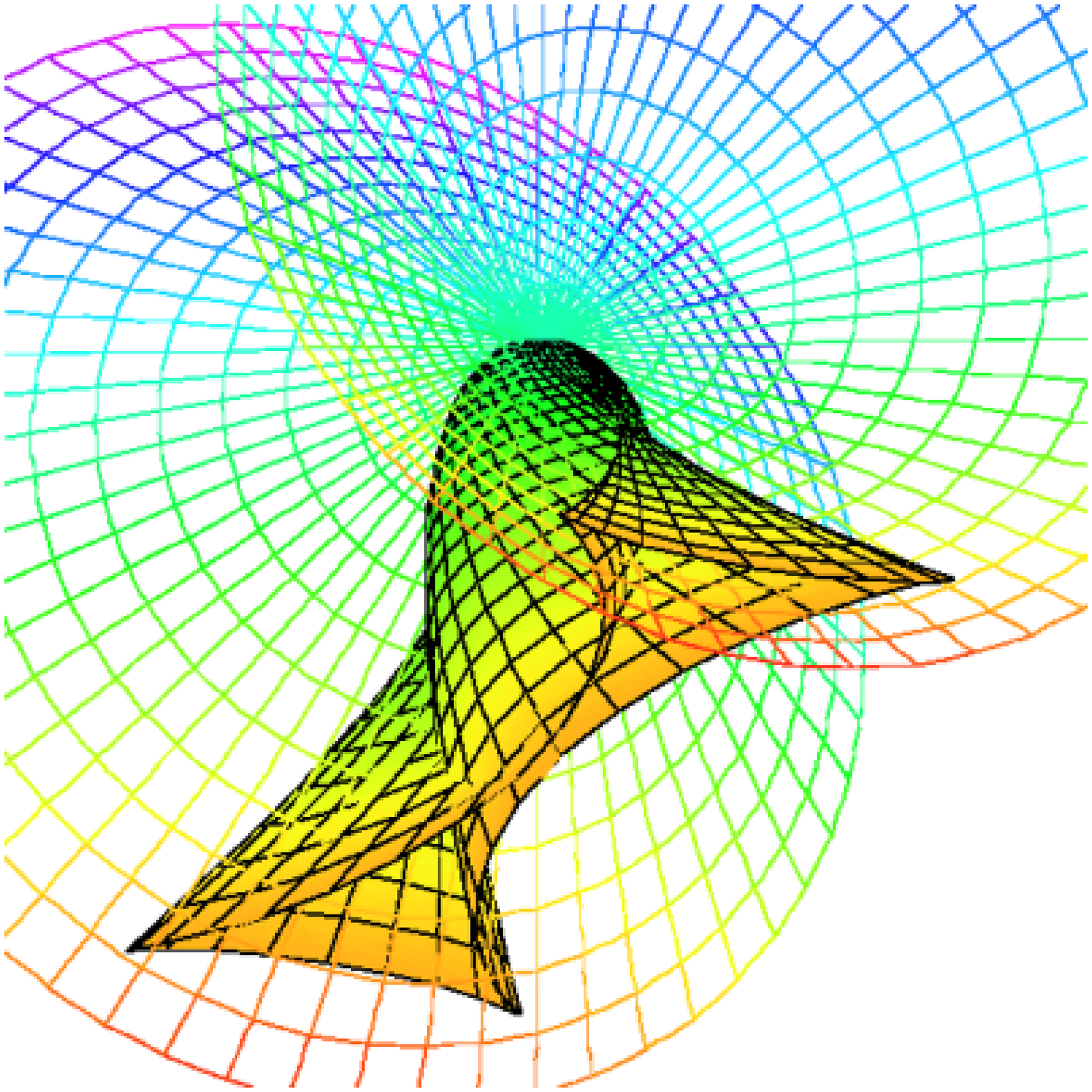,width=.48\textwidth}}}
\caption{Isothermic type timelike Enneper surface in
  $\E^3_1$ with light cone\label{fig:tlminenneper}}
\end{figure}
\end{example}
\begin{example}[Timelike Enneper Cousin in $\H^3_1(-1)$ (Anti-isothermic
    Type)]
Let $(q,r)=(-u,v)$. Then using the Bryant-Umehara-Yamada type
representation \eqref{eq:BUY}, we set up the following initial value
problem:
$$
F_1^{-1}dF_1=\begin{pmatrix}
-u & -u^2\\
1 & u
\end{pmatrix}du,\ F_2^{-1}dF_2=\begin{pmatrix}
v & -v^2\\
1 & -v
\end{pmatrix}dv$$
with the initial condition $F_1(0)=F_2(0)=\1$.
This initial value problem has a unique solution
\begin{align*}
F_1(u,v)&=\begin{pmatrix}
\cos u & -\sin u+u\cos u\\
\sin u & \cos u+u\sin u
\end{pmatrix},\\
F_2(u,v)&=\begin{pmatrix}
\cosh v & \sinh v-v\cosh v\\
\sinh v & \cosh v-v\sinh v
\end{pmatrix}
\end{align*}
which are Lorentz holomorphic and Lorentz antiholomorphic null
immersions into ${\rm SL}_2\R$. The Bryant type representation formula
\eqref{eq:bryant} then yields a timelike $\cmc$ $1$ surface in
$\H^3_1(-1)$. The resulting surface is a correspondent of anti-isothermic
type timelike Enneper surface in $\E^3_1$ under the Lawson-Guichard
correspondence. For this reason, the resulting surface is called
\emph{anti-isothermic type timelike Enneper cousin} in $\H^3_1(-1)$.

Figure \ref{fig:tlenneper2} shows different views of anti-isothermic
type timelike Enneper cousin in $\H^3_1(-1)$ projected via $\wp_+$
into the interior of the boundary $\S^2_1(1)$.

\begin{figure}[ht]
\centering
\mbox{\subfigure[]{\epsfig{figure=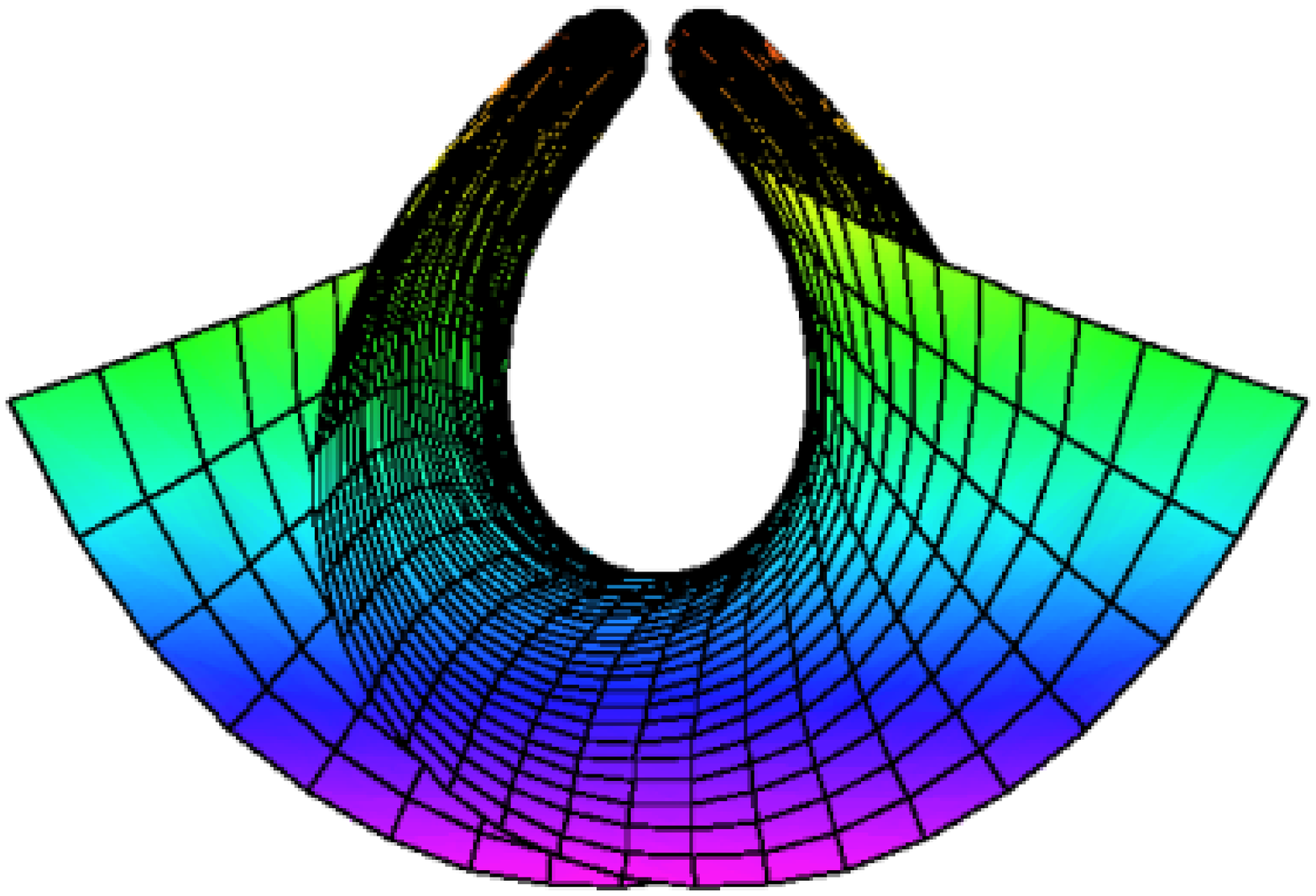,width=.48\textwidth}}\quad
    \subfigure[]{\epsfig{figure=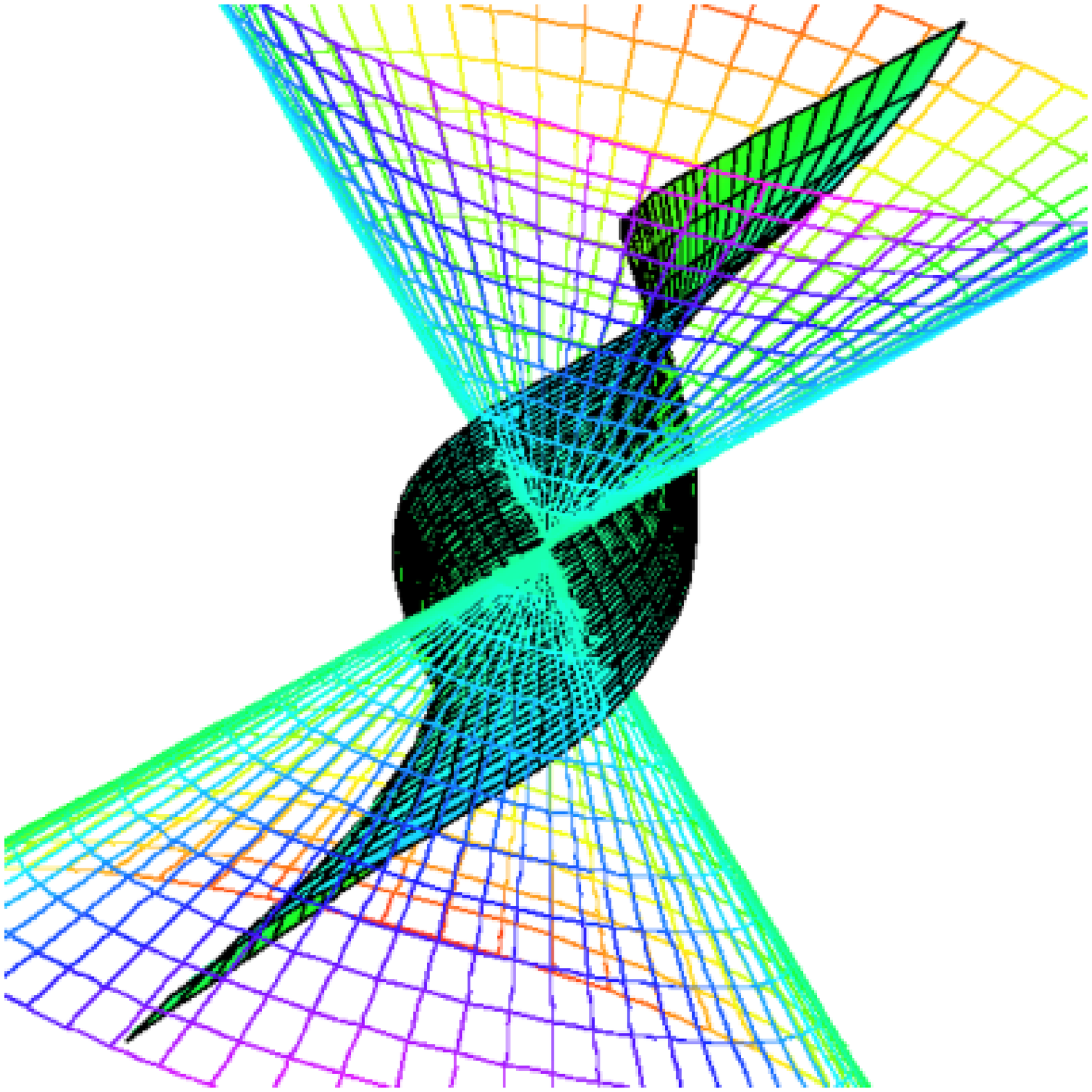,width=.48\textwidth}}}
\mbox{\subfigure[]{\epsfig{figure=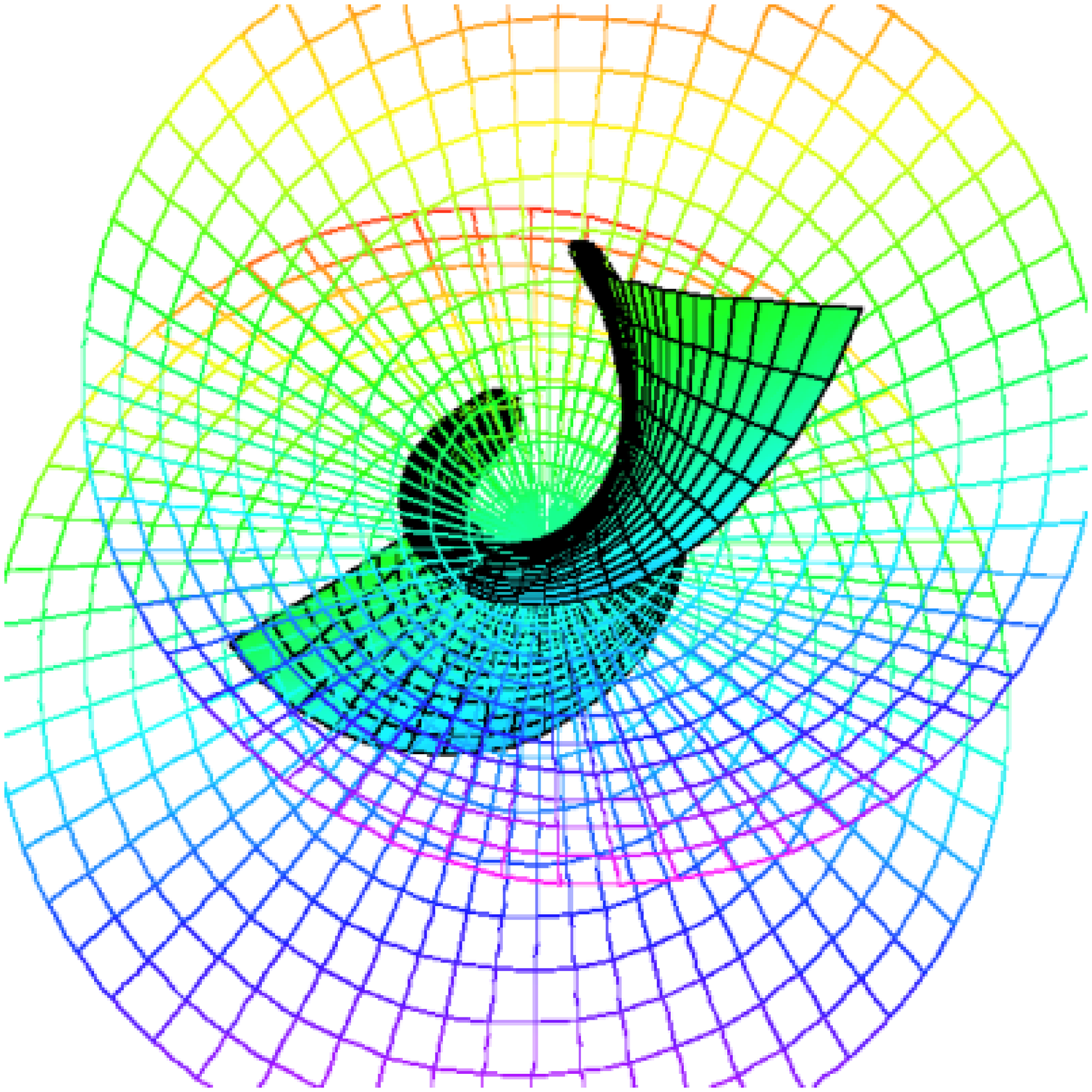,width=.48\textwidth}}\quad
    \subfigure[]{\epsfig{figure=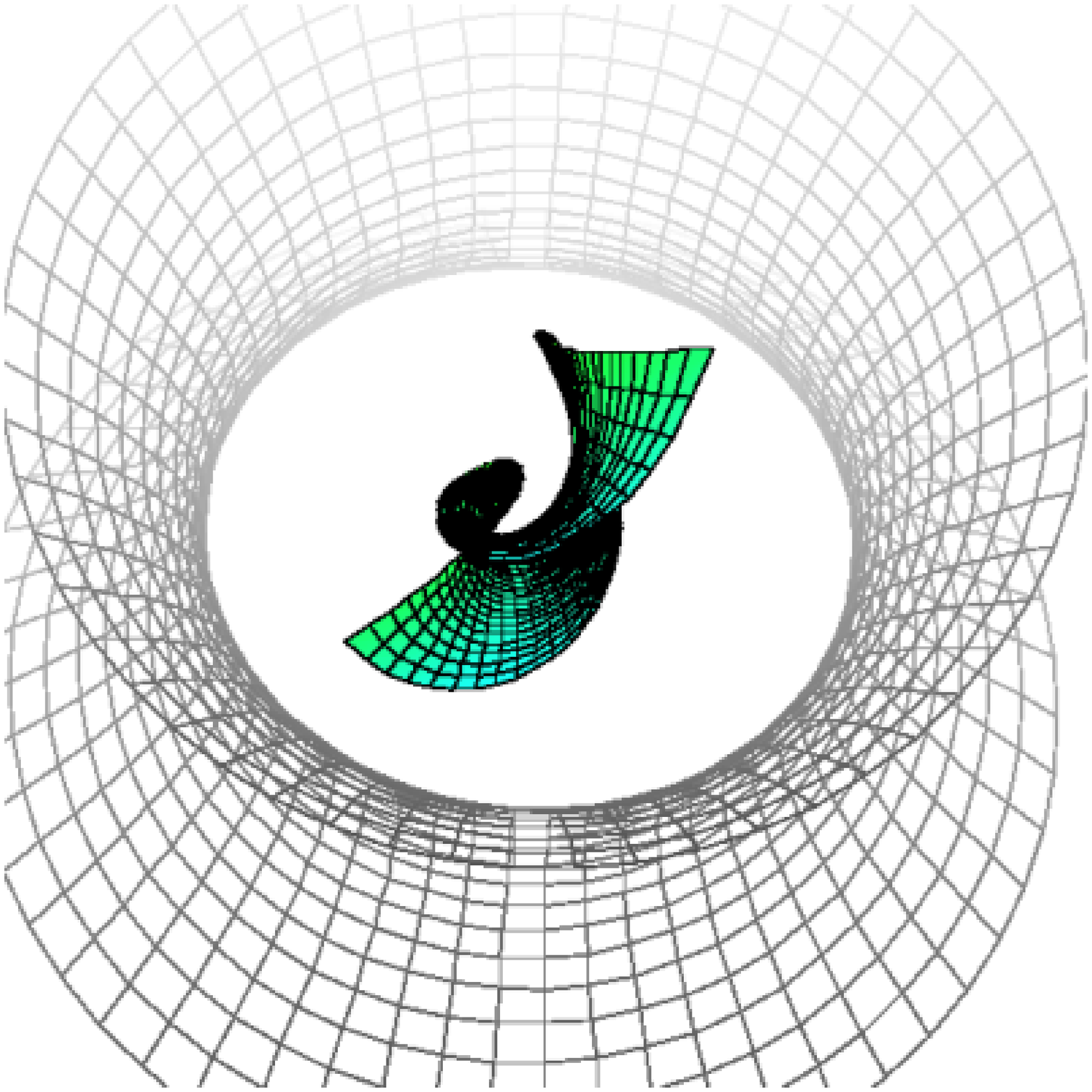,width=.48\textwidth}}}
\caption{Anti-isothermic type timelike Enneper cousin projected into ${\rm
  Int}\S^2_1(1)$ via
  $\wp_+$ with light cone and the boundary $\S^2_1(1)$ in
  $\E^3_1$\label{fig:tlenneper2}}
\end{figure}

Figures \ref{fig:tlminenneper2} and \ref{fig:tlminenneper3} show
different views of anti-isothermic type
timelike Enneper surface in $\E^3_1$.

\begin{figure}[ht]
\centering
\mbox{\subfigure[]{\epsfig{figure=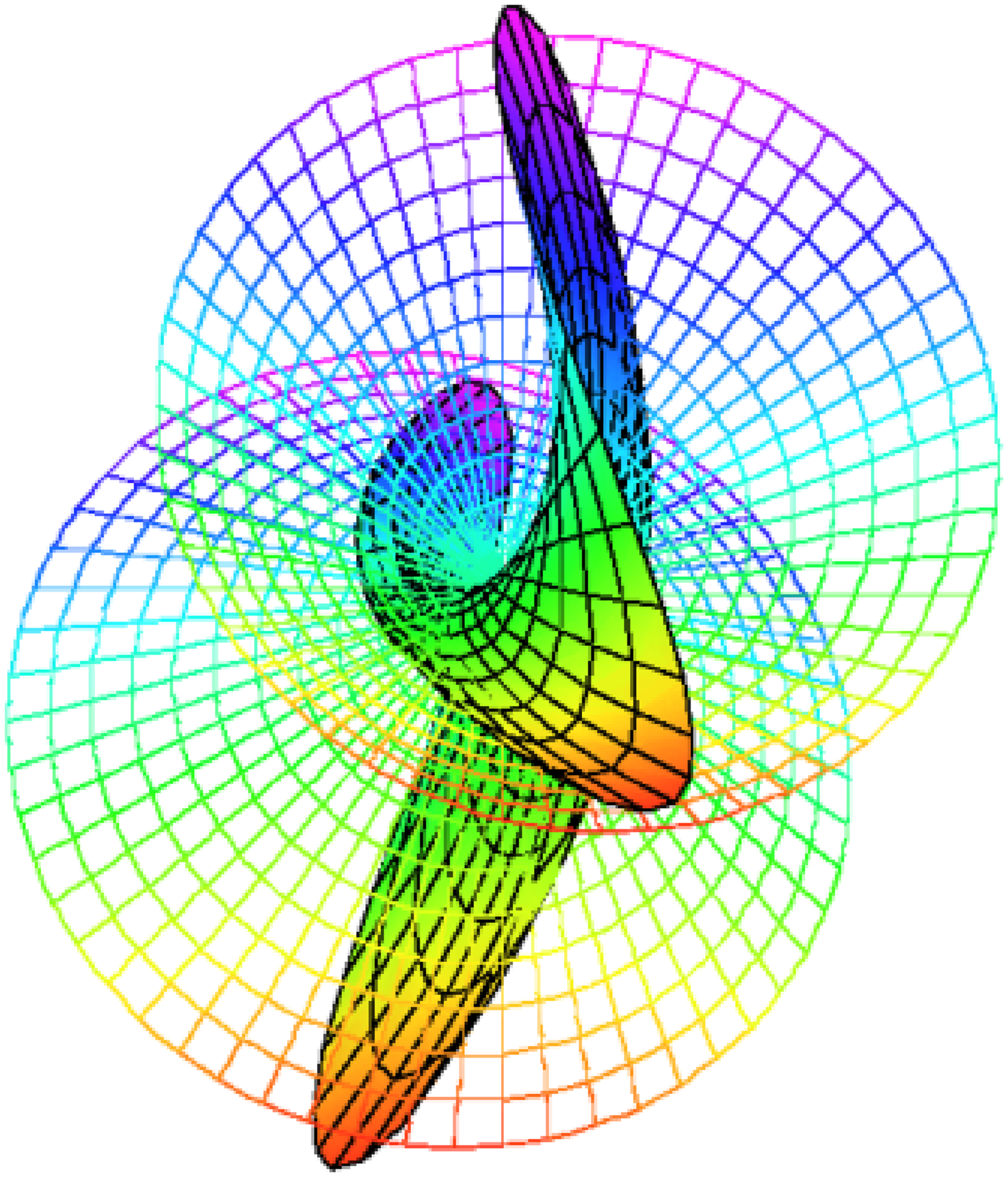,width=.48\textwidth}}\quad
    \subfigure[]{\epsfig{figure=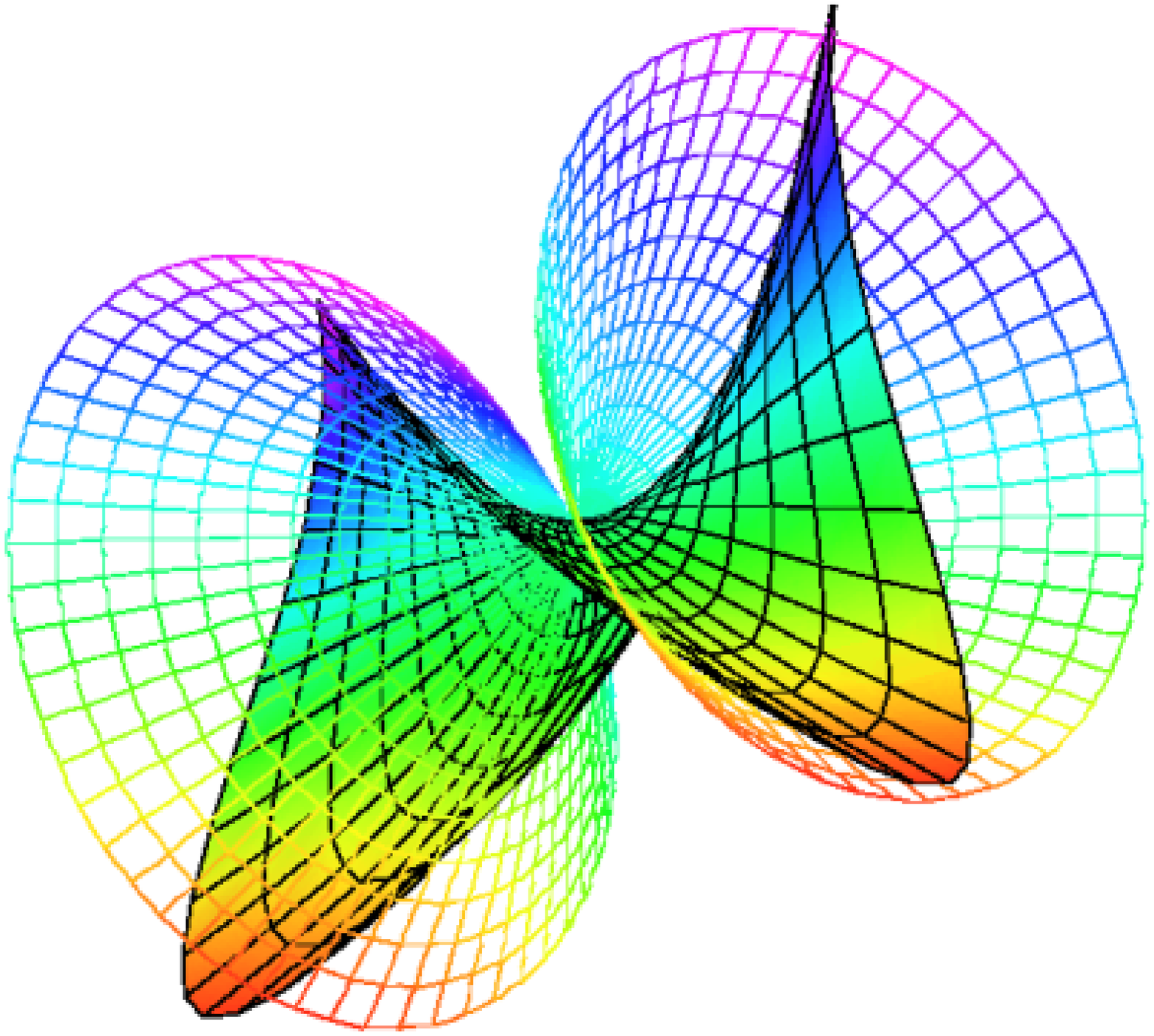,width=.48\textwidth}}}
\mbox{\subfigure[]{\epsfig{figure=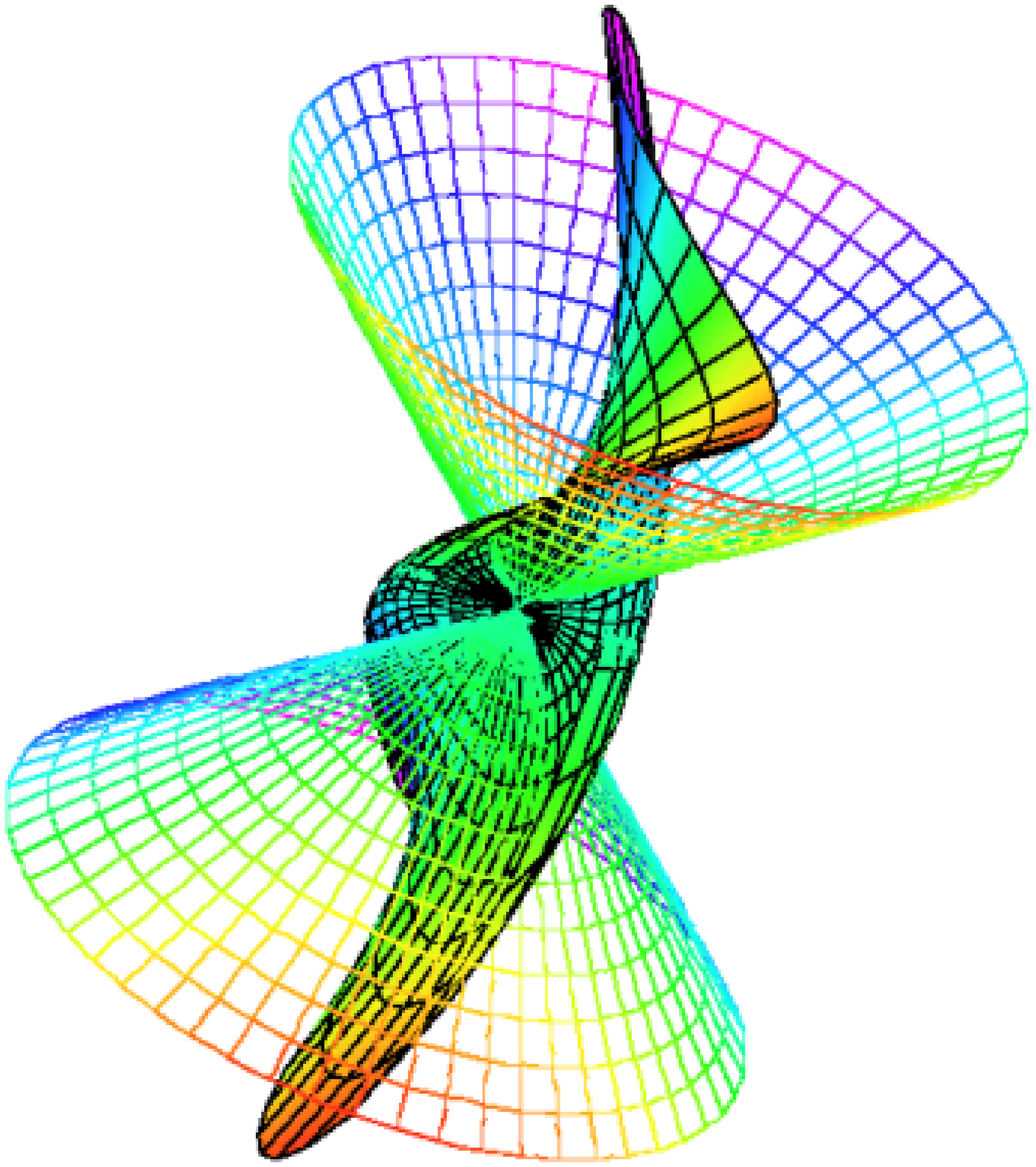,width=.48\textwidth}}}
\caption{Anti-isothermic type timelike Enneper surface in
  $\E^3_1$ with light cone\label{fig:tlminenneper2}}
\end{figure}

\begin{figure}[ht]
\centering
\mbox{\subfigure[]{\epsfig{figure=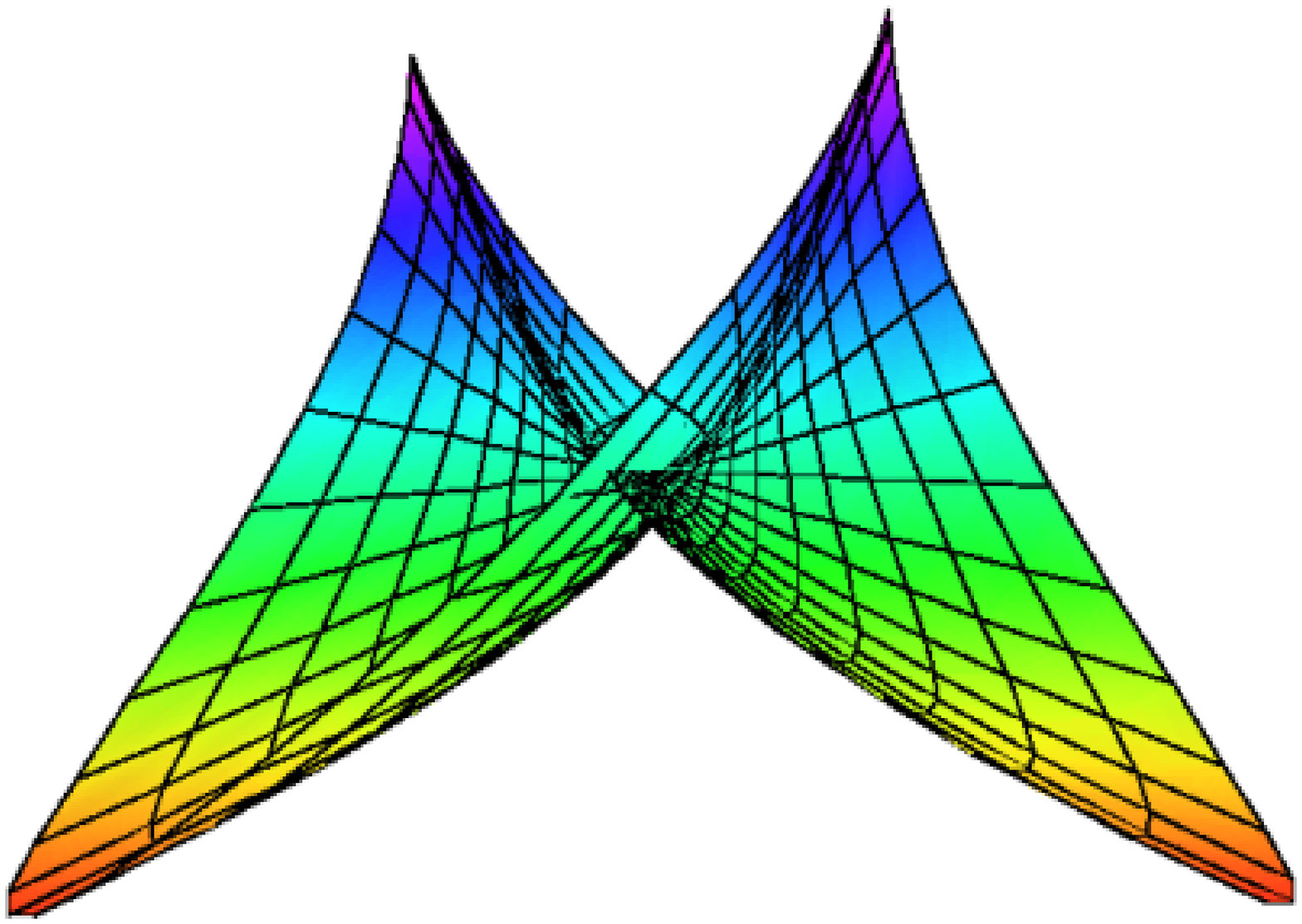,width=.48\textwidth}}\quad
    \subfigure[]{\epsfig{figure=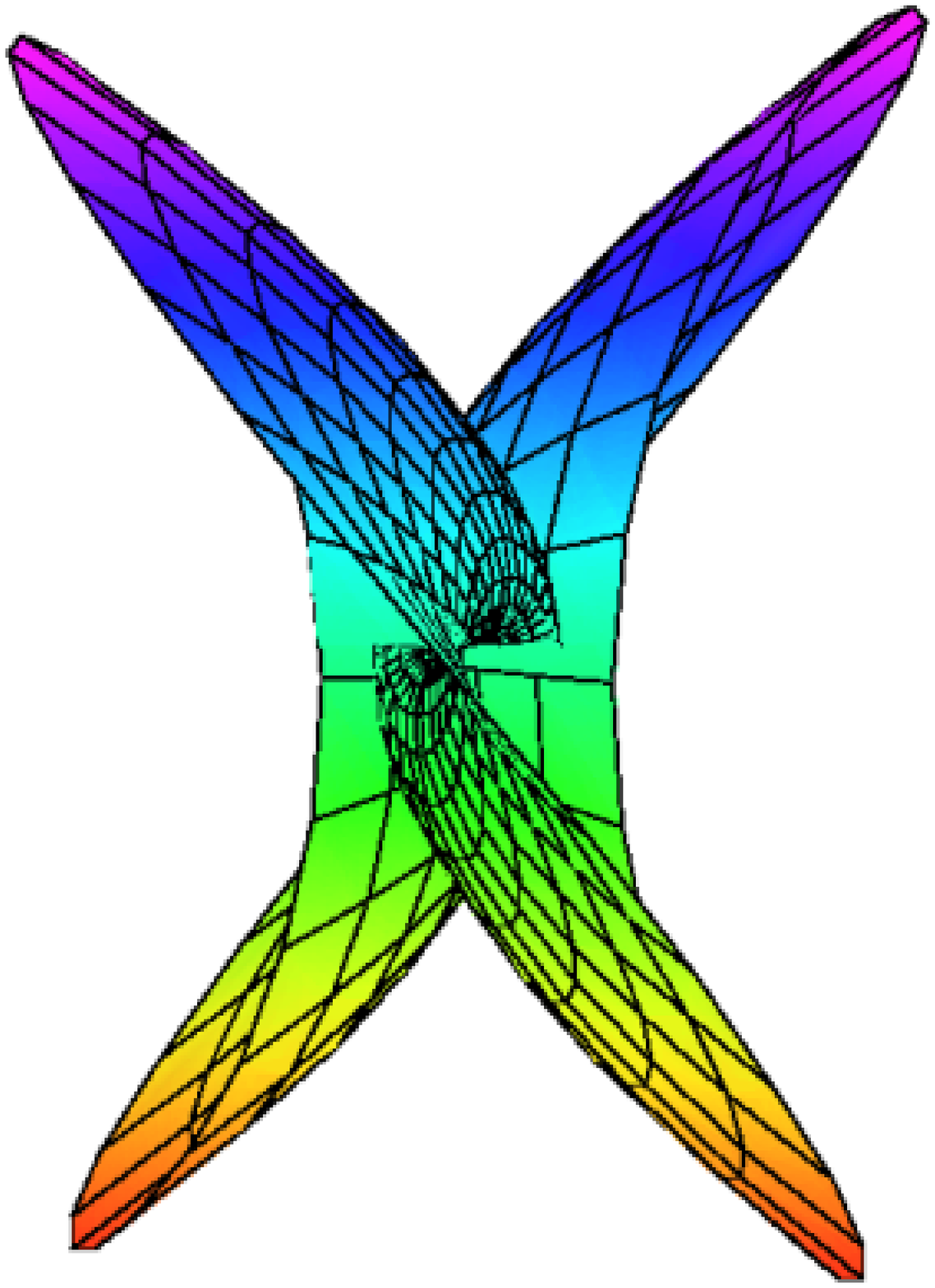,width=.48\textwidth}}}
\mbox{\subfigure[]{\epsfig{figure=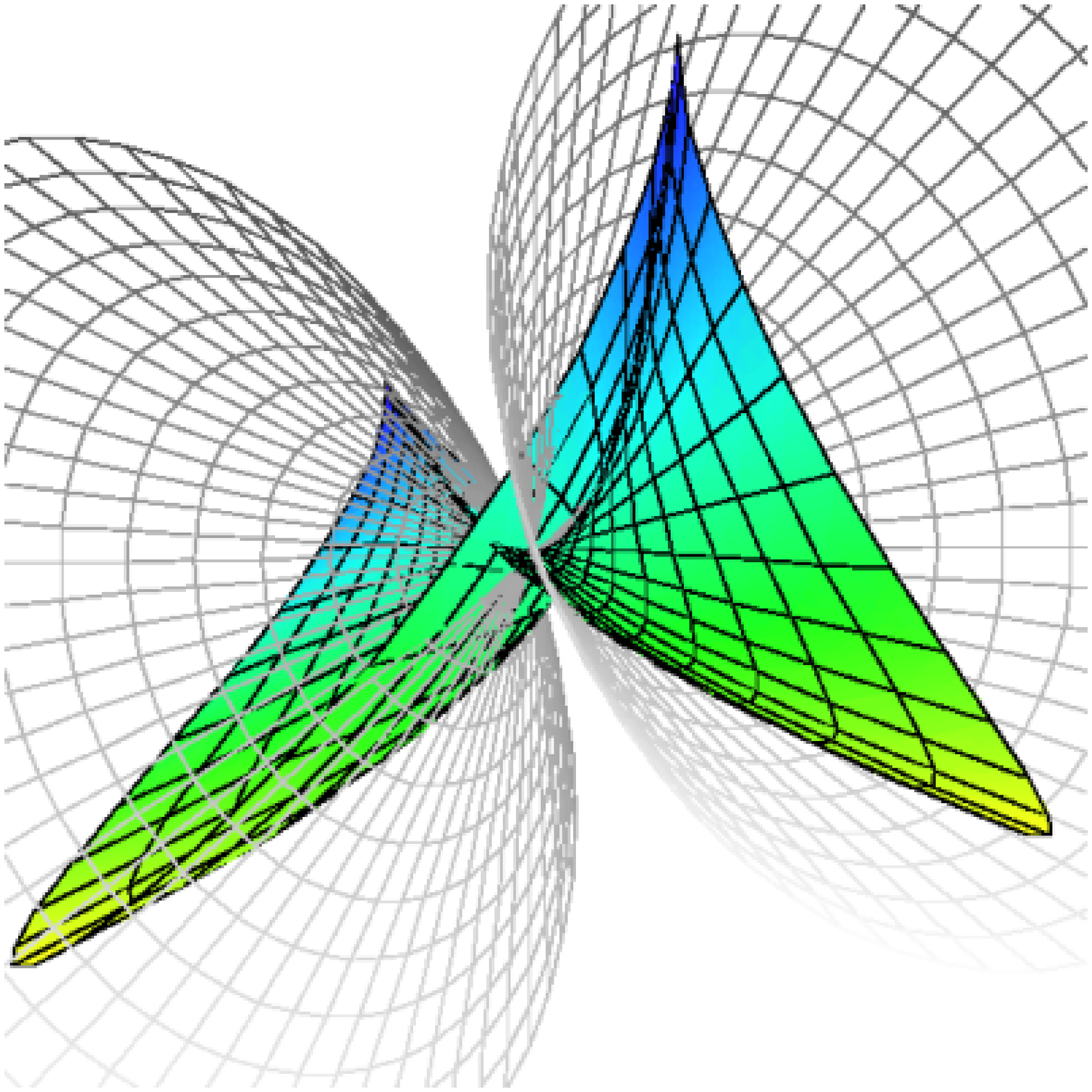,width=.48\textwidth}}}
\caption{Anti-isothermic type timelike Enneper surface in
  $\E^3_1$ with self-intersection\label{fig:tlminenneper3}}
\end{figure}
\end{example}
\begin{example}[$\cmc$ $1$ B-scroll in $\H^3_1(-1)$]

Let $(q,r)=(u,0)$. Then using the Bryant-Umehara-Yamada type
representation \eqref{eq:BUY}, we set up the following initial value
problem:
$$
F_1^{-1}dF_1=\begin{pmatrix}
u & -u^2\\
1 & -u
\end{pmatrix}du,\ F_2^{-1}dF_2=\begin{pmatrix}
0 & 0\\
1 & 0
\end{pmatrix}dv$$
with the initial condition $F_1(0)=F_2(0)=\1$.
This initial value problem has a unique solution
\begin{align*}
F_1(u,v)&=\begin{pmatrix}
\cosh u & \sinh u-u\cosh u\\
\sinh u & \cosh u-u\sinh u
\end{pmatrix},\\
F_2(u,v)&=\begin{pmatrix}
1 & 0\\
v & 1
\end{pmatrix}
\end{align*}
which are Lorentz holomorphic and Lorentz antiholomorphic null
immersions into ${\rm SL}_2\R$. The Bryant type representation formula
\eqref{eq:bryant} yields a timelike $\cmc$ $1$ surface
$$\varphi=F_1F_2^t=\begin{pmatrix}
\cosh u & -(u-v)\cosh u-\sinh u\\
\sinh u & -(u-v)\sinh u+\cosh u
\end{pmatrix}$$
in $\H^3_1(-1)$. The resulting surface is a correspondent of minimal
B-scroll in $\E^3_1$ under the Lawson-Guichard
correspondence.

Figure \ref{fig:b-scroll} shows different views of $\cmc$ $1$ B-scroll
projected in $\H^3_1(-1)$ via $\wp_+$ into the interior
of the boundary $\S^2_1(1)$.

\begin{figure}[ht]
\centering
\mbox{\subfigure[]{\epsfig{figure=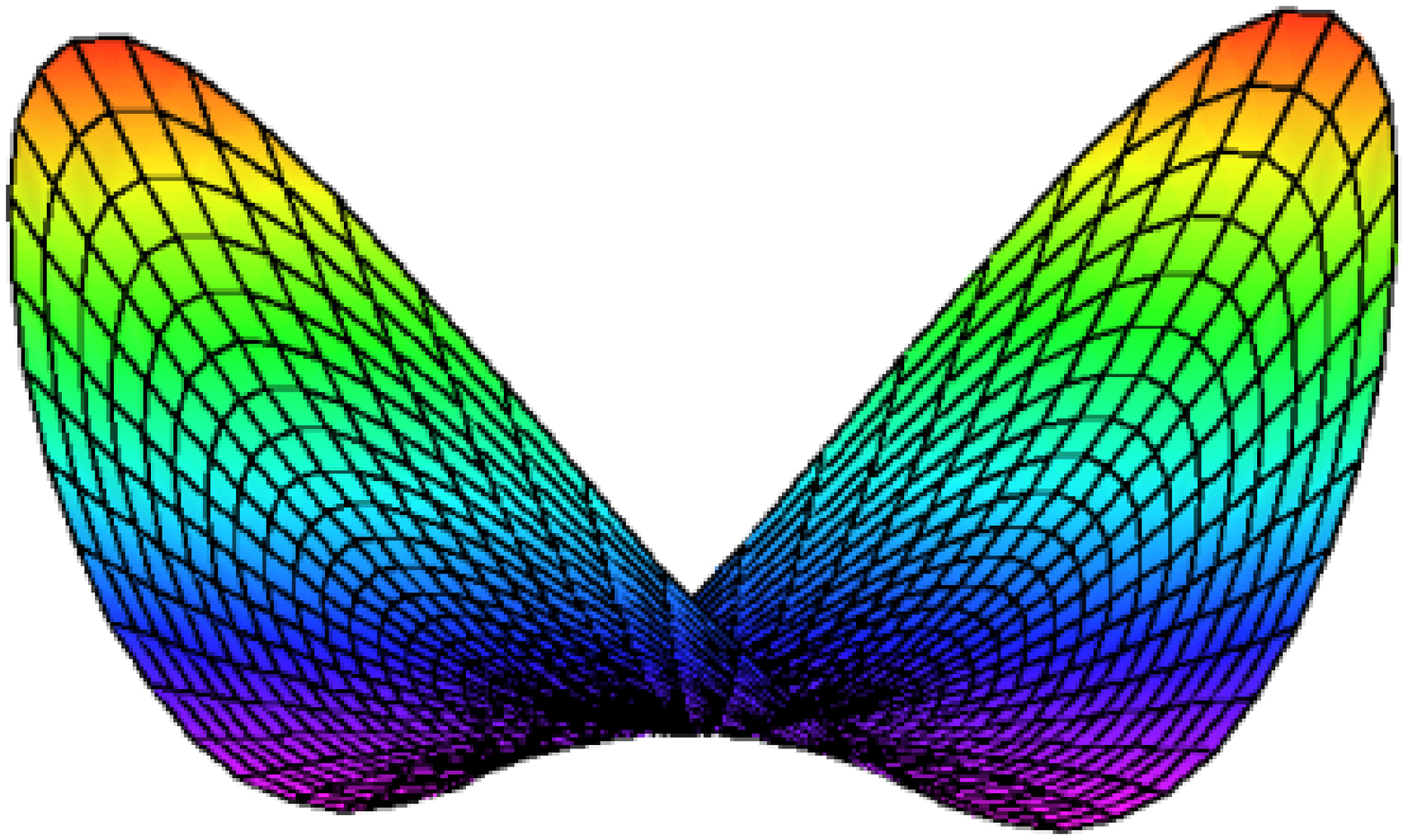,width=.48\textwidth}}\quad
    \subfigure[]{\epsfig{figure=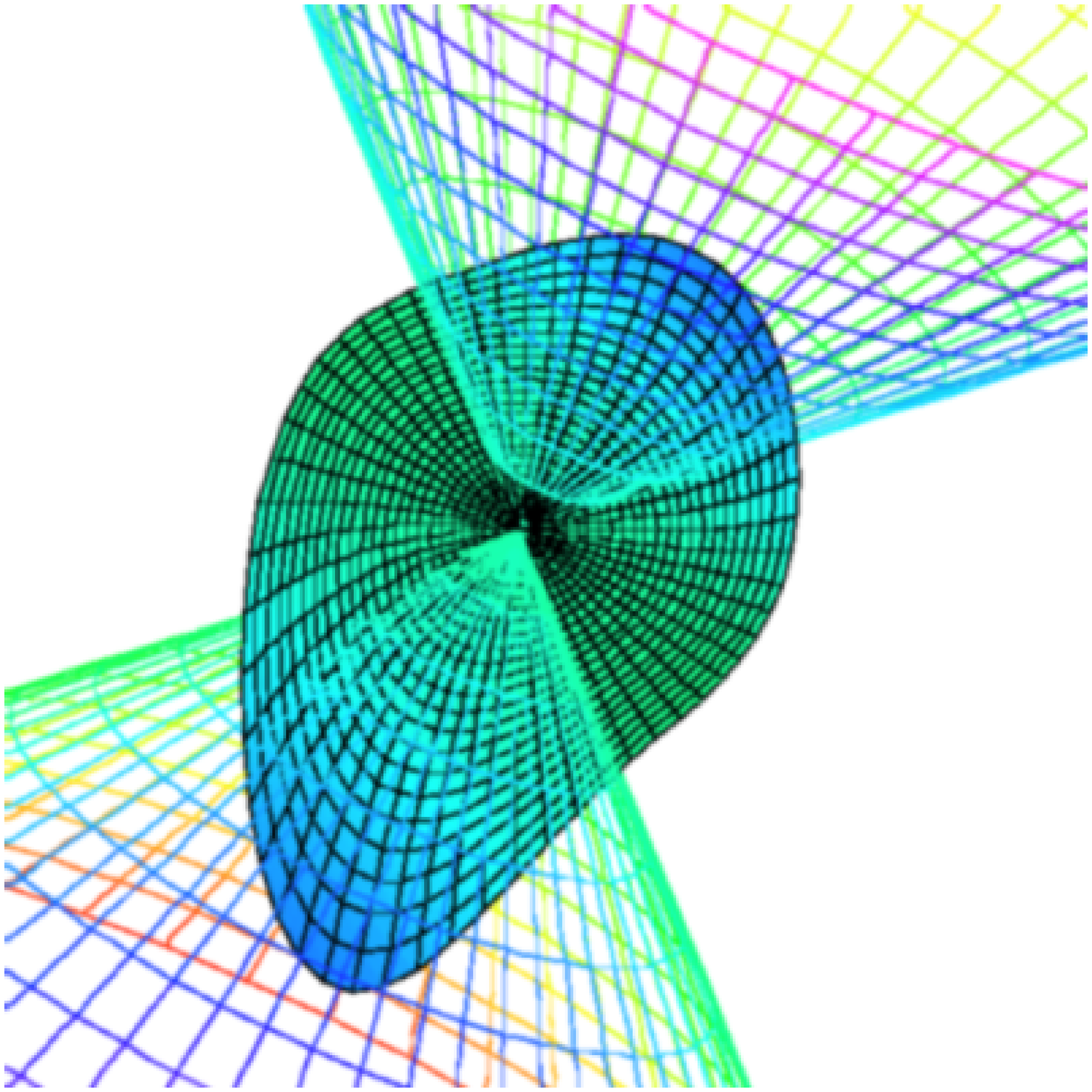,width=.48\textwidth}}}
\mbox{\subfigure[]{\epsfig{figure=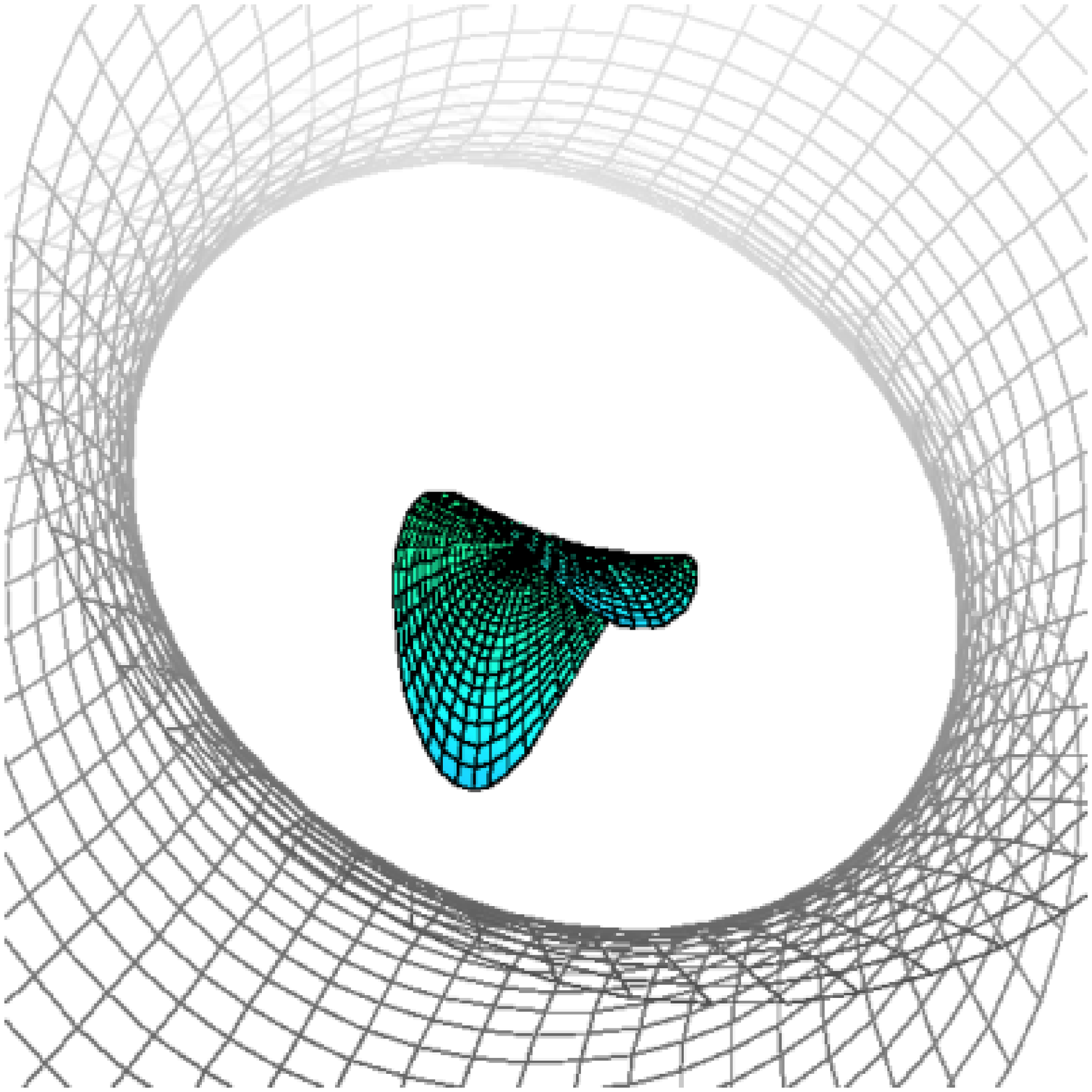,width=.48\textwidth}}}
\caption{$\cmc$ $1$ B-scroll projected into ${\rm
  Int}\S^2_1(1)$ via
  $\wp_+$ with light cone and the boundary $\S^2_1(1)$ in
  $\E^3_1$\label{fig:b-scroll}}
\end{figure}

Figure \ref{fig:minb-scroll} shows
different views of minimal B-scroll in $\E^3_1$.

\begin{figure}[ht]
\centering
\mbox{\subfigure[]{\epsfig{figure=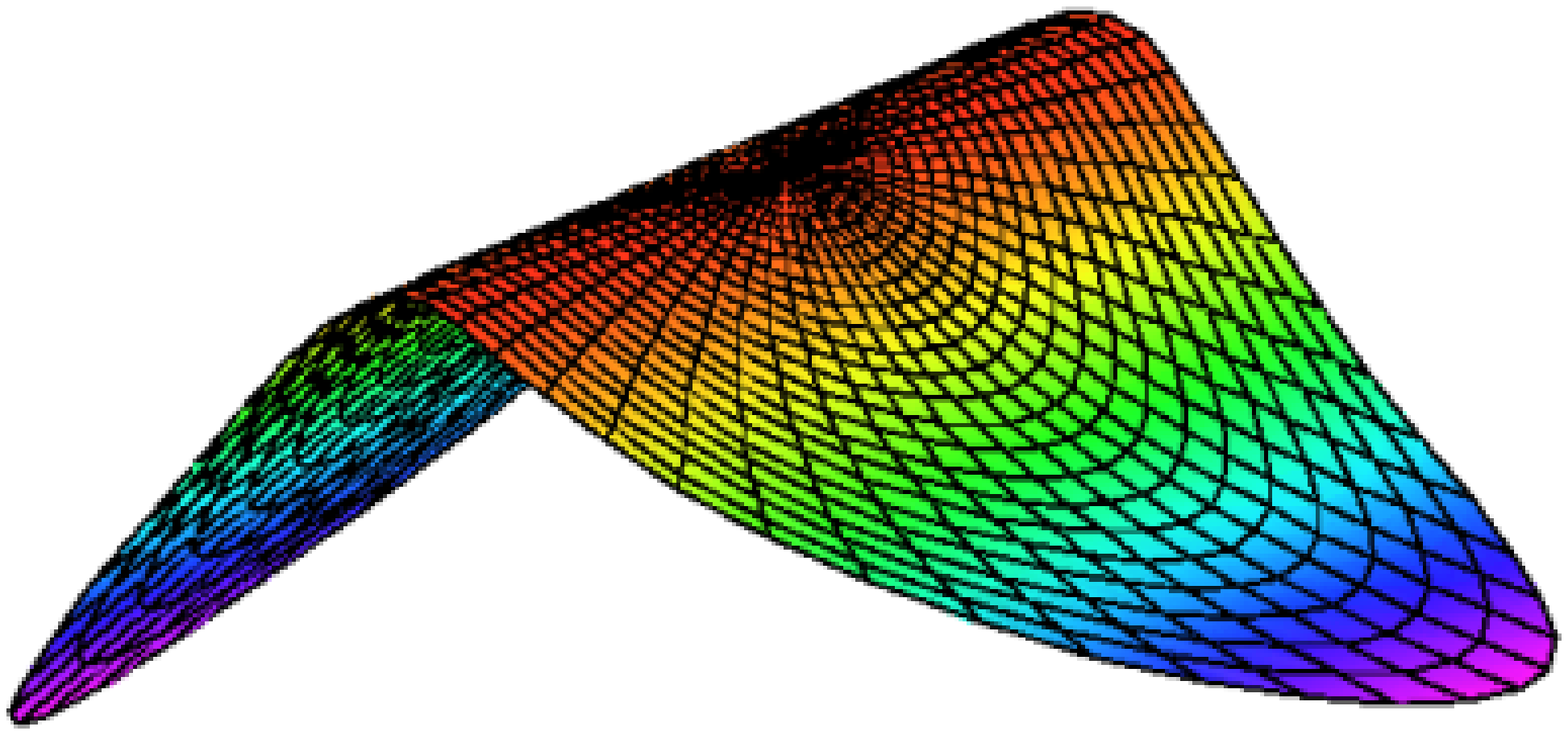,width=.48\textwidth}}\quad
    \subfigure[]{\epsfig{figure=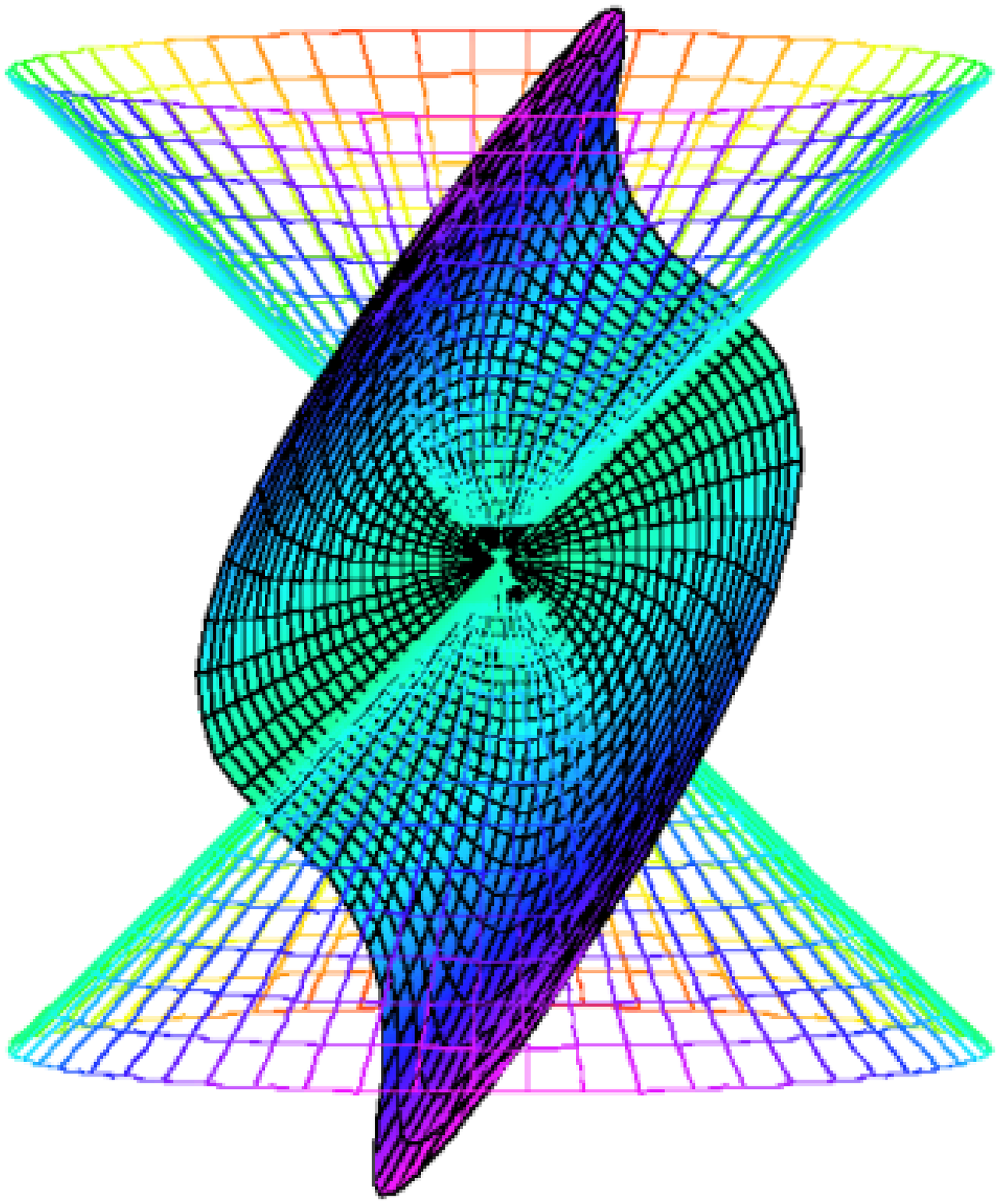,width=.48\textwidth}}}
\mbox{\subfigure[]{\epsfig{figure=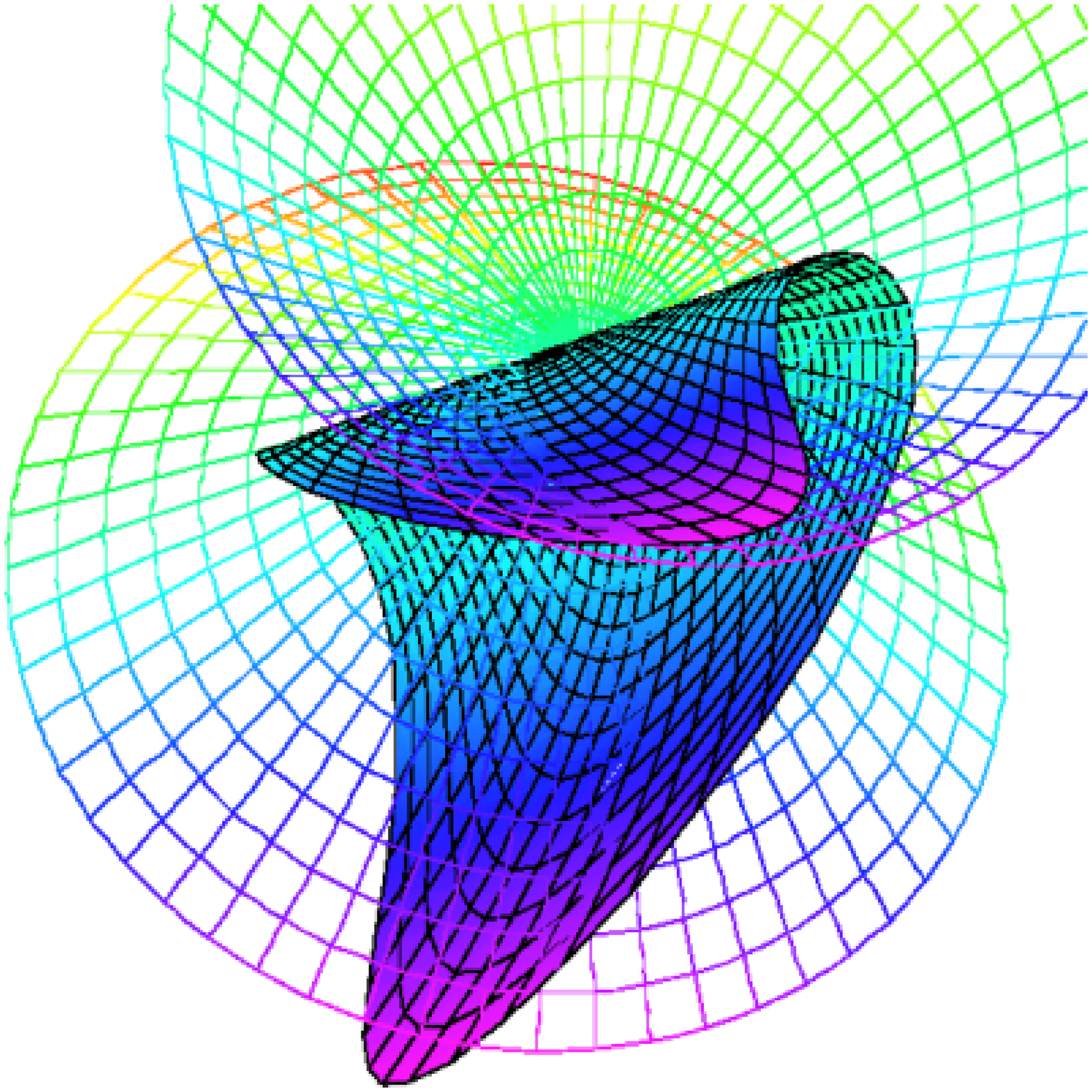,width=.48\textwidth}}}
\caption{Minimal B-scroll in
  $\E^3_1$ with light cone\label{fig:minb-scroll}}
\end{figure}
\end{example}

{\bf Acknowledgment:} The author wishes to thank Dr.~Jun-ichi
Inoguchi for drawing his attention to this subject and for many
invaluable suggestions. The author is grateful to Dr.~Magdalena Toda
for her hospitality during his short visit to Texas Tech University in
2003. Most of this paper grew out of discussions with
Dr.~Jun-ichi Inoguchi who was also visiting Texas Tech University during that time.

\end{document}